\documentclass[a4paper,12pt]{article}
\def\ti{\tilde}

\def\ov{\over}

\def\ni{\noindent}
\def\bi{\bigskip}
\def\bn{\bi\ni}

\def\a{\alpha}
\def\b{\beta}
\def\G{\Gamma}
\def\g{\gamma}
\def\r{\rho}
\def\d{\delta}
\def\e{\epsilon}
\def\t{\tau}
\def\w{\omega}
\def\x{\w\ov 2\pi}
\def\s{\sigma}
\def\l{\lambda}
\def\u{\mu}
\def\v{\nu}
\def\th{\theta}

\def\G{\Gamma}
\def\S{\Sigma}
\def\z{\zeta}
\def\x{\chi}

\def\T{\mathop{\bf T}}
\def\N{\mathop{\bf N}}
\def\L{\Lambda}
\def\R{\mathop{\bf R}}

\def\C{\mathop{\bf C}}
\def\Z{\mathop{\bf Z}}

\def\D{\Delta}
\def\Im{\mathop{\rm Im}}

\def\Id{\mathop{\rm Id}}
\def\pa{\partial}

\def\Im{\mathop{\rm Im}}
\def\ad{\mathop{\rm ad}}
\def\Ad{\mathop{\rm Ad}}
\def\Log{\mathop{\rm Log}}
\def\Conj{\mathop{\rm Conj}}

\def\lab{\label}

\def\Id{\rm Id}
\def\L{\rm L}
\def\<{\langle}
\def\>{\rangle}
\def\bm{\begin{pmatrix}}
\def\em{\end{pmatrix}}
\def\rot{{\rm rot}}
\def\Ree{{\rm {Re}}}
\usepackage{latexsym}
\usepackage{amsmath}
\usepackage{theorem}
\usepackage[english]{babel}

\renewenvironment{proof}{\ni{\ni\it Proof}\quad\ni}{\vskip .2cm\hfill$\Box$\vskip .2cm}

\newtheorem{theo}{Theorem}[section]
\newtheorem{lemme}{Lemma}[section]
\newtheorem{prop}{Proposition}[section]
\newtheorem{cor}{Corollary}[section]

\begin{document}
\author{Rapha\"el KRIKORIAN\\[.5cm] Centre de Math\'ematiques-UMR 7640 du CNRS \\ Ecole Polytechnique\\ 91128 PALAISEAU\\
and UMA/ENSTA, 32 bd Victor, 75015 PARIS}
\date{}
\title{Reducibility, differentiable rigidity and Lyapunov exponents \\ for quasi-periodic cocycles on ${\T}\times SL(2,{\R})$,}
 \maketitle
\abstract{\it Given $\a$ in some set $\Sigma$ of total (Haar) measure in ${\T}={\R}/{\Z}$, and $A\in C^{\infty}({\T},SL(2,{\R}))$ which is homotopic to the identity, we prove that if the fibered rotation number of the skew-product system $(\a,A):{\T}\times SL(2,{\R})\to {\T}\times SL(2,{\R})$, $(\a,A)(\th,y)=(\th+\a,A(\th)y)$ is diophantine with respect to $\a$ and if the fibered products are uniformly bounded in the $C^0$-topology then the cocycle $(\a,A)$ is $C^\infty$-reducible --that is $A(\cdot)=B(\cdot+\a)A_0B(\cdot)^{-1}$, for some  $A_0\in SL(2,{\R})$, $B\in C^{\infty}({\T},SL(2,{\R}))$. This result which can be seen as a non-pertubative version of a theorem by L.H. Eliasson has two interesting corollaries: the first one is a result of differentiable rigidity: if $\a\in\Sigma$ and the cocycle $(\a,A)$ is $C^0$-conjugated to a constant cocycle $(\a,A_0)$ with $A_0$ in a set of total measure in $SL(2,{\R})$ then the conjugacy is $C^\infty$; the second consequence is: if $\a\in \Sigma$ is fixed then the set of $A\in C^\infty({\T},SL(2,{\R}))$  for which $(\a,A)$ has positive Lyapunov exponent is $C^\infty$-dense. A similar result is true for the Schr\"odinger cocycle and for 2-frequencies conservative differential equations in the plane. } 

\section{Introduction}
In the following paper we shall be interested in  smooth quasi-periodic cocycles on ${\T}\times SL(2,{\R})$ where ${\T}$ denotes the circle ${\R}/{\Z}$ (the base) and $SL(2,{\R})$ the set of 2 by 2 real matrices with determinant 1 (the fibre). 
Such cocycles are by definition  diffeomorphisms on the product  ${\T}\times SL(2,{\R})$ of the form 
\begin{align*}(\a,A): {\T}\times SL(2,{\R})&\to {\T}\times SL(2,{\R})\\
(\th,y)&\mapsto(\th+\a,A(\th)y)\end{align*}
where $\a\in{\T}$ and $A\in C^\infty({\T},SL(2,{\R}))$.
Our aim is to understand the dynamics of $(\a,A)$ that is the iterated diffeomorphisms $(\a,A)^n$, $n\in{\Z}$.
Notice that this amounts to understand the quasi-periodic products ($n\geq 1$)
$$\begin{cases}&A_n(\cdot)=A(\cdot+(n-1)\a)\cdots A(\cdot)\\
&A_{-n}(\cdot)=A(\cdot-n\a)^{-1}\cdots A(\cdot-\a)^{-1}\end{cases}$$
since for $n\in{\Z}$ one has $(\a,A(\cdot))^n=(n\a,A_n(\cdot))$

The general definition of a quasi-periodic cocycle can be done along the same lines by replacing  the group $SL(2,{\R})$  by any other group, ${\T}$  by ${\T}^d$ and by taking $\a$  in ${\T}^d$. 
\subsection{Understanding the dynamics}
There are at least two ways of understanding the dynamics of ``generic'' $(\a,A)$ (we use the word ``generic'' in a very vague sense) . 
\subsubsection{Reducibility}
The first one is to (try to) {\it conjugate}  $(\a,A)$ to a simpler system e.g a constant one. By this we mean that we try to find a cocycle $(0,B)$ with $B\in C^\infty({\T},SL(2,{\R}))$ and a constant $A_0\in SL(2,{\R})$ such that
$$(\a,A)=(0,B)\circ(\a,A_0)\circ (0,B)^{-1},$$
which is equivalent to the following equality for all $\th\in{\T}$:
$$A(\th)=B(\th+\a)A_0B(\th)^{-1}.$$
In that case the dynamics of $(\a,A)$ is perfectly understood since for every $n\in{\Z}$
$$(\a,A)^n=(0,B)\circ(\a,A_0)^n\circ(0,B)^{-1}.$$
In case  the base is the circle and the fibre is the group $SU(2)$ of 2 by 2 unitary matrices with determinant 1, this attempt to understand the dynamics of $(\a,A)$  proved to be succesfull. Unfortunately the situation is more complicated in the $SL(2,{\R})$-case due to the fact that fibered products are not always bounded (Lyapunov exponents). Nonetheless when the fibre is the group $SL(2,{\R})$ there is another way to understand the dynamics of the cocycle of $(\a,A)$ in a looser sense (cf. section~\ref{sec:1.1.2} Lyapunov exponents). 

\subsubsection{Dynamical invariants}\lab{sec:1.1.2}
\paragraph{The rotation number}
Let us  denote by ${\bf S}^1$ the set of vectors of ${\R}^2$ of euclidian norm 1 and by $\pi:{\R}\to {\bf S}^1$ the projection $\pi(x)=e^{2\pi ix}$ (we have identified ${\R}^2$ with ${\C}$). Assume that $A(\cdot):{\T}\to SL(2,{\R})$ is continuous and {\it homotopic to the identity}; then the same is true for the map 
\begin{align*}F:{\T}\times{\bf S}^1&\to {\T}\times{\bf S}^1\\
(\th,v)&\mapsto (\th+\a,{A(\th)v\ov\|A(\th)v\|})\end{align*}
and  therefore it admits  a continuous lift $\ti F:{\T}\times{\R}\to {\T}\times{\R}$ of the form $\ti F(\th,x)=(\th+\a,x+f(\th,x))$ such that $f(\th,x+1)=f(\th,x)$ and $\pi (x+f(\th,x))=A(\th)\pi(x)/\|A(\th)\pi(x)\|$. In order to simplify the terminology we shall say that $\ti F$ is a lift for $(\a,A)$. The map $f$ is independent of the choice of  the lift up to the addition of a constant rational integer $p\in{\Z}$. Now, since $\th\mapsto\th+\a$ is uniquely ergodic on ${\T}$, we can invoque  a theorem by  M.R. Herman and Johnson-Moser (cf.~\cite{He}, \cite{J-M}): for {\it every} $(\th,x)\in {\T}\times {\R}$ the  limit 
$$\lim_{n\to\pm\infty}{1\ov n}\sum_{k=0}^{n-1}f({\ti F}^k(\th,x)),$$
exists, is independent of $(\th,x)$ and the convergence is uniform in $(\th,x)$; the class of this number in ${\R}/{\Z}$ (which is independent of the chosen lift) is called the {\it fibered} rotation number of $(\a,A)$. Moreover $\r_f(\a,A)$ is continuous with respect to $A$. In fact for any $\ti F$-invariant measure $\u$ on ${\T}\times {\bf S}^1$ 
$$\r=\int_{{\T}\times{\bf S}^1}f(\th,x)d\u(\th,x).$$

As we shall see this rotation number is also almost invariant by $C^0$-conjugacy. We recall that the first homotopy group of $SL(2,{\R})$ is isomorphic to ${\Z}$ with generator $E_1(\cdot)$
$$\forall \th\in{\T}\ E_1(\th)=\begin{pmatrix}\cos(2\pi\th)&-\sin(2\pi\th)\\ \sin(2\pi\th)&\cos(2\pi\th)\end{pmatrix}.$$ 
Then a continuous map $B:{\T}\to SL(2,{\R})$ is said to be of degree $r\in{\Z}$ if it is homotopic to $E_r(\cdot)$=$E_1(\cdot)^r$. It is equivalent to say that $B(\cdot)E_r(\cdot)^{-1}$ is homotopic to the identity.

\begin{prop}If $A:{\T}\to SL(2,{\R})$ is continuous and homotopic to the identity and if $B:{\T}\to SL(2,{\R})$ is continuous of degree $r$, then
$$\r_f((0,B)\circ(\a,A)\circ (0,B)^{-1})=\r_f((\a,A))+r\a\mod 1.$$
\end{prop}
\begin{proof}

We first consider the case when $B$ is homotopic to the identity. With the previous notations there exists a map $\ti H:{\T}\times {\R}\to {\R}$ which is a lift for $(0,B)$. Consequentely $\ti H\circ \ti F\circ \ti H^{-1}$ is a lift for $(0,B)\circ (\a,A)\circ (0,B)^{-1}$ and it is clear that
$$\r_f((0,B)\circ(\a,A)\circ (0,B)^{-1})=\r_f((\a,A))\mod 1.$$
Now if $B$ is of degree $r\in{\Z}$ then it can be written $B(\cdot)=E_r(\cdot).\bar B(\cdot)$ with $\bar B(\cdot)$ homotopic to the identity. It is then enough to
check the proposition for $E_r(\cdot)$ in place of $B(\cdot)$. But in that case it is not difficult to check that a lift for $(\a,E_r(\cdot+\a)A(\cdot)E_r(\th)^{-1})$ is given by
$$\ti G(\th,x)=(\th+\a,x+\ti f(\th,x-r\th)+r\a).$$
If we define $\ti g(\th,x)=\ti f(\th,x-r\th)+r\a$ and if we call $\ti g_k$ and $\ti f_k$ the $k$-th Birkhoff sums of $\ti g$ and $\ti f$ respectively along $\ti G$ and $\ti F$, a simple computation shows that
$$\ti g_k(\th,x)=kr\a+\ti f_k(\th,x-r\th),$$
which gives the conclusion of the proposition since the convergence of the associated Birkhoff means to the corresponding rotation numbers is uniform.

\end{proof}

\bn{\bf Remark 1:} when 
\begin{itemize}
\item $B(\cdot)$ is only defined on ${\R}/{2 {\Z}}$  and the degree of $\th\mapsto B(2\th)$ is $r$ and at the same time
\item $B(\cdot+\a)A_0B(\cdot)^{-1}$ is defined on ${\R}/{\Z}$
\end{itemize}
a similar result is true:
$$\r_f((0,B)\circ(\a,A)\circ (0,B)^{-1})=\r_f((\a,A))+{1\ov 2}\a\mod {1}.$$
(Use the fact that $\r_f(\a,A(\cdot))= \r_f(\a/2,A(2\cdot))$).

\bn{\bf Remark 2:} when $A_0$ is a constant elliptic matrix conjugated to 
$$\begin{pmatrix}\cos 2\pi\psi&-\sin2\pi\psi\\
\sin2\pi\psi &\cos2\pi\psi\end{pmatrix}$$
($\psi\in{\R}$) the fibered rotation number of $(\a,A_0)$ is $\psi\mod 1$

\paragraph{The Lyapunov exponent}
It is defined as the limit
$$\l(\a,A):=\lim_{n\to\infty}{1\ov n}\int_{\th\in{\T}}\Log \|A_n(\th)\|,$$
which by Kingman's subbadditive theorem always exists (similarly the limit when $n$ goes to $-\infty$ exists and is equal to $\l(\a,A)$). Moreover,  for almost every $\th \in{\T}$ with respect to the Haar measure, the following limit exists and is equal to $\l(\a,A)$:
$$ \lim_{n\to\infty}{1\ov n}\Log \|A_n(\th)\|$$
exists and is equal to $\l(\a,A)$. Now, when the Lyapunov exponent is positive (non zero)  Oseledec's theorem tells us that for Lebesgue-almost $\th\in{\T}$ there exists a  splitting
$${\R}^2=E^s(\th)\oplus E^u(\th),$$
which is invariant by $(\a,A)$ which means that: 
\begin{itemize}
\item i) $E^s(\th)$, $E^u(\th)$ (viewed as elements of $P^1{\R}^2$) are measurable with respect to $\th$ and 
\item ii) for almost every $\th$ 
$$E^{s,u}(\th+\a)=A(\th)E^{s,u}(\th);$$
\item iii) for almost every $\th$ and every $v\in E^{s}(\th)$ (resp. $v\in E^{u}(\th)$) we have
$$\lim_{n\to\infty}{1\ov n}\Log\|A_n(\th)v\|=-\l\ \ ({\rm resp.}\ \lim_{n\to\infty}{1\ov n}\Log\|A_n(\th)v\|=\l)$$
\end{itemize}
This is the meaning we give to ``understand the dynamics'' in that case though it is abusive.

As we shall see later this two approaches (via reducibility and Lyapunov exponents) are in fact complementary.

\subsection{The local situation}
The link between reducibility and the two dynamical invariants we have introduced is  understood at least in  the local situation, that is when $A(\cdot)\in C^\infty({\T},SL(2,{\R})$ is $C^\infty$-close to some constant matrix $A_0$. 
%Historically, this appeared when people try to understand the spectrum of the one-dimensional Schr\"odinger operator with a quasi-periodic potential. The rotation number and the Lyapunov exponent appeared to have some spectral interpretation (in terms of the integrated density of states). 
One has to make some extra-assumption namely that the frequency vector on the base satisfies some diophantine property. We shall say that $\a\in{\T}^d$ satisfies a condition $CD(\g,\s)$ ($\g,\s>0$) if there exist constants $K$ and $\t>0$ such that
$$\forall k\in{\Z}^d-\{0\},\ \ \min_{l\in{\Z}}|k\a-l|\geq{\g^{-1}\ov |k|^\s}.$$

Moreover if $\r$ is  in ${\T}$ we say that $\r$ is {\it diophantine with respect to $\a$} if
$$\min_{l\in{\Z}}|\r-{1\ov 2}\<k,\a\>-{1\ov 2}l|\geq \frac{K^{-1}}{|k|^\t}.$$
If  there exists $k_0$ such that
$$\r={1\ov 2}\<k_0,\a\>\mod \frac{1}{2}$$
we say that $\r$ is {\it rational} with respect to $\a$. 
\begin{theo}Let $\a\in{\T}^d$ fixed and satisfying a diophantine condition $CD(\g,\s)$. Let also $A_0$ be in $SL(2,{\R})$. Then there exists $\e_0>0$ and $s_0\in{\N}$ such that for any $A(\cdot)\in C^\infty({\T}^d,SL(2,{\R}))$ satisfying
$$\|A(\cdot)-A_0\|_{C^{s_0}}\leq \e_0$$
the following is true:

i) If, $\r_f(\a,A)$ is either  diophantine or rational w.r.t to $\a$ then $(\a,A)$ is $C^\infty$-reducible in ${\R}/{2{\Z}}$ (i.e there exists a smooth $B:{\R}/2{\Z}\to SL(2,{\R})$ such that $B(\cdot+\a)^{-1}A(\cdot)B(\cdot)$ is constant). 

ii) Moreover, if $\r_f(\a,A)$ is diophantine w.r.t $\a$ then the conjugacy can be chosen so that it is defined on ${\R}/{\Z}$ (this changes the value of $A_0$). This is also  true if  $(\a,A)$ has bounded fibered products and $\r_f(\a,A)$ is rational w.r.t $\a$ and in that case $A_0=\pm Id$.
\end{theo}
In a continuous time set-up (differential equations) and for Schr\"odinger equations with analytic potentials  this theorem, the proof of which is based on a KAM scheme, is due to  L.H Eliasson (cf.~\cite{El}) but the proof in the general case is the same (local results in other groups are also true; see for example~\cite{K_ast}, \cite{K_annens}).  

To complete the picture in the local situation we give the following theorem:
\begin{prop}\label{loc.sit}Any (elliptic) constant $A_0$ in $SL(2,{\R})$ is $C^\infty$-accumulated by functions $A(\cdot)\in C^\infty({\T},SL(2,{\R}))$ such that $(\a,A(\cdot))$ is hyperbolic (in the fiber).
\end{prop}
\begin{proof}

We can assume that $A_0$ is a rotation matrix $R_\phi$. Now let $\e>0$ and $s\in{\N}$ be arbitrary and consider $k\in{\Z}$ such that
$$\min_{l\in{\Z}}|\phi-k\a|\leq {\e\ov 2}.$$
We then have
$$\|R_\phi-R_{k\a}\|\leq {\e\ov 2},$$
and also if we set $B(\th)=R_{-k\th}$,
$$B(\th+\a)R_{k\a}B(\th)^{-1}=Id.$$
Now choose $H$ a hyperbolic matrix in $SL(2,{\R})$ such that
$$\|H-Id\|\leq {\e\ov 20|k|^s},$$
and define
$$A(\th)=B(\th+\a)^{-1}HB(\th).$$
A simple calculation shows that
$$\|A(\cdot)-R_{k\a}\|_{C^s}\leq 10|k|^s.{\e\ov 20 |k|^s},$$
and therefore
$$\|A(\cdot)-A_0\|_{C^s}\leq \e,$$
while $(\a,A(\cdot))$ is $C^\infty$-conjugated to a hyperbolic system.
\end{proof}
\section{The main theorem and its corollaries}
The main theorem (A$+$B) of this paper is to provide an extension of Eliasson's theorem to the global case. Since we are dealing with a non-perturbative situation ($A(\cdot)$ is no longer close to a constant) it is not surprizing that we have to make an extra assumption.
\begin{theo}[Main theorem A]\lab{main:theo}There exists a set $\S$ of total Haar measure in ${\T}$ (for the definition of $\S$ we refer to section~\ref{sec:3}) such that for any fixed  $\a\in\S$ and any $A(\cdot)\in C^\infty({\T},SL(2,{\R}))$ satisfying
\begin{itemize}
\item $A(\cdot)$ is homotopic to the identity,
\item $\r_f(\a,A)$ is diophantine w.r.t $\a$ and 
\item the fibered products of $A(\cdot)$ are $C^0$-bounded i.e:
$$\sup_{k\in{\Z}}\max_{\th\in{\T}}|A_k(\cdot)|<\infty$$
\end{itemize} 
then $(\a,A)$ is $C^\infty$-reducible on ${\R}/{\Z}$.
\end{theo}
The case of cocycles non-homotopic to the identity is also interesting:
\begin{theo}[Main theorem B]\lab{main:theo} If $\a\in\S$ and  $A(\cdot)\in C^\infty({\T},SL(2,{\R}))$ is such that 
\begin{itemize}
\item $A(\cdot):{\T}\to SL(2,{\R})$ is of degree $r\in{\Z}-\{0\}$ and
\item the fibered products of $A(\cdot)$ are bounded
\end{itemize} 
then $(\a,A)$ is $C^\infty$-conjugated  on ${\R}/{\Z}$ to a cocycle $(\a,E_r(\cdot))$ where 
$$E_r(\th)=\begin{pmatrix}\cos(2\pi r\th)&-\sin(2\pi r\th)\\\sin(2\pi r\th)&\cos(2\pi r\th)\end{pmatrix}.$$
\end{theo}
The  boundedness assumption on the iterates of the cocycle  is of course very strong, but surprisingly enough the theorem still  has some very interesting consequences. The first one is a {\it differentiable rigidity} theorem

\begin{cor}[Differentiable rigidity]\lab{cor:diff_rig}Given $\a\in\S$, there exists a set ${\cal A}_\a$ of matrices in $SL(2,{\R})$ which is of total Haar measure such that for any $A_0\in {\cal A}$ the following is true: if for $A(\cdot)\in C^\infty({\T},SL(2,{\R}))$ the cocycle $(\a,A)$ is $C^0$-conjugated to the constant cocycle $(\a,A_0)$ then the conjugacy is $C^\infty$. 
\end{cor}
\begin{proof}

Let us call $B\in C^0({\T},SL(2,{\R}))$ the corresponding conjugacy so that for every $\th\in{\T}$
\begin{equation}B(\th+\a)A_0B(\th)^{-1}=A(\th)\lab{1}\end{equation} 
Up to a set of Haar measure zero (corresponding to parabolic matrices) there is basically two cases to consider: either $A_0$ is hyperbolic and the result which  is a consequence of the (uniform) hyperbolicity in the fibers is  well known; or $A_0$ is elliptic, and we can assume that its fibered rotation number is diophantine w.r.t $\a$ (this is a total Haar measure condition in the set of elliptic matrices). The hypothesis of Main theorem A are then satisfied. Therefore there exists $\ti B\in C^\infty({\T},SL(2,{\R}))$ and $\ti A_0$ such that for every $\th\in {\T}$
\begin{equation}\ti B(\th+a)\ti A_0 {\ti B}(\th)^{-1}=A(\th).\lab{2}\end{equation}
We can assume (we just have to make a conjugacy by a constant matrix in $SL(2,{\R})$) that $A_0$ is a rotation matrix of angle $\r_f(\a,A_0)$
The rotation numbers of $(\a,A_0)$ and $(\a,\ti A_0)$ satisfy
$$\r_f(\a,\ti A_0)=\r_f(\a,A_0)+r\a\mod 1$$
where $r$ is the degree of $(\ti B)^{-1}B$, so that there exists $P\in SL(2,{\R})$ such that
$$\ti A_0=PE_{r}(\th+\a) A_0E_{r}(\th)^{-1} P^{-1}$$ with $P\in SL(2,{\R})$. Consequently, if
$$D(\th)=B^{-1}(\th)\ti B(\th)PE_{r}(\th)$$
(\ref{1}) and (\ref{2}) show that
$$D(\th+\a)=A_0D(\th)A_0^{-1},$$
and using Fourier coefficients
$$\left(\Ad(A_0)-e^{2\pi ik\a}Id\right).\hat D(k)=0.$$
But since $\r(A_0)$ is diophantine w.r.t $\a$ this means that $D(\cdot)$ is constant (hence $C^\infty$), and $B$ is then $C^\infty$.  
\end{proof}

\bn{\bf Question:} Is corollary~\ref{cor:diff_rig} same theorem true if $B$ is only assumed to be $L^2$ or $L^\infty$...

The second application of theorem~\ref{main:theo} (A$+$B)   is maybe more striking and concerns density of systems having {\it positive Lyapunov exponent}. Questions concerning positivity of Lyapunov exponents are usually difficult and most of the time deal with perturbative situations and with particular one-parameter families (see for example the very nice paper~\cite{B-G} to have an idea of the difficulty involved in that kind of problems).  We emphasize on the fact that our theorem is {\it non-pertubative}.

\begin{cor}[Density of positive Lyapunov exponent]\lab{dple}If $\a\in\S$ is fixed, then there is a dense set of    $A(\cdot)$ in $C^\infty({\T},SL(2,{\R}))$  for which the fibered Lyapunov exponent of $(\a,A)$ is positive. 
\end{cor}
\begin{proof}

We refer the reader to the appendix for the notion of complex rotation number and its properties.
For $z\in{\bf D}-\{0\}$ we define
$$C_z=\begin{pmatrix}{z+z^{-1}\ov 2}& -{z-z^{-1}\ov 2i}\\
{z-z^{-1}\ov 2i}& {z+z^{-1}\ov 2}\end{pmatrix}$$
When $z=e^{i\b}\in \pa {\bf D}$ we denote by $R_\b$ the corresponding matrix $C_z$; it  is in $SL(2,{\R})$ and it is a rotation matrix.
The complex fibered rotation number $w(z)$ of $(\a,A_z(\cdot))$
where $A_z(\cdot)=A(\cdot)C(z)$ satisfies:
\begin{itemize}
\item a) $z\mapsto w(z)$ is holomorphic for $z\in{\bf D}-\{0\}$;
\item b) for Lebesgue-a.e. $\b\in{\bf S}^1$ , ${\rm Re} w(e^{i\b})$ is the fibered-Lyapunov exponent of $(\a,A(\cdot)R_\b)$.
\item for all $\b\in\pa{\bf D}$, $\Im w(e^{i\b})$ is the fibered rotation number of   $(\a,A(\cdot)R_\b)$ and is continuous w.r.t $\b$;
\item c) if ${\rm Re} w(z)=0$ a.e on an open arc of $\pa {\bf D}$ then $w(\cdot)$ can be holomorphically extended through the arc  and then 
$$\sup_{n\in{\Z}}\max_{\th\in{\T}}\|A_n(\th)\|<\infty$$
(the iterates of $(\a,A)$ are $C^0$-bounded);
\item d) if $A(\cdot)$ is homotopic to the identity and  c) holds then on the open arc where ${\rm Re} w(\cdot)\equiv 0$, $\Im w(\cdot)$ cannot be constant.
\end{itemize}

So given $(\a,A)$ we consider $(\a,A(\cdot)R_\b)$ for $\b$ in an arbitrary small interval $(-\d,\d)$. Either there exixts $\b_0\in (-\d,\d)$ for which the Lyapunov exponent of $(\a,A(\cdot)R_\b)$ is positive and we are done, or the iterates of $(\a,A)$ are $C^0$-bounded (cf. item c)). The same is also true for the iterates of $(\a,A(\cdot)R_\e)$ for any small $\e$. 
At this point we distinguish two cases:
\begin{itemize}
\item  i) either  $A(\cdot)$ is homotopic to the identity: and 
by item d) the map $\e\mapsto \Im w(e^{i\e})$ is never constant on any interval centered at 0 and continuous. Therefore one can choose $\e$ small enough so that the fibered rotation number of $(\a,A(\cdot)R_\e)$ is diophantine w.r.t $\a$ and moreover the iterates of $(\a,A(\cdot)R_\e)$ are $C^0$-bounded. The Main Theorem A then applies and proves that $(\a,A(\cdot)R_\e)$ is $C^\infty$-conjugated to a constant elliptic  system on ${\R}/{\Z}$. Now by proposition~\ref{loc.sit}, arbitrarily $C^\infty$-close to $A(\cdot)R_\e$ there are systems with positive fibered Lyapunov exponent.
\item ii) or $A(\cdot)$ is of degree $r\in{\Z}-\{0\}$: and we know from the Main theorem part B that $(\a,A(\cdot))$ is conjugated to $(\a,E_r(\cdot))$. But now using a result by M. Herman (or an improved version by A. Avila and J. Bochi) see~\cite{He},~\cite{A-B} for any matrix $D\in SL(2,{\R})$ arbitrarily close to the identity and which does not commute with the rotations $E_r(\th)$, the Lyapunov exponent of $(\a,E_r(\cdot)D)$ is positive. We now conjugate back to find a system $(\a,\ti A(\cdot))$ with positive Lyapunov exponent, $\ti A$ being as close as we want from $A(\cdot)$, to conclude the proof.
\end{itemize} 
\end{proof}

Let $V:{\T}\to {\R}$ be a smooth function. We introduce  the  Schr\"odinger cocycle $(\a,A_V(\cdot))$ where
$$A_V(\cdot)=\begin{pmatrix}V(\cdot) &1\\-1&0\end{pmatrix}.$$
\begin{cor}[Schr\"odinger cocycles]If $\a\in\S$ is fixed, the set of $V:{\T}\to{\R}$ for which $(\a,A_V(\cdot))$ has positive Lyapunov exponent is $C^\infty$-dense.
\end{cor}
For a proof see section~\ref{sec:schrod}.

The preceding results for discrete cocycles have a counterpart for differential equations. Given $u:{\R}^2/{\Z}^2\to sl(2,{\R})$ a smooth function and $\a\in{\R}$ we introduce the differential equation on ${\R}^2/{\Z}^2\times SL(2,{\R})$ ($(\th,X)\in{\R}^2/{\Z}^2\times SL(2,{\R})$) 
$$\begin{cases}&\dot X=u(\th)X\\
&\dot\th=(1,\a)
\end{cases}$$
Also if $v:{\T}^2\to {\R}$ and if 
$$u(\th)=\begin{pmatrix}0&-1\\v(\th)&0\end{pmatrix}$$
we call this O.D.E a Schr\"odinger (or Hill) equation with potential $v$.
\begin{cor}[Differential equations]If $\a\in\S$ is fixed, the set of $u:{\T}^2\to sl(2,{\R})$ for which the previous O.D.E has positive Lyapunov exponent is $C^\infty$-dense. Also the set of $v:{\T}^2\to{\R}$ for which the corresponding Schr\"odinger equation with potential $v$ has positive Lyapounov exponent is $C^\infty$-dense.
\end{cor}
\begin{proof}

We just give the proof of the first part of the corollary (the proof of the second part goes along the same lines). Let us denote by ${\cal Z}^t$ the flow of the differential equation on ${\T}^2\times SL(2,{\R})$ which is of the form: ${\cal Z}^t(\th_1,\th_2,y)=(\th_1+t,\th_2+t\a,Z^t(\th_1,\th_2)y)$. We introduce $A(\th)=Z^1(0,\th)$ which is clearly defined on ${\T}$. To say that ${\cal Z}^t$ has positive Lyapunov exponent is equivalent to the same statement concerning $(\a,A)$. The theory of the complex rotation number we have mentionned is available in the ODE context (see~\cite{C-J}) and proceding as in the proof of corollary~\ref{dple} we can assume by contradiction that $A(\cdot)$ has bounded fibred products and diophantine rotation number w.r.t $\a$ and hence is $C^\infty$-reducible. This means that there is a complex section $\s:{\T}\to {\bf C}^2-\{0\}$ and a real number $\r$ such that
$$A(\cdot)\s(\cdot)=e^{2\pi i \r}\s(\cdot+\a).$$
We now extend the section $\s(\cdot)$ to the torus ${\T}^2$ the following way: for $(x_1,x_2)\in {\R}^2$ we define:
$$\ti \s(x_1,x_2)=e^{-2\pi i x_1\r}Z^{x_1}(0,x_2-\a x_1)\s(x_2-\a x_1).$$
This function is clearly smooth and 1-periodic in $x_2$. Let us check that it is 1-periodic in $x_1$:
\begin{align*}\ti \s(x_1+1,x_2)&=e^{-2\pi i (x_1\r+\r)}Z^{x_1+1}(0,x_2-\a x_1-\a)\s(x_2-\a x_1-\a)\\
&=e^{-2\pi i x_1\r}e^{-2\pi i\r}Z^{x_1}(1,x_2-\a x_1)\\
&\ \ \ Z^1(0,x_2-\a x_1-\a)\s(x_2-\a x_1-\a)\\
&=e^{-2\pi i\r x_1}Z^{x_1}(1,x_2-\a x_1)\s(x_2-\a x_1)\\
&=\ti \s(x_1,x_2).
\end{align*}
Also by definition we have
$$Z^t(\th_1,\th_2)\ti\s(\th_1,\th_2)=e^{2\pi i\r}\ti \s(\th_1+t,\th_2+t\a).$$
The fibered flow $(t,t\a,Z^t(\cdot))$ is then reducible which means that for some $B:{\T}^2\to SL(2,{\R})$ and $U_0\in sl(2,{\R})$ one has ($(\w_1,\w_2):=(1,\a)$):
$$(\sum_{i=1}^2\w_i\pa_{\th_i}B(\th_1,\th_2))B(\th_1,\th_2)^{-1}+\Ad(B(\th_1,\th_2).U_0=U(\th_1,\th_2).$$
 and it is not difficult to see that there are abritrary $C^\infty$-small  perturbations of $U_0$ that make the fibered flow $(1,\a,Z^t(\cdot))$ hyperbolic

\end{proof}

\section{${\bf Z}^2$-actions and renormalization}\lab{sec:3}
\paragraph*{Actions}
Assume we are given a ${\Z}^2$-action ${\cal A}$ on ${\R}\times SL(2,{\R})$; we denote the action of $(n,m)\in{\Z}^2$ on $(t,y)\in{\R}\times SL(2,{\R})$ by
$(n,m)_{\cal A}.(t,y)$. We say that the action ${\cal A}$ is a {\it fibered ${\Z}^2$-action} if for every $(n,m)\in{\Z}^2$ there exist a real number $\g^{\cal A}_{n,m}$ and a smooth map $D_{n,m}^{\cal A}:{\R}\to SL(2,{\R})$ such that
$$\forall (t,y)\in {\R}\times SL(2,{\R})),\ (n,m)_{\cal A}.(t,y)=(t+\g_{n,m},D_{n,m}^{\cal A}(t)y).$$
When all the maps $D_{n,m}^{\cal A}$ are constant maps, we say that the action ${\cal A}$ is {\it constant}.
Two fibered ${\Z}^2$ actions ${\cal A},{\cal A'}$ are said to be {\it fibered-conjugated} (or {\it $f$-conjugated} or for short {\it conjugated}) if there exists a smooth map $B:{\R}\to SL(2,{\R})$ such that
$$\forall(n,m)\in{\Z}^2\ (n,m)_{\cal A'}=(0,B)\circ(n,m)_{\cal A}\circ(0,B)^{-1}$$
that is if
$$\begin{cases}&D_{n,m}^{\cal A'}(t)=B(t+\g_{n,m}^{\cal A})D_{n,m}^{\cal A}(t)B(t)^{-1}\\
&\g_{n,m}^{\cal A'}=\g_{n,m}^{\cal A}\end{cases}$$
A fibered action is said to be {\it reducible} if it is conjugated to a constant action. 

There is another notion of conjugation (we shall call $d$-conjugation to avoid confusion with $f$-conjugation): if $\u,\l\in{\R}$ with $\l\ne 0$ we define the dilatation
$$L_\l.(t,y)=(\l t,y),$$
and the translation
$$T_\u.(t,y)=(t-\u,y).$$
If we define
 $$\forall(n,m)\in{\Z}^2\ (n,m)_{\cal A'}=L_\l^{-1}\circ(n,m)_{\cal A}\circ L_\l$$
we get
$$\begin{cases}D_{n,m}^{\cal A'}(t)&=D_{n,m}^{\cal A}(\l t)\\ 
\g_{n,m}^{\cal A'}&=\l^{-1}\g_{n,m}^{\cal A}\end{cases},$$

\paragraph*{Bases}
If $(e_1,e_2)$ is a basis of the ${\Z}$-module ${\Z}^2$ any other basis is of the form $(ae_1+be_2,ce_1+de_2)$ where the matrix
$$\begin{pmatrix}a&b\\c&d\end{pmatrix}\in GL(2,{\Z}).$$
We say that a fibered action ${\cal A}$ is {\it normalized} if 
$$(1,0)_{\cal A}.(t,y)=(1,Id).(t,y)=(t,y).$$
We have the following proposition the proof of which is similar to that of~\cite{K_annsmaths} Proposition 4.1:
\begin{prop}\lab{prop:norm}Every action is conjugated to a normalized one.
\end{prop}

Also, if there exists a basis $(e_1,e_2)$ of the ${\Z}$-module ${\Z}^2$ and a smooth map $B:{\R}\to SL(2,{\R})$ such that   {\it simultaneously} 
$$(0,B)\circ e_{1,{\cal A}}\circ (0,B)^{-1},\qquad (0,B)\circ e_{2,{\cal A}}\circ (0,B)^{-1}$$
are constant cocycles then the action ${\cal A}$ is conjugated to a constant one.

\paragraph*{Module of frequencies}
For further purpose we define the module of frequencies of the action ${\cal A}$ as being the sub ${\Z}$-module of ${\R}$ generated by $\g^{\cal A}_{e_1},\g^{\cal A}_{e_2}$ which is the same as the one generated by $\g^{\cal A}_u,\g^{\cal A}_v$ where $(u,v)$ is  any other base of ${\Z}^2$.
We denote it by $M_{\cal A}$. 
The {\it half frequency module} ${1\ov 2}M_{\cal A}$ will be important later.
Notice that if $\ti{\cal A}=(0,B)\circ{\cal A}\circ (0,B)^{-1}$
then $M_{\cal A}=M_{\ti{\cal A}}$.
Moreover since $L_\l:{\R}\to {\R}$ is an isomorphism and $L_\l:M_{\cal A}\to M_{\ti{\cal A}}$ we can define an isomorphism of ${\Z}$-module $l_\l:{\R}/M_{\cal A}\to {\R}/M_{\ti{\cal A}}$ (and also  $l_\l:{\R}/{1\ov 2}M_{\cal A}\to {\R}/{1\ov 2}M_{\ti{\cal A}}$).

We say that an element $\b\in{\R}$ is {\it $\t$-diophantine w.r.t} to the half frequency module ${1\ov 2}M_{\cal A}$ if there exist a basis $(e_1,e_2)$ of ${\Z}^2$ and  a constant  $K>0$ such that for  any $(k_1,k_2)\in{\Z}^2-\{(0,0)\}$
$$|\b-{1\ov 2}k_1\g_{e_1}-{1\ov 2}k_2\g_{e_2}|\geq \frac{K^{-1}}{(|k_1|+|k_2|)^\t}.$$
Notice that if this holds for one basis $(e_1,e_2)$ then it holds for any other basis (the constant $K$ has then to be changed).

\paragraph*{Example, link with quasi-periodic cocycles}
If $A,C:{\R}\to SL(2,{\R})$ satisfy for every $t\in{\R}$
$$C(t+\a)A(t)=A(t+1)C(t)$$
one has
$$(1,C)\circ (\a,A)=(\a,A)\circ(1,C)$$
and one can define the fibered-${\Z}^2$ action ${\cal A}$ on ${\R}\times SL(2,{\R})$ by
$$(n,m)_{\cal A}.(t,y)=(1,C)^n\circ(\a,A)^m.(t,y).$$
We shall then write ${\cal A}=((1,C),(\a,A))$.
Notice that when $C(\cdot)\equiv Id$ the map $A(\cdot)$ is automatically ${\Z}$-periodic.

In particular if $(\a,A)$ is a quasi-periodic cocycle  with $A:{\R}/{\Z}\to SL(2,{\R})$ and $\a$ irrational one can define the fibered ${\Z}^2$-action on ${\R}\times SL(2,{\R})$ defined by
\begin{align*}(n,m)_{\cal A}.(t,y)=&(1,Id)^n\circ(\a,A)^m(t,y)\\
=&(\a,A)^m\circ (1,Id)^n(t,y).\end{align*}

The notion of reducibility is illustrated as follows:
The action ${\cal A}=((1,C),(\a,A))$ is reducible if there exist $A_0,C_0$ in $SL(2,{\R})$ and a smooth map $B:{\R}\to SL(2,{\R})$ such that {\it simultaneously}
\begin{align*}C(t)&=B(t+1)C_0B(t)^{-1}\\
A(t)&=B(t+\a)A_0B(t)^{-1}.\end{align*}
Solving {\it one} of these two equations (for example the first one) is an easy task (at least in the smooth category), and in that case one can even choose $C_0=Id$: this is proposition~\ref{prop:norm}.

The link between reducibility of quasi-periodic cocycles and reducibility of fibered actions is just the following trivial but fundamental observation: the quasi-periodic cocycle $(\a,A)$ (defined on ${\T}\times SL(2,{\R})$) is reducible if and only if the action ${\cal A}=((1,C),(\a,A))$ (defined on ${\R}\times SL(2,{\R})$) is reducible.

Since every ${\Z}^2$-action can be normalized, the reducibility problems for cocycles or for actions are equivalent.

\paragraph*{Renormalization}
The gain in proceding so is that one can perform new operations namely:
\begin{itemize}
\item change of basis: given an action $((1,C),(\a,A))$ and a matrix 
$$\begin{pmatrix}a&b\\c&d\end{pmatrix}\in GL(2,{\Z})$$ one can look at the action $\ti{\cal A}$ defined by
$$(n,m)_{\ti{\cal A}}=(an+bm,cn+dm)_{\cal A};$$
\item $d$-conjugacies
\item $f$-conjugacies
\end{itemize} 
\paragraph{Continued fraction expansion}
We refer to~\cite{Y} for a more detailled exposition of the following facts.

\bn Define as usual for $0<\a<1$,
$$ a_0=0,\qquad \a_{0}=\a,$$
and inductively for $k\geq 1$,
$$a_k=[\a_{k-1}^{-1}],\qquad \a_k=\a_{k-1}^{-1}-a_k=G(\a_{k-1})=\{{1\ov \a_{k-1}}\},$$
where $[\ ]$ denotes the integer part and $G(\cdot)$  the fractional part (the gauss map).
We also set,
$$\b_k=\prod_{j=0}^k\a_j,$$
\begin{align*}&p_0=0\qquad q_1=a_1\\
&q_0=1\qquad p_1=1,\end{align*}
and inductively,
\begin{equation}\lab{3.1}\begin{cases}&p_k=a_kp_{k-1}+p_{k-2}\\
&q_k=a_kq_{k-1}+q_{k-2}.\end{cases}\end{equation}
Then, we have the following useful relations,
\begin{equation}\forall k\geq 2,\qquad \b_{k-2}=a_k\b_{k-1}+\b_k,\end{equation}
\begin{equation}\forall k\geq 1,\qquad q_kp_{k-1}-q_{k-1}p_k=(-1)^k,\end{equation}
and for all $k\geq 0$,
\begin{equation}\lab{beta-k}\b_k=(-1)^k(q_k\a-p_k),\end{equation}
\begin{equation}{1\ov q_{k+1}+q_k}<\b_k<{1\ov q_{k+1}},\end{equation}
\begin{equation}\b_k={1\ov q_{k+1}+\a_{k+1}q_k}.\end{equation}

We state a classical lemma wich will be useful later:
\begin{lemme}\lab{m_1...}Let $k\geq l\geq 0$ and assume to simplify that $2|k-l$. Set
$$\begin{pmatrix}m_1&m_2\\
m_3&m_4\end{pmatrix}=\begin{pmatrix}0&1\\1&-a_k\end{pmatrix}\cdots\begin{pmatrix}0&1\\1&-a_{l+1}\end{pmatrix}.$$
Then $m_1,m_4\geq 0$, $m_3,m_2\leq 0$ and
$$\begin{pmatrix}|m_4|&|m_2|\\
|m_3|&|m_1|\end{pmatrix}\begin{pmatrix}1\\\a_k\end{pmatrix}=\frac{\b_{l-1}}{\b_{k-1}}\begin{pmatrix}1\\\a_l\end{pmatrix}$$ 
If $l=0$ then
$$\begin{pmatrix}m_1&m_2\\
m_3&m_4\end{pmatrix}=\begin{pmatrix}p_{k-1}&-q_{k-1}\\
-p_k&q_k\end{pmatrix}.$$
\end{lemme}
\begin{proof}

If we take the inverse we get 
$$\begin{pmatrix}m_4&-m_2\\
-m_3&m_1\end{pmatrix}=\begin{pmatrix}a_{l+1}&1\\1&0\end{pmatrix}\cdots\begin{pmatrix}a_k&1\\1&0\end{pmatrix},$$
and this matrix has non negative coefficients. 
To conclude just observe that
$$\begin{pmatrix}a_j&1\\1&0\end{pmatrix}\begin{pmatrix}1\\\a_j\end{pmatrix}=\a_{j-1}^{-1}\begin{pmatrix}1\\\a_{j-1}\end{pmatrix}.$$
The second part of the lemma is also classical and checked by induction

\end{proof}

\paragraph{The set $\S$}\lab{sec:thesetsigma}
If $\g$ is a positive number and $\s>1$, we define the set,
$$CD(\g,\s)=\{\a\in{\R}/{\Z},\ \forall k\in{\Z}-\{0\},\ \min_{l\in{\Z}}|k\a-l|\geq {\g^{-1}\ov |k|^{\s}}\}.$$
These are classical facts  that if $\g>0$ is fixed, such a set is of positive Haar measure on ${\T}^1$, that this measure goes to 1 when $\g$ goes to infinity and that the union,
$$\bigcup_{\g>0}CD(\g,\s),$$
is of total Haar measure. We shall fix $\s>1$.  Since the map  $G:(0,1)\to (0,1)$ preserves an absolutely continuous measure $m$ (the Gauss measure $\frac{1}{\Log 2}\frac{dx}{1+x}$), we have, provided $\g$ is chosen large enough,
$$m(CD(\g,\s)\cap G^{-1}(CD(\g,\s)))>1-m(({1\ov 5},{1\ov 4}]\cap G^{-1}(({1\ov 5},{1\ov 4}])),$$
(since the set in the right hand side has positive $m$-measure)
and consequently the set $\S_1$ equal to
$$(CD(\g,\s)\cap((1/5),(1/ 4)])\cap G^{-1} (CD(\g,\s)\cap((1/5),(1/ 4)])$$
 is of positive $m$-measure. Since $G$ is $m$-ergodic, the set $\S$ of $\a\in(0,1)$ such that for any $N\in{\N}$ there exists $k\geq N$ for which $G^k(\a)\in\S_1$, that is for which
$$(\a_k,\a_{k+1})\in (CD(\g,\s)\cap({1\ov 5},{1\ov 4}])^2,$$
is of total $m$-measure and hence of total Lebesgue measure (since $m$ is absolutely continuous).

Notice that for any $\a$ such that $\a_k\in CD(\g,\s)\cap((1/5),(1/ 4)]$ one has $a_{k+1}=4$ where $\a_k=G^k(\a)$.

\subsection{Degree and rotation numbers of q.p.f ${\Z}^2$-actions}
In this section we intend to define the notion of degree and of rotation number of a q.p.f ${\Z}^2$-action. Given such an action ${\cal A}$,
there is a corresponding action on ${\R}\times {\bf S}^1$
$$(n,m)_{\cal A}:(t,\bar x)\mapsto (t+\g^{\cal A}_{n,m},f_{n,m}^{\cal A}(t,\bar x)).$$
For each $(n,m)\in{\Z}^2$ we can lift $(n,m)_{\cal A}$ to ${\R}\times{\R}$
$$(n,m)_{\cal A}.(t,x)=(t+\g^{\cal A}_{n,m},x+d_{n,m}^{\cal A}(t,x))$$
and we have
$$\pi (x+d_{n,m}^{\cal A}(t,x))=f_{n,m}^{\cal A}(t,\pi(x)).$$ 
Notice that the map $d_{n,m}^{\cal A}(t,x)$ is not uniquely defined, every other lift being of the form $d_{n,m}^{\cal A}+k_{n,m}$ where $k_{n,m}$ can be any rational integer. Notice also that since for each fixed $t$ the map $\bar x\to f(t,\bar x)$ is of degree 1 (these are homographies), a lift $x+d_{(n,m)}^{\cal A}(t,x)$ enjoys the property
$$d_{(n,m)}^{\cal A}(t,x+1)=d_{(n,m)}^{\cal A}(t,x).$$

\paragraph*{Degree}
Assume $(e_1,e_2)$ is a base of the ${\Z}$-module ${\Z}^2$ and set 
$$(e_1)_{\cal A}=(\g_1,C(\cdot)),\qquad (e_2)_{\cal A}=(\g_2,A(\cdot)).$$
In order to simplify the notations we define
$$c(t,x)=d_{e_1}(t,x),\qquad a(t,x)=d_{e_2}(t,x)$$
$$F_c(t,\bar x)=(e_1)_{\cal A}(t,\bar x)\qquad F_c(t,\bar x)=(e_2)_{\cal A}(t,\bar x).$$
Since $(\g_1,C(\cdot))$ and $(\g_2,A(\cdot))$ commute we can say that
\begin{align*}\deg(e_1,e_2)&=(a\circ F_c+c)-(c\circ F_a+a)\\
&=(a\circ F_c-a)-(c\circ F_a-c)\end{align*}
is a constant rational integer  which does not depend on the choice of the lifts  for $a,c$

We now show that $|\deg(e_1,e_2)|$ does not depend on the choice of the base $(e_1,e_2)$. Since $GL(2,{\Z})$ is generated by the matrices
\begin{equation}\lab{sl2gen}\begin{pmatrix}1&0\\0&-1\end{pmatrix},\quad\begin{pmatrix}0&-1\\1&0\end{pmatrix},\quad \begin{pmatrix}1&1\\0&1\end{pmatrix}\end{equation}
it is enough to check  the following three equalities: $\deg(-e_1,-e_2)=-\deg(e_1,e_2)$  and $$\deg(e_1,e_2)=\deg(-e_2,e_1),\quad \deg(e_1,e_2)=\deg(e_1+e_2,e_2).$$
If we introduce the notation $a^{(-1)}=d_{-e_2}$, one sees that 
$$a\circ F_a^{-1}+a^{(-1)}$$
is a constant rational integer and hence
\begin{align*}\deg(e_1,-e_2)&=(a\circ F_c^{-1}+c^{(-1)})-(c^{(-1)}\circ F_a+a)\\
&=a\circ F_c^{-1}-c\circ F_c^{-1}+c\circ F_c^{-1}\circ F_a-a\\
&=(c\circ F_a+a)\circ F_c^{-1}-(a\circ F_c+c)\circ F_c^{-1}\\
&=-deg(e_1,e_2),\end{align*}
(we have used the fact that $F_a$ and $F_c$ commute) which proves the first equality. The second one is obtained the same way:
\begin{align*}\deg(-e_2,e_1)&=(c\circ F_a^{-1}-c)-(a^{(-1)}\circ F_c-a^{(-1)})\\
&=(c\circ F_a^{-1}-c)-(-a\circ F_a^{-1}\circ F_c+a\circ F_a^{-1})\\
&=(c\circ F_a^{-1}-c)-(-a\circ F_c\circ F_a^{-1}+a\circ F_a^{-1})\\
&=[-(c\circ F_a-c)+(a\circ F_c-a)]\circ F_a^{-1}\\
&=\deg(e_1,e_2).\end{align*}
Similarly, up to some integer additive constant
$$d_{e_1+e_2}^{\cal A}=a\circ F_c+c$$
\begin{align*}\deg(e_1+e_2,e_2)&=(a\circ F_a\circ F_c-a)-((a\circ F_c+c)\circ F_a-(a\circ F_c+c))\\
&=c-c\circ F_a-a+a\circ F_c\\
&=\deg(e_1,e_2)
\end{align*}
and the third equality is satisfied.

We can thus define the degree of the action ${\cal A}$ as being  $|\deg(e_1,e_2)|$. Notice that if the action is normalized ($(1,0)_{\cal A}=(1,Id)$) we recover the usual notion of the degree of a map ${\R}/{\Z}\to SL(2,{\R})$ (up to the sign).
Notice also that if one allows only changes of basis of determinant 1 then the absolute value is not needed.
\begin{lemme}The degree of a q.p.f action is invariant by conjugation.
\end{lemme}
\begin{proof}

Let us first consider the case of $f$-conjugacies.
Let us denote by ${\cal A}'$ the conjugate action by $(0,B)$ that is $(n,m)_{{\cal A}'}=(0,B)\circ(n,m)_{{\cal A}}\circ (0,B)^{-1}.$ We denote $\ti a=d_{e_2}^{{\cal A}'}$, $\ti c=d_{e_1}^{{\cal A}'}$...
Take a lift $(0,b(t,x)).$ Then 
$$\deg({\cal A}')=(\ti a\circ F_{\ti c}-\ti a)-(\ti c\circ  F_{\ti a}-\ti c)$$
and using the fact that up to some integer constant 
$$\ti a=b\circ F_a\circ F_b^{-1}+a\circ F_b^{-1}-b\circ F_b^{-1}$$
$$\ti c=b\circ F_c\circ F_b^{-1}+c\circ F_b^{-1}-b\circ F_b^{-1}$$
and the fact that $F_{\ti c}=F_b\circ F_c\circ F_b^{-1}$, $ F_{\ti a}=F_b\circ F_a\circ F_b^{-1}$
we get after some calculation (and the fact that $F_a,F_c$ commute) the desired conclusion.

The case of $d$-conjugacies is treated the same way.

\end{proof}

\paragraph*{Rotation number}
We assume in that section that the degree of the action ${\cal A}$ is zero.
With the previous notations this implies that the following equality holds
$$a\circ F_c+c=c\circ F_a+a.$$
We introduce the set ${\cal M}_{Leb}$ of measures on ${\R}\times {\bf S}^1$ that project on the first factor to Lebesgue measure on ${\R}$. We claim that there exists a non trivial invariant measure $\u\in{\cal M}_{Leb}$ defined on ${\R}\times {\bf S}^1$ invariant by $F_a,F_c$ since  any q.p.f ${\Z}^2$-action can be normalized (and thus we are reduced to the classical case of the existence of an invariant measure for a diffeomorphism on a compact space). 
If $\g_1,\g_2$ are real numbers and $f$ is a real valued map defined on ${\R}\times{\bf S}^1$ we introduce the notation
\begin{equation*}I_{t_1}^{t_2}f=\begin{cases}+\int_{[t_1,t_2]\times {\bf S}^1}f(t,x)d\u(t,x)\qquad{\rm if}\quad t_1\leq t_2\\
-\int_{[t_2,t_1]\times {\bf S}^1}f(t,x)d\u(t,x)\qquad{\rm if}\quad t_1> t_2
.\end{cases}\end{equation*}
Notice that for any $t_1,t_2,t_3$
$$I_{t_1}^{t_3}f=I_{t_1}^{t_2}f+I_{t_2}^{t_3}.$$
Also, since $\u$ is invariant by $F_a,F_c$ we get
\begin{equation}\lab{form}I_{t_1}^{t_2}f\circ F_a=I_{t_1+\g_a}^{t_2+\g_a}f,\qquad I_{t_1}^{t_2}f\circ F_c=I_{t_1+\g_c}^{t_2+\g_c}f.\end{equation}
We then define
$${\rm rot}_\u(e_1,e_2)=I_0^{\g_c}a-I_0^{\g_a}c.$$
Notice that if we choose other lifts for $a,c$ this number changes by the addition of an element of ${\Z}\g_a+{\Z}\g_c$ which means that we can define the element ${\rm rot}_\u(e_1,e_2)$ as an element of ${\R}/M_{\cal A}$ where ${\cal A}$ is the module of frequencies of the action ${\cal A}$. But if we do so the rotation number is only defined up to the sign. It is hence  more convenient to introduce its class modulo $(1/2)M_{\cal A}$ the {half frequency module}.
We  now check that its class  does not depend on the choice of the base $(e_1,e_2)$ by proving like in the previous section 
$$\rot_\u(e_1,-e_2)\equiv -\rot_\u(e_1,e_2)\mod M_{\cal A},$$
$$\rot_\u(-e_2,e_1)\equiv \rot_\u(e_1,e_2),\qquad \rot_\u(e_1+e_2,e_2)\equiv \rot_\u(e_1,e_2)\mod M_{\cal A}.$$
Let us check the first equality, that is, using the fact that
$a^{(-1)}\equiv -a\circ F_a^{-1}\mod{\Z}$ and that $\u$ projects to Lebesgue measure on ${\R}$
\begin{align*}\rot_\u(e_1,-e_2)&\equiv I_0^{\g_c}a^{(-1)}-I_0^{-\g_a}c\mod M_{\cal A}\\
&\equiv -I_{-\g_a}^{\g_c-\g_a}a+I_{-\g_a}^0c\mod M_{\cal A}\\
&\equiv I_{-\g_a}^0(c-a)-I_0^{\g_c-\g_a}a\mod M_{\cal A}\end{align*}
and since 
$$a(t)-c(t)=a\circ F_c-c\circ F_a$$
this gives
\begin{align*}\rot_\u(e_1,-e_2)&\equiv -I_{\g_c-\g_a}^{\g_c}a+I_0^{\g_a}c-I_0^{\g_c-\g_a}a\mod M_{\cal A}\\
&\equiv -\rot_\u(e_1,e_2)\mod M_{\cal A}.\end{align*}

For the second equality the computations are the same:
\begin{align*}\rot_\u(-e_2,e_1)&\equiv I_0^{-\g_a}c-I_0^{\g_c}a^{(-1)}\mod M_{\cal A}\\
&\equiv -I_{-\g_a}^0c+I_{-\g_a}^{\g_c-\g_a}a\mod M_{\cal A}\\
&\equiv I_{-\g_a}^0(a-c)+I_0^{\g_c-\g_a}a\mod M_{\cal A}\end{align*}
and since 
$$a(t)-c(t)=a\circ F_c-c\circ F_a$$
this gives
\begin{align*}\rot_\u(-e_2,e_1)&\equiv I_{\g_c-\g_a}^{\g_c}a-I_0^{\g_a}c+I_0^{\g_c-\g_a}a\mod M_{\cal A}\\
&\equiv \rot_\u(e_1,e_2)\mod M_{\cal A}.\end{align*}
Similarly
\begin{align*}\rot_\u(e_1+e_2,e_2)&\equiv I_0^{\g_a}(a\circ F_c+c)-I_0^{\g_c+\g_a}a\mod M_{\cal A}\\
&\equiv I_0^{\g_a}c+I_{\g_c}^{\g_c+\g_a}a-I_0^{\g_c+\g_a}\mod M_{\cal A}\\
&\equiv \rot_\u(e_1,e_2)\mod M_{\cal A}.\end{align*}
The previous computations show that we can define  $\rot_\u({\cal A})$ in ${\R}/{1\ov 2}M_{\cal A}$.

We now investigate the behavior of this rotation number with respect to conjugacies.
\begin{lemme}If $\u\in{\cal M}_{Leb}({\cal A})$ and $\ti{\cal A}=(0,B)\circ {\cal A}\circ (0,B)^{-1}$ is an  $f$-conjugated action then $\ti\u=(0,B)_*\u$ is in ${\cal M}_{Leb}(\ti {\cal A})$ and the following equality holds in ${\R}/{1\ov 2}M_{{\cal A}}$
$$\rot_\u({\cal A})=\rot_{\ti\u}(\ti{\cal A}).$$
Also if $\ti{\cal A}=L_\l\circ {\cal A}\circ L_\l^{-1}$ is a  $d$-conjugated action then $\ti\u=(1/\l)(L_\l)_*\u$ is in ${\cal M}_{Leb}(\ti {\cal A})$ and the following equality holds:
$$\rot_\u(\ti{{\cal A}})=l_\l(\rot_{\ti\u}({\cal A})),$$
where $l_\l:{\R}/{1\ov 2}M_{{\cal A}}\to {\R}/{1\ov 2}M_{\ti{\cal A}}$ is the map we have already defined.
\end{lemme}
\begin{proof}

Denote by $b$ a lift for $B$, $\ti a$, $\ti c$ lifts for $(0,B)\circ (\g_a,A)\circ (0,B)^{-1}$ and  $(0,B)\circ (\g_c,C)\circ (0,B)^{-1}$.
We have modulo ${\Z}$
$$\ti a\equiv b\circ F_a\circ F_b^{-1}+a\circ F_b^{-1}+b^{(-1)},\qquad \ti c\equiv b\circ F_c\circ F_b^{-1}+c\circ F_b^{-1}+b^{(-1)}$$
that is
$$\ti a\equiv (b\circ F_a+a-b)\circ F_b^{-1},\qquad \ti c\equiv(b\circ F_c+c-b)\circ F_b^{-1}$$
and modulo ${1\ov 2}M_{\cal A}$:
$${\rot}_{\ti\u}(\ti{\cal A})\equiv {\rot}_\u({\cal A})+I_0^{\g_c}b\circ F_a-I_0^{\g_c}b-I_0^{\g_a}b\circ F_c+I_0^{\g_a}b;$$ 
using formulae~(\ref{form}) we get the result.
The second part of the proof is obvious.
\end{proof}

We can now prove that $\rot_\u({\cal A})$ does not depend on the choice of the invariant measure $\u\in{\cal M}_{Leb}$.  Any q.p.f ${\Z}^2$-action ${\cal A}$ can be normalized which means conjugated to an action ${\ti{\cal A}}$ such that $(1,0)_{\ti{\cal A}}=(1,Id)$  and the measure $\ti\u$ thus obtained   is invariant by $(1,Id)$. Therefore it defines a measure on ${\R}/{\Z}$ invariant by the diffeomorphism $(\th,y)\mapsto(\th+\ti\a,\ti A(\th)y)$ where $(\ti\a,\ti A)=(0,1)_{\ti{\cal A}}$. Since the degree of a q.p.f ${\Z}^2$-action is invariant by conjugation, and since we have assumed the degree of ${\cal A}$ to be 0, the degree of $\ti A(\cdot):{\R}/{\Z}\to SL(2,{\R})$ is also 0. But in that case $\rot_{\ti\u}({\ti{\cal A}})$ coincide modulo an element of ${1\ov 2}M_{\ti{\cal A}}$ with the usual notion of fibered rotation number modulo $(1/2){\Z}$, which does not depend on the choice of the invariant measure. Consequentely $\rot_\u({\cal A})\in{\R}/{1\ov 2}M_{\cal A}$ does not depend on $\u$. 

As a corollary:
\begin{lemme}For actions of degree 0,  the rotation number  is invariant by conjugation.
\end{lemme}

\subsection{The algorithm}
Let $A\in C^\infty({\R}/{\Z},SL(2,{\R}))$ and let  ${\cal A}$  the corresponding ${\Z}^2$-action. 
\begin{align*}&U=U_0=(1,0)_{{\cal A}}=(1,Id)=((\a,A)^{q_{-1}}\circ (1,Id)^{p_{-1}})^{(-1)^{-1}}\\
&V=V_0=(0,1)_{{\cal A}}=(\a,A)=((\a,A)^{q_{0}}\circ (1,Id)^{p_{0}})^{(-1)^{0}},\end{align*}
and for $k\geq 1$,
\begin{equation}\lab{3.2}\begin{cases}&U_k=V_{k-1}\\
&V_k=U_{k-1}V_{k-1}^{-a_k}.\end{cases}\end{equation}
We have just performed a change of base.

We can replace in the above construction $A(\cdot)$ by $A(\cdot-\nu)$ ($\nu\in{\T}^1$) and we shall then denote $U_k,V_k$ by $U_k(\nu),V_k(\nu)$.

It is clear from ~(\ref{3.1}) and~(\ref{beta-k}) that,
\begin{align*}&U_k=(\b_{k-1},A^{(k-1)}(\cdot))=(\a,A)^{(-1)^{k-1}q_{k-1}}=(V^{q_{k-1}}U^{-p_{k-1}})^{(-1)^{k-1}}\\
&V_k=(\b_{k},A^{(k)}(\cdot))=(\a,A)^{(-1)^{k}q_k}=(V^{q_{k}}U^{-p_{k}})^{(-1)^{k}},\end{align*}
where,
$$A^{(k)}(\cdot)=A_{(-1)^kq_{k}}(\cdot),$$
(we have used the notation $(\a,A(\cdot))^n=(n\a,A_n(\cdot))$).

Finally we set,
\begin{align*}&\ti U_k={\L}_{\b_{k-1}}^{-1}\circ U_k\circ{\L}_{\b_{k-1}}=(1,{\ti C}^{(k)}(\cdot))\\
&\ti V_k={\L}_{\b_{k-1}}^{-1}\circ V_k\circ{\L}_{\b_{k-1}}=(\a_k,{\ti A}^{(k)}(\cdot)),\end{align*}
with,
\begin{align*}&{\ti C}^{(k)}(t):=C^{(k)}(\b_{k-1}t):=A^{(k-1)}(\b_{k-1}t)\\
&{\ti A}^{(k)}(t):=A^{(k)}(\b_{k-1}t):=A^{(k)}(\b_{k-1}t).\end{align*}

\section{Rough estimates}
Given a smooth map $u:{\T}\to SL(2,{\R})$ we define
$$Lu(t)=(\pa u(t))(u(t))^{-1}$$
which is $sl(2,{\R})$-valued.

For $\g\in C^\infty({\T}^1,sl(2,{\R}))$ we define its  $C^k({\T})$ norm (where $k\geq 0$ is some integer) by,
$$\|\g\|_k=\max_{0\leq j\leq k}\max_{\th\in{\T}}\|\pa_\th^j\g(\th)\|,$$
and we shall make an extensive use of the so called {\it convexity inequalities} (or Hadamard inequalities): {\it For $0\leq j\leq k$,
\begin{equation}\|\pa_\th^j\g\|_0\leq C_k\|\g\|_0^{1-{j\ov k}}.\|\g\|_k^{j\ov k},\lab{eq:hadamard}\end{equation}
where $C_k$ is some positive constant depending on $k$.
}

The next simple proposition is one of the key observation in our way to the proof of the Main Theorem (see also~\cite{Ry}).
\begin{prop}\lab{prop:1.3.4} If $u_1,\ldots, u_N\in C^{\infty}({\T}, SL(2,{\R}))$ then,
$$L(u_1\cdots u_N)=L u_1+\Ad(u_1).L u_2+\cdots+\Ad(u_1\cdots u_{N-1}).L u_N,$$
and,
\begin{multline*}\pa(L(u_1\cdots u_N))=\sum_{i=1}^N\Ad(u_1\cdots u_{i-1}).\pa(L u_{i})+\\
+\sum_{1\leq i<j\leq N}[\Ad(u_1\cdots u_{i-1}).(L u_i),\Ad(u_1\cdots u_{j-1}).(L u_j)].\end{multline*}
Moreover if 
$$\sup_{1\leq i\leq N}\|u_1\cdots u_i\|\leq M$$
then for any integer $r\geq 0$ we have the following inequality,
$$\|\pa^r L(u_1\cdots u_N)\|_0\leq C_rM^{r+1}m_0^{r+1}(1+m_r)N^{r+1},$$
where,
$$m_0=\max_{1\leq i\leq N} \|L u_i\|_0,\ \ \ \ m_r=\max_{1\leq i\leq N\atop 0\leq k\leq r}\|\pa^k (L u_i)\|_0.$$
\end{prop}
\begin{proof}

The first two formulae are obtained by simple computations. In a similar way it is possible to give analogous expressions for the higher derivatives $\pa^k (L(u_1\cdots u_N))$, which using the fact that $\|Ad(u)\|\leq CM$ and the convexity inequalities~(\ref{eq:hadamard}), gives the last inequality. 
\end{proof}
We define
$$\g_k(t)=L(C^{(k)})(t),\qquad \eta_k(t)=L(A^{(k)})(t),$$
and
$$\ti\g_k(t)=L({\ti C}^{(k)})(t),\qquad \ti\eta_k(t)=L({\ti A}^{(k)})(t).$$

\section{Simple geometrical facts}
We denote by $sl(2,{\R})$ the Lie-algebra of $SL(2,{\R})$ that is the set of 2 by 2 matrices with real coefficients of the form
$$\begin{pmatrix}x& y+z\\
y-z& -x\end{pmatrix}$$
where $x,y,z\in{\R}$. We shall often denote such a matrix by $\{x,y,z\}$.
The quadratic form 
$$q(\{x,y,z\})=\det(\{x,y,z\})=-x^2-y^2+z^2,$$
is invariant under the $\Ad$-action of $SL(2,{\R})$ that is
$$q(\Ad(A).\{x,y,z\})=q(A\{x,y,z\}A^{-1})=q(\{x,y,z\}),\quad (A\in SL(2,{\R}))$$
and the same is true for $\kappa$ the associated bilinear form:
$$\kappa(v_1,v_2)=-x_1x_2-y_1y_2+z_1z_2,$$
($v_i=\{x_i,y_i,z_i\}$, $i=1,2$).
We  introduce ${\cal E}^+$ (resp. ${\cal E}^-$) the set of $\{x,y,z\}\in sl(2,{\R})$ such that $q(\{x,y,z\})\geq 0$ and $z\geq 0$ (resp. $z\leq 0$). 
The set ${\cal E}^{+}$ (resp. ${\cal E}^{-}$) is a  cone and is preserved by 
the $\Ad$-action (take a path connecting $A\in SL(2,{\R})$ to the identity and check that the $z$ component of $\Ad(A).\{x,y,z\}$ cannot be zero). For any $v\in {\cal E}^+$ one can then define
$$N(v)=\sqrt{q(v)}.$$
Notice that when $v,w$ are in ${\cal E}^+$, $\kappa(v,w)\geq 0$.
We also introduce the euclidian norm
$$\|\{x,y,z\}\|=x^2+y^2+z^2,$$
which is {\it not} $\Ad$-invariant.

Closely ralated to the group $SL(2,{\R})$ is the isomorphic group $SU(1,1)$ of matrices of the form
$$\begin{pmatrix}a&\bar b\\ b &\bar a\end{pmatrix}$$
with $a,b\in{\C}$ satisfying $|a|^2-|b|^2=1.$
The isomorphism between $SL(2,{\R})$ and $SU(1,1)$ is given by $A\mapsto PAP^{-1}$ where 
$$P=\begin{pmatrix}-1&-i\\-1&i\end{pmatrix}$$
The Lie algebra $su(1,1)$ of $SU(1,1)$ is just the set of matrices
$$\begin{pmatrix}it&\nu\\\bar\nu&-it\end{pmatrix}$$
with $t\in{\R}$ and $\nu\in{\C}$. We then  denote such a matrix by $\{t,\nu\}$.The associated $\Ad$-invariant quadratic form on $SU(1,1)$ is then
$$q(\{t,\nu\})=t^2-|\nu|^2.$$ 

We now recall some basic facts.

For any $v,w$ in $sl(2,{\R})$
\begin{equation}\lab{acs}q(v)q(w)=\kappa(v,w)^2+q([v,w]);\end{equation}
consequentely the following anti-Cauchy-Schwarz inequality holds for 
$v,w\in{\cal E}^+$
$$\kappa(v,w)\geq N(v)N(w),$$
where equality holds only if $v=tw$ for some $t\geq 0$ (if $w\ne 0$).
The proof of~(\ref{acs}) is just the equality
$$(x_1^2+y_1^2+y_1^2)(x_2^2+y_2^2+z_2^2)=(x_1x_2+y_1y_2+z_1z_2)^2+|y\wedge z|^2+|z\wedge x|^2+|x\wedge y|^2$$
applied with $ix_k$, $iy_k$, $z_k$ in place of $x_k,y_k,z_k$.
We can be more precise
\begin{lemme}
For any $v,w$ in ${\cal E}^+$ one has
$$\|[v,w]\|\leq 2\sqrt{\kappa(v,w)^2-q(v)q(w)}.\frac{\|v\|^2}{q(v)}.$$
\end{lemme}
\begin{proof}

For $t\in{\R}$,
$$q(tv+w)=t^2q(v)+2t\kappa(v,w)+q(w),$$
and the discriminant of this polynomial is
$$\D'=\kappa(v,w)^2-q(v)q(w).$$
Since $v,w$ lie in the cone ${\cal E}^+$, $q(tv+w)$ has always two real roots $t_+,t_-$ such that 
$$t_+-t_-=\frac{\D'}{q(v)},$$
and where $t_{\pm}v+w$ lies on ${\cal C}^{\pm}$ (the boundary of ${\cal E}^+$). If $(x_\pm,y_\pm,z_\pm)$ are the coordinates of $t_\pm v+w$ we have
$$z_+=\sqrt{x_+^2+y_+^2},\qquad z_-=-\sqrt{x_-^2+y_-^2}$$
so that
$$z_+-z_-=\sqrt{x_+^2+y_+^2}+\sqrt{x_-^2+y_-^2}$$
But this is the $z$-coordinate of $(t_+-t_-)v$ and thus
$$0<\max(z_+,|z_-|)\leq z_+-z_-\leq (t_+-t_-)\|v\|.$$
We then get
\begin{align*}\|t_\pm v+w\|&\leq 2(t_+-t_-)\|v\|\\
&\leq 2\sqrt{\kappa(v,w)^2-q(v)q(w)}.\frac{\|v\|}{q(v)}\end{align*}
Since
$$\|[v,w]\|=\|[t_\pm v+w,v]\|\leq \|t_\pm v+w\|\|v\|$$
this proves the lemma.

\end{proof}
If we define
$${\cal E}_\d^+=\{(x,y,z),\ z\geq(1+\d)\sqrt{x^2+y^2}\}$$
we have for any $v\in{\cal E}_\d^+$
$$\frac{\|v\|^2}{q(v)}\leq \frac{2}{\d},$$
so that
for any $v,w$ in ${\cal E}_\d^+$ 
$$\|[v,w]\|\leq C_\d\sqrt{\kappa(v,w)-q(v)q(w)},$$
with $C_\d=(2/\d)$. Notice that ${\cal E}_\d^+$ is a cone but it is not preserved by the $Ad$-action. However, for any $M\geq 0$ there exists a $\d'>0$ (depending on $M$) such that for any $v\in{\cal E}_\d^+$ and any $A\in SL(2,{\R})$ the norm of which is bounded by $M$, $\Ad(A).v$ lies in ${\cal E}_{\d'}^+$.

\begin{lemme}For any $v_1,\ldots,v_p$ in ${\cal E}_\d^+$
\begin{multline*}\sum_{1\leq i<j\leq p}\|[v_i,v_j]\|\leq (4/\d)\biggl(N(v_1+\cdots v_p)-(N(v_1)+\cdots+N(v_p))\biggr)^{1/2}\\.\biggl(\sum_{i=1}^p(N(v_i)+\|v_i\|\biggr)^{3/2}.\end{multline*}
\end{lemme}
\begin{proof}

We compute
\begin{align*}&N(\sum_{i=1}^pv_i)-\sum_{i=1}^pN(v_i)\\
&=\biggl(q(\sum_{i=1}^pv_i)\biggr)^{1/2}-\sum_{i=1}^p\sqrt{q(v_i)}\\
&=\biggl(\sum_{i=1}^pq(v_i)+2\sum_{1\leq i<j\leq p}\kappa(v_i,v_j)\biggr)^{1/2}-\sum_{i=1}^p\sqrt{q(v_i)}\\
&=\biggl(\biggl(\sum_{i=1}^p\sqrt{q(v_i)}\biggr)^2+2\sum_{1\leq i<j\leq p}\biggl(\kappa(v_i,v_j)-\sqrt{q(v_i)q(v_j)}\biggr)\biggr)^{1/2}-\sum_{i=1}^p\sqrt{q(v_i)}\\
&=\frac{2\sum_{1\leq i<j\leq p}\kappa(v_i,v_j)-\sqrt{q(v_i)q(v_j)}}{\sum_{i=1}^pN(v_i)+\biggl((\sum_{i=1}^pN(v_i))^2+2\sum_{1\leq i<j\leq p}\kappa(v_i,v_j)-\sqrt{q(v_i)q(v_j)}\biggr)^{1/2}}.\end{align*}
Since  $v_1,\ldots,v_p$ are in ${\cal E}_\d^+$
$$\|[v_i,v_j]\|\leq C_\d\sqrt{\kappa(v_i,v_j)^2-q(v_i)q(v_j)},$$
with $C_\d=(2/\d)$
and 
\begin{align*}&\sum_{1\leq i<j\leq p}\|[v_i,v_j]\|\\
&\leq C_\d\sum_{1\leq i<j\leq p}\sqrt{\kappa(v_i,v_j)^2-q(v_i)q(v_j)}\\
&\leq C_\d\sum_{1\leq i<j\leq p}\biggl(\sqrt{\kappa(v_i,v_j)-N(v_i)N(v_j)}.\sqrt{\kappa(v_i,v_j)+N(v_i)N(v_j)}\biggr)\\
&\leq C_\d\biggl(\sum_{1\leq i<j\leq p}\kappa(v_i,v_j)-N(v_i)N(v_j)\biggr)^{1/2}\biggl(\sum_{1\leq i<j\leq p}\kappa(v_i,v_j)+N(v_i)N(v_j)\biggr)^{1/2},\end{align*}
so that
\begin{multline*}\sum_{1\leq i<j\leq p}\|[v_i,v_j]\|\leq C_\d\biggl(N(v_1+\cdots+v_p)-(N(v_1)+\cdots+N(v_p)\biggr)^{1/2}\\
.\biggl(\sum_{i=1}^pN(v_i)+\biggl((\sum_{i=1}^pN(v_i))^2+2\sum_{1\leq i<j\leq p}\kappa(v_i,v_j)-N(v_i)N(v_j)\biggr)^{1/2}\biggr)^{1/2}\\
.\biggl(\sum_{1\leq i<j\leq p}\kappa(v_i,v_j)+N(v_i)N(v_j)\biggr)^{1/2} \end{multline*}

Observe that $0\leq \kappa(v_i,v_j)\leq \|v_i\|\|v_j\|$; thus
\begin{multline*}\sum_{1\leq i<j\leq p}\|[v_i,v_j]\|\leq C_\d\biggl(N(v_1+\cdots+v_p)-(N(v_1)+\cdots+N(v_p)\biggr)^{1/2}\\
.\biggl(\sum_{i=1}^pN(v_i)+\biggl((\sum_{i=1}^pN(v_i))^2+(\sum_{i=1}^p\|v_i\|)^2\biggr)^{1/2}\biggr)^{1/2}\\
.\biggl((\sum_{i=1}^p\|v_i\|)^2+(\sum_{i=1}^pN(v_i))^2\biggr)^{1/2} \end{multline*}
which gives the conclusion of the lemma
\end{proof}

If we use Cauchy-Schwarz inequality we then get the integrated version of the previous proposition:
\begin{prop}\lab{prop:4.1}If $\g_1(\cdot),\ldots,\g_p(\cdot)$ are $C^0$ maps from $I=[0,1]$ to ${\cal E}_\d^+\subset sl(2,{\R})$ we have
\begin{multline}\sum_{1\leq i<j\leq p}\int_I\|[\g_i(t),\g_j(t)]\|dt\leq \\(4/\d)\biggl(\int_IN(\g_1(t)+\cdots+\g_p(t))dt-\sum_{i=1}^p\int_IN(\g_i(t))dt\biggr)^{1/2}\\
.\biggl(\int_I(\sum_{i=1}^pN(\g_i(t))+\|\g_i(t)\|)^3dt\biggr)^{1/2}\end{multline}
\end{prop}

\paragraph*{A remark on a theorem by P. Thieullen}
At this point we can make the following comment. In the paper~\cite{Th} P. Thieullen gives a nice measurable description of general skew-products on $X\times SL(2,{\R})$ (with any dynamics on the base) via the notion of conformal barycenter (cf~\cite{D-E}). We intend to give an application of our techniques in this context. To simplify the notations we assume the dynamics on the base to be an irrational translation but in fact the result described in this paragraph are valid for any ergodic dyanamics on the base. 
We still assume that the fibered products of $(\a,A)$ are $C^0$-bounded. Consider the set ${\cal H}$ (resp. ${\cal H}_\d$) of $L^2$-maps from the circle ${\T}$ to the cone ${\cal E}^+\subset sl(2,{\R})$ (resp. ${\cal E}_\d^+$) and define the operator ${\cal U}$ defined on ${\cal H}$ by
$${\cal U}\s(\cdot)=\Ad(A(\cdot-\a))\s(\cdot-\a).$$
Given an element $\s$ of ${\cal H}$ we set
$$N(\s)=\int_{{\T}}N(\s(t))dt.$$
Of course for any $k\in{\Z}$
\begin{equation}\lab{5.1}N({\cal U}^k\s)=\s,\end{equation}
and since the fibered products of $A$ are assumed to be bounded
\begin{equation}\lab{5.2}\|{\cal U}^k\s\|_{L^2}\leq C.\|\s\|_{L^2},\end{equation}
uniformly in $k$. For the same reason, given any continuous $\s\in{\cal H}$ there exists a $\d>0$ such that for all $k$
\begin{equation}{\cal U}^k\s\in{\cal H}_\d.\lab{5.3}\end{equation}
For any set $I\subset {\Z}$ of $p$ consecutive integers we set $M(I)=\max I$ and we define
$$ {\cal L}_I\s=\sum_{k\in I}{\cal U}^k\s,\qquad \bar{\cal L}_I\s=\frac{1}{|I|}{\cal L}_I\s.$$
If $I,J$ are adjacent intervals of ${\Z}$ we clearly have
$${\cal L}_{I\cup J}={\cal L}_I+{\cal U}^{M(I)}{\cal L}_J.$$
Also,
$$\|{\cal U}\bar{\cal  L}_I-\bar{\cal L}_I\|_{L^1}\leq C|I|^{-1}.$$
Equality~(\ref{5.1}) shows that $N({\cal L}_I\s)$ depends only on  $\s$ and on the size $|I|$ of $I$. By the  anti-Minkowsky inequality 
the sequence $a_n(\s)=N({\cal L}_{\{1,\ldots, n\}}\s)$ is then super-additive and, since inequality~(\ref{5.2}) holds, the limit
$\bar a(\s)$ of the sequence $\bar a_n={1\ov n}a_n(\s)$
exists and is finite. Now take $\s\in{\cal H}$ any continuous function and take a sequence of positive integers $r_k$ increasing fast enough so that
$$\sum_{k=1}^{\infty}|\bar a(\s)-\bar a_{n_k}(\s)|^{1/2}<\infty,$$
where $n_k=r_k\cdots r_1$. We define inductively a sequence of intervals $I_k$ of length $n_k$ the following way: take for $I_1$ any interval of size $n_1$ and assume $I_l$ has been defined for $l<k$. Then apply proposition~\ref{prop:4.1} with $p=r_k$ and: 
$$\g_i={\cal L}_{I_{k-1}+in_{k-1}}={\cal U}^{M(I_{k-1})+in_{k-1}}{\cal L}_{I_{k-1}}\s;$$
we get:
$$\sum_{1\leq i<j\leq r_k}\int_{{\T}}\|[{\cal L}_{I_{k-1}+in_{k-1}}\s,{\cal L}_{I_{k-1}+jn_{k-1}}\s]\|dt\leq C_\d n_k^2\biggl(\bar a_{n_k}(\s)-\bar a_{n_{k-1}}(\s)\biggr)^{1/2}.$$
But this implies that for at least one index $1\leq i\leq r_k$
$$\sum_{1\leq j\leq r_k}\int_{{\T}}\|[{\cal L}_{I_{k-1}+in_{k-1}}\s,{\cal L}_{I_{k-1}+jn_{k_1}}\s]\|dt\leq C_\d\frac{n_k^2}{r_k}\biggl(\bar a_{n_k}(\s)-\bar a_{n_{k-1}}(\s)\biggr)^{1/2},$$
hence
$$\sum_{1\leq j\leq r_k}\int_{{\T}}\|[{\cal L}_{I_{k-1}}\s,{\cal U}^{(j-i)n_{k-1}+M(I_{k-1})}{\cal L}_{I_{k-1}}\s]\|dt\leq C_\d\frac{n_k^2}{r_k}\biggl(\bar a_{n_k}(\s)-\bar a_{n_{k-1}}(\s)\biggr)^{1/2},$$
and so if we define 
$$I_k=-in_{k-1}+(\bigcup_{j=1}^{r_k}(j+n_{k-1}I_{k-1})$$
we get
$$\int_{{\T}}\|[\bar{\cal L}_{I_{k-1}}\s,\bar{\cal L}_{I_{k}}\s]\|dt\leq C_\d \biggl(\bar a_{n_k}(\s)-\bar a_{n_{k-1}}(\s)\biggr)^{1/2}.$$

This implies that in the projective space the $L^1$ section so defined is invariant. 
A consequence of the previous comment is 
\begin{prop}If $(\a,A)$ has bounded fibered products then there exists a $B:{\T}\to SL(2,R)$ which is $L^\infty$-bounded and such that
$$B(\th+\a)^{-1}A(\th)B(\th)$$
are for a.e $\th$ rotation matrices.
\end{prop}

\section{Definition of $\g_k^{\pm},\eta_k^{\pm}$}
\subsection{Free monoids, free groups}
Let $\ti\Pi$ be the free monoid $M(D({\T}))$ constructed on the double $D({\T})={\T}\times\{-1,1\}$ of ${\T}$, that is, the set of words of finite length of the form $t_1^{\e_1}\cdots t_p^{\e_p}$, $\e_i\in\{-1,1\}$, $t_i\in{\T}$ ($i=1,\ldots,p$). The composition law is then concatenation:
$$(t_1^{\e_1}\cdots t_p^{\e_p})*(t_{p+1}^{\e_{p+1}}\cdots t_s^{\e_s})=t_1^{\e_1}\cdots t_s^{\e_s},$$
which is an associative composition law; the identity element is the empty word we shall denote $e$. The set $\{t_1,\ldots,t_p\}$ is called the support of $w$ and we denote it by $|w|$. Notice that on $D({\T})$ one can define an involution: $\t(t^{\e})=t^{-\e}$ and that $\t$ can be extended to $M(D({\T}))$ by
$\t(t_1^{\e_1}\cdots t_p^{\e_p})=t_p^{-\e_p}\cdots t_1^{-\e_1}$.
We  introduce the notion of {\it elementary reduction}: given a word $w=t_1^{\e_1}\cdots t_k^{\e_k}t_{k+1}^{\e_{k+1}}\cdots t_n^{\e_n}$  we define $\r(w)=t_1^{\e_1}\cdots t_{k-1}^{\e_{k-1}}t_{k+2}^{\e_{k+2}}\cdots t_n^{\e_n}$ if $k$ is the largest integer such that  $t_k=t_{k+1}$,$\e_k=-\e_{k+1}$. The sequence $\r^n(w)$ clearly stabilizes and we define the {\it reduced} word $r(w)=\lim_{n\to\infty}\r^n(w)$. The free group $\Pi$ constructed on ${\T}$ is then the set of words $w\in M(D({\T}))$ such that $r(w)=w$. The composition law $w_1.w_2=r(w_1*w_2)$ makes $(\Pi,.)$ a (free) group (the main point here is to check the associativity).
Also if $w$ is a reduced word we say that it is {\it square free} if $w=t_1^{\e_1}\cdots t_p^{\e_p}$ with $t_i\ne t_j$ for every $1\leq i<j\leq p$. We then say that two reduced words $w_1,w_2$ are {\it mutually prime } if $|w_1|\cap |w_2|=\emptyset$.

Let ${\Z}[Y]$ the set of polynomials in the variable $Y$  let
$(X_t^{\pm})_{t\in{\T}}$ be a set of formal variables and denote by ${\cal \ti T}_Y$ (resp. ${\cal T}_Y$) the set of finite sums with integer coefficients
$$n_1(Y)w_1.X_{t_1}^{\e_1}+\cdots+n_p(Y)w_p.X_{t_p}^{\e_p},$$
with $\e_i\in\{-1,1\}$, $w_i\in{\ti\Pi}$ (resp. $w_i\in\Pi$), $n_i(Y)\in{\Z}[Y]$, $t_i\in{\T}$ (this is the set of  almost zero maps from ${\ti \Pi}\times{\T}\times\{-1,1\}$ to ${\Z}[Y]$). We can define a map $R:{\cal\ti T}_Y\to{\cal T}_Y$ by
\begin{multline*}R(n_1(Y)w_1.X_{t_1}^{\e_1}+\cdots+n_p(Y)w_p.X_{t_p}^{\e_p})=\\n_1(Y)r(w_1).X_{t_1}^{\e_1}+\cdots+n_p(Y)r(w_p).X_{t_p}^{\e_p}.\end{multline*}
The sets ${\cal\ti T}_Y$ and ${\cal T}_Y$ can be given the structure of a ${\Z}[Y]$-module.
We denote similarly ${\cal\ti T}_Y^+$ the set of such sums with coefficients $n_i(Y)$ being polynomials in ${\Z}[Y]$ with nonnegative integer coefficients (this is only a monoid). 
There is a natural left action of $\ti\Pi$ on ${\cal\ti T}_Y$  defined by
$$w*(\sum_{i=1}^pn_i(Y)w_i.X_{t_i}^{\e_i})=\sum_{i=1}^p n_i(Y)(w*w_i).X_{t_i}^{\e_i}.$$ Observe that for $T\in{\cal\ti T}$, $w\in\ti\Pi$ the following relation holds  $R(w*T)=r(w).R(T).$
Also an action of $\Pi$ on ${\cal T}_Y$ is defined by:
$$w.(\sum_{i=1}^pn_i(Y)w_i.X_{t_i}^{\e_i})=\sum_{i=1}^p n_i(Y)(w.w_i).X_{t_i}^{\e_i}.$$
Given $T\in{\cal\ti T}_Y$ of the form $\sum_{i=1}^pn_i(Y)w_i.X_{t_i}^{\e_i}$ we introduce its positive and negative parts $T^+$, $T^-$
$$T^\pm=\sum_{i=1}^pn_i^\pm(Y)w_i.X_{t_i}^{\e_i}$$
where $n_i(Y)$ is decomposed uniquely in $n_i(Y)=n_i^+(Y)+Yn_i^-(Y)$ where $n_i^\pm(Y)$ are polynomial wich monomials in $Y$ are of even degrees.

Similarly we define ${\cal\ti T}$, ${\cal T}$, ${\cal T}^+$ where the $n_i(Y)$ are replaced by $n_i\in{\Z}$, $n_i\in{\Z}_+$.
 
Assume now we are given elements $(T_t)_{t\in{\T}}$ of ${\cal \ti T}_Y$ and define the map
$\ti F_Y:\ti\Pi\to {\cal \ti T}_Y$ by
\begin{align}\lab{def:F1} &\ti F_Y(t)=T_t\\
&\ti F_Y(t^{-1})=Yt^{-1}.T_t\lab{def:F2}\\
&\lab{def:F3}\ti F_Y(t_1^{\e_1}\cdots t_r^{\e_r})=\ti F_Y(t_1^{\e_1})+t_1^{\e_1}*\ti F_Y(t_2^{\e_2})+\cdots+(t_1^{\e_1}\cdots t_{r-1}^{\e_{r-1}})*\ti F_Y(t_r^{\e_r}).\end{align} 
for all $t\in{\T}$, $t_1,\ldots,t_r\in{\T}$, $\e_i\in\{-1,1\}$. For any $w\in\ti\Pi$ one can decompose the polynomial $\ti F_Y(w)$ uniquely in
$$\ti F_Y(w)={\ti F}_Y^+(w)+Y{\ti F}^-_Y(w).$$
We then have
\begin{lemme}\lab{lemme:a1}For all $w_1,w_2\in\ti\Pi$
\begin{align*}&i)\quad \ti F_Y(w_1*w_2)=\ti F_Y(w_1)+w_1*\ti F_Y(w_2)\\
&ii)\quad {\ti F}_Y^{\pm}(w_1*w_2)={\ti F}_Y^{\pm}(w_1)+w_1*{\ti F}_Y^{\pm}(w_2). 
\end{align*}
Moreover, if for each $t\in{\T}$, $T_t\in{\cal\ti T}_Y^+$ then ${\ti F}_Y^\pm$ are in ${\cal\ti T}_Y^+$.
\end{lemme}
\begin{proof}

The first equality is proven by staightforward calculation and the second one by identifying the odd and even parts of the polynomials. The last claim is obvious.
\end{proof}
We define for $w\in{\Pi}$
$$F_Y(w)=R\circ \ti F_Y(w),\qquad F_Y^{\pm}(w)=R\circ {\ti F}_Y^{\pm}(w).$$
We have
\begin{lemme}\lab{lemme:fpm'}If $w_1,w_2\in\Pi$ 
\begin{align*}&F_Y(w_1.w_2)\equiv F_Y(w_1)+w_1.F_Y(w_2)\mod (1+Y)\\
&F_Y^{\pm}(w^{-1})\equiv w^{-1}.F_Y^{\mp}(w))\mod (1+Y)
\end{align*}
and if $w_1,w_2$ are mutually prime
$$F^\pm_Y(w_1.w_2)=F^\pm_Y(w_1)+w_1.F^\pm_Y(w_2)$$
\end{lemme} 
\begin{proof}

For the proof of the first equality just observe that for $x\in D({\T})$ ($x=t^{\pm 1}$), $R(\ti F_Y(x)+x*\ti F_Y(\t(x)))$ that is $F_Y(x)+x.F_Y(\t(x))$ is divisible by $(1+Y)$. This implies that $F_Y(\r(w))\equiv F_Y(w)\mod(1+Y)$ and hence $F_Y(r(w))\equiv F_Y(w)\mod (1+Y)$.

For the second equality observe that this is true for words of length 1 since
$${\ti F}_Y^+(t^{-1})=Y^2t^{-1}*{\ti F}_Y^-(t),\quad {\ti F}_Y(t^{-1})=t^{-1}{\ti F}_Y^+(t)$$
and 
proceed by induction on the length of $w$ (which is reduced): assume $w=w_1w_2$ with lengths of $w_1$ and $w_2$ smaller than the length of $w$ then
\begin{align*}F_Y^\pm((w_1w_2)^{-1})&=R({\ti F}_Y^{\pm}(w_2^{-1})+w_2^{-1}{\ti F}_Y(w_1^{-1}))\\
&=R((w_1w_2)^{-1}.((w_1w_2).{\ti F}_Y^\pm(w_2^{-1})+w_1.{\ti F}_Y^\pm(w_1^{-1})))\\
&=(w_1w_2)^{-1}.R((w_1w_2).{\ti F}_Y^\pm(w_2^{-1})+w_1.{\ti F}_Y^\pm(w_1^{-1}))\\
&\equiv w^{-1}.R((w_1w_2).w_2^{-1}{\ti F}^\pm_Y(w_2)+w_1.w_1^{-1}{\ti F}^\pm_Y(w_1))\mod (1+Y)\\
&\equiv w^{-1}.F^\pm_Y(w)\mod (1+Y)
\end{align*}

The third one is a consequence of equality ii) of lemma~\ref{lemme:a1} and the fact that $r(w_1.w_2)=w_1.w_2$ if $w_1,w_2$ are mutually prime. 
\end{proof}

As a particular case we choose for $(T_t)_{t\in{\T}}$
$$\forall t\in{\T},\ T_t:=X_t^++YX_t^{-},$$
which are in ${\cal\ti T}_Y^+$.
We then define  for $w\in\Pi$ the elements $F(w),F^\pm(w)\in {\cal T}$ by making $Y=-1$:
$$F^\pm(w)=F_Y^\pm(w)_{|Y=-1},\qquad F(w)=F_Y(w)_{|Y=-1},$$
and we get
$$F(w)=F^+(w)-F^-(w).$$
\begin{lemme}\lab{lemme:fpm}If $w_1,w_2\in\Pi$ 
\begin{align*}&F(w_1.w_2)= F(w_1)+w_1.F(w_2)\\
&F^{\pm}(w^{-1})= w^{-1}.F^{\mp}(w))
\end{align*}
and if $w_1,w_2$ are mutually prime
$$F^\pm(w_1.w_2)=F^\pm(w_1)+w_1.F^\pm(w_2).$$
Also $F^\pm(w)\in{\cal T}^+$.
\end{lemme} 

\subsection{Applications}
Our aim now is to define $\g_k^{\pm},\eta_k^{\pm}$. We assume that we are given a smooth $A:{\T}\to SL(2,{\R})$ and $\a$ {\it irrational}. We denote by $\eta_0(\cdot)$ 
$$\eta_0(t)=L(A)(t):=(\pa_t A(t))A(t)^{-1},$$
and we {\it choose} a decomposition of $\eta_0(t)$ of the form
\begin{equation}\lab{etapm}\eta_0(t)=\eta_0^+(t)-\eta^-_0(t),\end{equation}
where $\eta_0^+,\eta_0^-:{\T}\to {\cal E}^+$ are smooth. Such a decomposition always exists: it is enough to take for example $\eta_0^+$ constant and very large so that $\eta_0^+-\eta(t)$ lies in ${\cal E}_\d$ ($\d>0$) for any $t$.

Denote by $\pi$ the map
\begin{align*}\Pi&\to SL(2,{\R})\\
(t_1^{n_1}\cdots t_r^{n_r})&\mapsto A(t_1)^{n_1}\cdots A(t_r)^{n_r}\end{align*}
This is a morphism of group. Also for each $t\in{\T}$ and $n\geq 1$ we define the words in $\Pi$ 
\begin{align*}&w(n,t)=(t+(n-1)\a)\cdots (t)\\
&w(0,t)=e\\
&w(-n,t)=(t-n\a)^{-1}\cdots(t-\a)^{-1}\end{align*}
and
\begin{align*}a_k(t)&=w((-1)^{k}q_k,t)\\
c_k(t)&=w((-1)^{{k-1}}q_{k-1},t)\end{align*} 

We have
$$\pi(a_k(t))=A^{(k)}(t),\qquad \pi(c_k(t))=C^{(k)}(t),$$
and also
\begin{align}a_k(t)&=c_{k-1}(t-a_k\b_{k-1})a_{k-1}(t-a_k\b_{k-1})^{-1}\cdots a_{k-1}(t-\b_{k-1})^{-1}\lab{a}\\
c_k(t)&=a_{k-1}(t)\lab{c}\end{align}

Since $\a$ is irrational, each word $w(n,t)$ is square free and so are $a_k(t),c_k(t)$. This proves that in the previous inequality~(\ref{a}) the words are mutually prime (since the length of $a_k(t)$ is $q_k=a_kq_{k-1}+q_{k-2}$ which is the sum of the lengths in the RHS of equality~(\ref{a}).)

As we have seen if $w=t_1^{\e_1}\cdots t_p^{\e_p}$ is squarefree, $F(w)$ (resp. $F^{\pm}(w)$) is a sum of the form
$$\sum_{i=1}^rw_{i}.T_{t_i}$$
(and in fact $w_1=e$ and $w_{i+1}=t_1^{\e_1}\cdots t_i^{\e_i}$ for $p>i\geq 1$). We define the element $f(w)$ of $sl(2,{\R})$ by
$$f(w)=\sum_{i=1}^p\Ad(\pi(w_{i})).\eta_0(t_i).$$
Similarly  $F^\pm(w)$ is a sum of the form
$$F^\pm(w)=\sum_i m_i^\pm w_{i}.X_{t_i}^{\nu_i},$$
where $m_i^\pm$ are nonnegative integers and $\nu_i\in\{-1,1\}$.
We define the element  $f^\pm(w)$ by
$$f^\pm(w)=\sum_{i}m_i^\pm\Ad(\pi(w_{i})).\eta_0^{\nu_i}(t_i),$$
where $\eta_0^{\nu_i}$ is defined by~(\ref{etapm}).
Since $T_t=X^+_t-X_t^-$ and $\eta_0=\eta_0^+-\eta_0^-$ we have
$$f(w)=f^+(w)-f^-(w).$$
It is clear that lemma~\ref{lemme:fpm} remains  true with $F$ replaced with $f$ and $F^\pm$ replaced with $f^\pm$.

It is also clear that for every $t$
$$f(a_k(t))=\eta_k(t),\qquad f(c_k(t))=\g_k(t),$$
where
$$\eta_k(t)=L(A^{(k)})(t),\qquad \g_k(t)=L(C^{(k)})(t),$$ 
(more generally $L(\pi(w(t)))=f(w(t))$).
We can now define $\eta_k^\pm(t)$ and $\g_k^\pm(t)$ by
$$\eta_k^\pm(t)=f^\pm(a_k(t)),\qquad \g_k^\pm(t)=f^\pm(c_k(t)),$$
and we have
$$\eta_k(t)=\eta_k^+(t)-\eta_k^-(t),\qquad \g_k(t)=\g_k^+(t)-\g_k^-(t).$$
Notice that $\g_k^\pm(\cdot),\eta_k^\pm(\cdot)$ are smooth from ${\T}$ to the cone ${\cal E}^+$.
Since lemma~\ref{lemme:fpm} is true with $F^\pm$ replaced with $f^\pm$ and since~(\ref{a}), (\ref{c}) hold and that the factors are mutually prime,  we have
\begin{equation}\g_k^{\pm}(t)=\eta_{k-1}^\pm(t)\lab{200}\end{equation}
\begin{multline}\eta_k^\pm(t)=\Ad\biggl(A^{(k-1)}(t+\b_{k-2}-a_k\b_{k-1})^{-1}\biggr).\eta_{k-1}^\mp(t+\b_{k-2}-a_k\b_{k-1})+\cdots\\
+\Ad\biggl(A^{(k-1)}(t+\b_{k-2}-a_k\b_{k-1})^{-1}\cdots \\
\ \ \ \ \ A^{(k-1)}(t+\b_{k-2}-\b_{k-1})^{-1}\biggr).\eta_{k-1}^\mp(t+\b_{k-2}-\b_{k-1})\\
+\Ad\biggl(A^{(k-1)}(t+\b_{k-2}-a_k\b_{k-1})^{-1}\cdots A^{(k-1)}(t+\b_{k-2}-\b_{k-1})^{-1}\biggr).\g_{k-1}^\pm(t)\lab{201}\end{multline}

Assume that $2|(k-l)$; then using obvious notations and lemma~\ref{m_1...}
\begin{align*}\begin{pmatrix}(\b_{k-1},c_k)\\(\b_k,a_k)\end{pmatrix}&=\begin{pmatrix}0&1\\1&-a_k\end{pmatrix}\cdots\begin{pmatrix}0&1\\1&-a_{l+1}\end{pmatrix}\begin{pmatrix}(\b_{l-1},c_l)\\(\b_l,a_l)\end{pmatrix}\\
&=\begin{pmatrix}m_1&m_2\\m_3&m_4\end{pmatrix}\begin{pmatrix}(\b_{l-1},c_l)\\(\b_l,a_l)\end{pmatrix}\end{align*}
that is
\begin{align*}&c_k(t)=c_l(t+(m_1-1)\b_{l-1}-|m_2|\b_l)\cdots c_l(t-|m_2|\b_l)\\&\ \ \ \ \ \ \ \ \ \ \ \ \ \ a_l(t-|m_2|\b_l)^{-1}\cdots a_l(t-\b_l)^{-1}\\
&a_k(t)=c_l(t+(m_4-1)\b_l-|m_3|\b_{l-1})^{-1}\cdots c_l(t+(m_4-1)\b_l-\b_{l-1})^{-1}\\
&\ \ \ \ \ \ \ \ \ \ \ \ \ \  a(t+(m_4-1)\b_l)\cdots a(t)\end{align*}
that we shall write
\begin{align*}&c_k(t)=c_l(s_1)\cdots c_l(s_{m_1})a_l(s_{m_1+1})^{-1}\cdots a_l(s_{m_1+|m_2|})^{-1}\\
&a_k(t)=c_l(t_1)^{-1}\cdots c_l(t_{|m_3|})^{-1}a(t_{|m_3|+1})\cdots a(t_{|m_3|+m_4})\end{align*}
where $s_i=t-|m_2|\b_l+\b_{l-1}(m_1-i)$ if $1\leq i\leq m_1$ and $s_i=t-(m_1+|m_2|-i+1)\b_l$ if $m_1+1\leq i\leq m_1+|m_2|$ and  
$t_i=t+(m_4-1)\b_l-(|m_3|-i+1)\b_{l-1}$ if $1\leq i\leq |m_3|$ and $t_i=t+(|m_3|+m_4-i)\b_l$ if $|m_3|\leq i\leq |m_3|+m_4$.
Consequently
\begin{multline}\lab{gkl}\g_k^\pm(t)=\Ad(U_1).\g_l^\pm(s_1)+\cdots+\Ad(U_1\cdots U_{m_1}).\g_l^\pm(s_{m_1})+\\
\Ad(U_1\cdots U_{m_1+1}).\eta_l^\mp(s_{m_1+1})+\cdots+\Ad(U_1\cdots U_{m_1+|m_2|}).\eta_l^\mp(s_{m_1+|m_2|}),\end{multline}
\begin{multline}\lab{ekl}\eta_k^\pm(t)=\Ad(V_1).\g_l^\mp(t_1)+\cdots+\Ad(V_1\cdots V_{m_3}).\g_l^\mp(t_{|m_3|})+\\
\Ad(V_1\cdots V_{|m_3|+1}).\eta_l^\pm(t_{|m_3|+1})+\cdots+\Ad(V_1\cdots V_{|m_3|+m_4}).\eta_l^\pm(t_{|m_3|+m_4}),\end{multline}
where on the one hand $U_1=Id$, $U_i=C^{(l)}(s_{i-1})$ if $2\leq i\leq m_1$, $U_{m_1+1}=C^{(l)}(s_{m_1})A^{(l)}(s_{m_1+1})^{-1}$ and   $U_i=A^{(l)}(s_i)^{-1}$ if $m_1+2\leq i\leq m_1+|m_2|$  and on the other hand $V_i=C^{(l)}(t_i)^{-1}$ if $1\leq i\leq |m_3|$, $V_{m_3+1}=Id$ and $V_i=A^{(l)}(t_{i-1})$ if $|m_3|+2\leq i\leq |m_3|+m_4$.
Also,
\begin{multline}\lab{gkl1}
(\g_k^+ \pm\g_k^-)(t)=
\Ad(U_1).(\g_l^+ \pm\g_l^-)(s_1)+\cdots+\Ad(U_1\cdots U_{m_1}).(\g_l^+ \pm\g_l^-)(s_{m_1})+\\
\Ad(U_1\cdots U_{m_1+1}).(\eta_l^- \pm\eta_l^+)(s_{m_1+1})+\\
\cdots+\Ad(U_1\cdots U_{m_1+|m_2|}).(\eta_l^- \pm\eta_l^+)(s_{m_1+|m_2|}),
\end{multline}
\begin{multline}\lab{ekl1}
(\eta_k^+ \pm\eta_k^-)(t)=
\Ad(V_1).(\g_l^-\pm\g_l^+)(t_1)+\cdots+\Ad(V_1\cdots V_{m_3}).(\g_l^- \pm\g_l^+)(t_{|m_3|})+\\
\Ad(V_1\cdots V_{m_3+1}).(\eta_l^+ \pm\eta_l^-)(t_{|m_3|+1})+\\
\cdots+\Ad(V_1\cdots V_{m_3+m_4}).(\eta_l^+ \pm\eta_l^-)(t_{|m_3|+m_4}),
\end{multline}

Formulas~(\ref{gkl1}),~(\ref{gkl1}) where the $\pm$ sign  is the $-$ sign, just express $\g_k,\eta_k$ in terms of $\g_l,\eta_l$. 

\subsection{Estimates for $\eta_k,\g_k,\eta_k^\pm, \g_k^\pm$}
\begin{lemme}\lab{u}There exists $\d'>0$ depending only on $\d$ and on $\sup_{k\in{\Z}}\|A_k(\cdot)\|_0$ such that for any $k\in{\N}$, $t\in{\T}$
$$\eta_k^\pm(t)\in{\cal E}_{\d'}^+,\qquad \g_k^\pm(t)\in{\cal E}_{\d'}^+.$$
\end{lemme}
\begin{proof}

From the definition of $\eta_k^\pm(t)$ it is clear that $\eta_k^\pm$ is a sum of terms of the form $\Ad(U_j).\eta_0^\pm(t_j)$ where $U_j$ is a fibered product of $A(\cdot)$ over the rotation $x\mapsto x+\a$. Since these products are by assumption bounded uniformly and since $\eta_0^\pm(t_j)$ is in ${\cal E}_\d^+$ the remark made just after the definition of ${\cal E}_\d$  shows that there exists $\d'>0$ such that  each $\Ad(U_j).\eta_0^\pm(t_j)$ is in ${\cal E}_{\d'}^+$. But this set is convex and the sum has only positive coefficients: the proof of the lemma is complete 
\end{proof}

Similarly
\begin{lemme}\lab{u:1}For any $k\geq 0$
$$\b_{k-1}\max(\|\g_k^\pm\|_0,\|\eta_k^\pm\|_0,\|\g_k\|_0,\|\eta_k\|_0)\leq M$$
where $M$ is a constant depending only on  $\sup_{k\in{\Z}}\|A_k(\cdot)\|_0$ (and on the chosen decomposition  $\eta_0=\eta_0^+-\eta_0^-$)
\end{lemme}
\begin{proof}

Just estimate the (euclidean) norm $\|\cdot\|$ in formulas~(\ref{gkl}),~(\ref{ekl}), (\ref{ekl1}), (\ref{gkl1}) where $l=0$ and use lemma~\ref{m_1...}.
\end{proof}
In a similar way
\begin{lemme}\lab{u:2}For any $k\geq 0$
$$\b_{k-1}^2\max(\|\pa\g_k\|_0,\|\pa\eta_k\|_0)\leq M'$$
where $M'$ is a constant depending only on  $\sup_{k\in{\Z}}\|A_k(\cdot)\|_0$ and on the $C^0$-norm of $\pa(\eta_0^\pm)$ ($\eta_0=\eta_0^+-\eta_0^-$)
Similarly,
$$\b_{k-1}^2\max(\|\pa\g_k^\pm\|_0,\|\pa\eta_k^\pm\|_0)\leq M'.$$
More generally for any integer $r$ there exists a constant $M_r$ (depending on $r$,  $\sup_{k\in{\Z}}\|A_k(\cdot)\|_0$ and on $\|\eta_0^\pm\|_{r+1}$)
such that
$$\b_{k-1}^{r+1}\max(\|\pa\g_k^\pm\|_r,\|\pa\eta_k^\pm\|_r)\leq M_r.$$
\end{lemme}
\begin{proof}

Apply proposition~\ref{prop:1.3.4}  to formulas~(\ref{gkl1}),~(\ref{ekl1}) where one makes $l=0$ and use lemma~\ref{m_1...}. The proof of the second statement,  as well as the one concerning higher order derivatives, is done the same way (use proposition~\ref{prop:1.3.4}).
\end{proof}

\section{Further estimates}
\subsection{Some integrated quantities}

We now introduce some important quantities:
\begin{align}
&e_k^\pm=\int_{\T}N(\eta_k^\pm(s))ds, & e_k=\int_{\T}N(\eta_k^+(s)+\eta_k^-(s))ds\\
&f_k^\pm=\int_{\T}N(\g_k^\pm(s))ds, & f_k=\int_{\T}N(\g_k^+(s)+\g_k^-(s))ds\\
&u_k^{\pm}=e_k^{\pm}+\a_k f_k^{\mp}, &  u_k=e_k+\a_k f_k.
\end{align}
Also, ${\bar u}^\pm_k=\b_{k-1}u_k^\pm$, ${\bar u}_k=\b_{k-1}u_k$.

If we define the new variables
\begin{equation}\lab{tildeC...}{\ti C}^{(k)}(t)=C^{(k)}(\b_{k-1}t),\qquad {\ti A}^{(k)}(t)=A^{(k)}(\b_{k-1}t)\end{equation}
\begin{equation}\lab{tigamma...}\ti\g_k(t)=L({\ti C}^{(k)})(t),\qquad \ti\eta_k(t)=L({\ti A}^{(k)})(t)\end{equation}
we have
$$\ti \g_k(t)=\b_{k-1}\g_k(\b_{k-1}t)\qquad \ti \eta_k(t)=\b_{k-1}\eta_k(\b_{k-1}t)$$
and we define
$${\ti \g}^\pm_k(t)=\b_{k-1}\g_k^\pm(\b_{k-1}t),\qquad {\ti\eta}^\pm_k(t)=\b_{k-1}\eta_k^\pm(\b_{k-1}t)$$
and
$${\ti u}_k^\pm=\int_0^1N({\ti \eta}_k^\pm(s))ds+\int_0^{\a_k}N({\ti\g}_k^\mp(s))ds.$$
\begin{lemme}\lab{u:ineq2}For $k\geq 1$,
\begin{align*}&2M\geq{\bar u}^\pm_k\geq {\bar u}^\mp_{k-1}\geq 0\\
&2M\geq{\bar u}_k\geq {\bar u}_{k-1}\geq 0\\
&2M\geq{\ti u}^\pm_k\geq {\ti u}^\mp_{k-1}\geq 0.\end{align*}
\end{lemme}
\begin{proof}

First notice that the estimations from above by $M$ come from lemma~\ref{u:1} since $N(v)\leq\|v\|$ if $v\in {\cal E}^+$.

Let us prove the first inequality: using formulas~(\ref{gkl})-(\ref{ekl1}) and the fact that $N$ satisfies the anti-Minkowsky inequality and is preserved by the $\Ad$-action we get
$$u_k^\pm\geq (a_k+\a_k)\int_0^1N(\eta_{k-1}^\mp(s))ds+\int_0^1N(\g_{k-1}^\pm(s))ds$$
and since $a_k+\a_k=\a_{k-1}^{-1}$ this proves the first inequality.

The proof of the second one is similar.

For the third inequality:
\begin{multline*}N(\eta_k^\pm(t))\geq N(\eta_{k-1}^\mp(t+\b_{k-2}-a_k\b_{k-1})+\cdots\\+ N(\eta_{k-1}^\mp(t+\b_{k-2}-\b_{k-1})+N(\g_{k-1}^\pm(t)),\end{multline*}
and consequently
\begin{multline*}N({\ti \eta}_k^\pm(t))\geq \a_{k-1}\biggl(N({\ti \eta}_{k-1}^\mp(\a_{k-1}(t-a_k)+1)+\cdots\\+ N({\ti \eta}_{k-1}^\mp(\a_{k-1}(t-1)+1)+N({\ti \g}_{k-1}^\pm(\a_{k-1}t))\biggr);\end{multline*}
thus
\begin{align*}&\int_0^{\a_k}N({\ti \g}_k^\mp(t))dt+\int_0^1N({\ti\eta}_k^\pm(t))dt\geq \\
&\int_0^{\a_k\a_{k-1}}N({\ti\eta}_{k-1}^\mp(s))ds+\int_0^{\a_{k-1}}N({\ti\g}_{k-1}^\pm(s))ds+\\
&\int_{1-a_k\a_{k-1}}^{1-(a_k-1)\a_{k-1}}N({\ti\eta}_{k-1}^\mp(s))ds+
\cdots+\int_{1-\a_{k-1}}^1N({\ti \eta}_{k-1}^\mp(s))ds\\
&\geq \int_0^{\a_{k-1}}N({\ti\g}_{k-1}^\pm(s))ds+\int_0^1N({\ti\eta}_{k-1}^\mp(s))ds\end{align*}
since $1-a_k\a_{k-1}=\a_k\a_{k-1}$. The lemma is proved.

\end{proof}

As a corollary of the preceding lemma we get
\begin{cor} The following limits exist and satisfy:
\begin{align*}&\lim_{k\to\infty}{\bar u}^\pm_{2k}=\lim_{k\to\infty}{\bar u}_{2k-1}^\mp\\
&\lim_{k\to\infty}{\bar u}_{2k}=\lim_{k\to\infty}{\bar u}_{2k-1}\\
&\lim_{k\to\infty}{\ti u}^\pm_{k}=\lim_{k\to\infty}{\ti u}_{k-1}^\mp.\end{align*}
\end{cor}

\subsection{Estimates for the $L^1$ norm of $\g_k,\eta_k$}

Assume $2|(k-l)$. Remembering  formulae~(\ref{gkl}), (\ref{ekl}) let us  introduce the nonnegative numbers
$$w^{(\g)}_{\g^\pm,\g^\pm}(t)=\sum_{1\leq i<j\leq |m_1|}\biggl\|\biggl[\Ad(U_1\cdots U_{i})\g_l^\pm(s_i),\Ad(U_1\cdots U_{j})\g_l^\pm(s_j)\biggr]\biggr\|$$
$$w^{(\eta)}_{\g^\pm,\g^\pm}(t)=\sum_{1\leq i<j\leq |m_3|}\biggl\|\biggl[\Ad(V_1\cdots V_{i})\g_l^\pm(t_i),\Ad(V_1\cdots V_{j})\g_l^\pm(t_j)\biggr]\biggr\|$$
$$w^{(\g)}_{\eta^\pm,\eta^\pm}(t)=\sum_{m_1+1\leq i<j\leq m_1+|m_2|}\biggl\|\biggl[\Ad(U_1\cdots U_{i})\eta_l^\pm(s_i),\Ad(U_1\cdots U_{j})\eta_l^\pm(s_j)\biggr]\biggr\|$$
$$w^{(\eta)}_{\eta^\pm,\eta^\pm}(t)=\sum_{|m_3|+1\leq i<j\leq |m_3|+m_4}\biggl\|\biggl[\Ad(V_1\cdots V_{i})\eta_l^\pm(t_i),\Ad(V_1\cdots V_{j})\eta_l^\pm(t_j)\biggr]\biggr\|$$
$$w^{(\g)}_{\g^\mp,\eta^\pm}(t)=\sum_{1\leq i\leq m_1< j\leq m_1+|m_2|}\biggl\|\biggl[\Ad(U_1\cdots U_{i})\g_l^\mp(s_i),\Ad(U_1\cdots U_{j})\eta_l^\pm(s_j)\biggr]\biggr\|$$
$$w^{(\eta)}_{\g^\mp,\eta^\pm}(t)=\sum_{1\leq i\leq |m_3|< j\leq |m_3|+m_4}\biggl\|\biggl[\Ad(V_1\cdots V_{i})\g_l^\mp(t_i),\Ad(V_1\cdots V_{j})\eta_l^\pm(t_j)\biggr]\biggr\|$$
and similarly  referring to formulae~(\ref{gkl1}), (\ref{ekl1})
\begin{multline*}w_{\g,\g}^{(\g)}(\pm)(t)=\sum_{1\leq i<j\leq m_1}\biggl\|\biggl[\Ad(U_1\cdots U_{i}).(\g_l^+(s_i)\pm\g_l^-(s_i)),\\\Ad(U_1\cdots U_{j}).(\g_l^+(s_j)\pm\g_l^-(t_j))\biggr]\biggr\|\end{multline*}
\begin{multline*}w_{\g,\g}^{(\eta)}(\pm)(t)=\sum_{1\leq i<j\leq |m_3|}\biggl\|\biggl[\Ad(V_1\cdots V_{i}).(\g_l^+(t_i)\pm\g_l^-(t_i)),\\\Ad(V_1\cdots V_{j}).(\g_l^+(t_j)\pm\g_l^-(t_j))\biggr]\biggr\|\end{multline*}
\begin{multline*}w_{\eta,\eta}^{(\g)}(\pm)(t)=\sum_{m_1+1\leq i<j\leq m_1+|m_2|}\biggl\|\biggl[\Ad(U_1\cdots U_{i}).(\eta_l^+(s_i)\pm\eta_l^-(s_i)),\\\Ad(U_1\cdots U_{j}).(\eta_l^+(s_j)\pm\eta_l^-(t_j))\biggr]\biggr\|\end{multline*}
\begin{multline*}w_{\eta,\eta}^{(\eta)}(\pm)(t)=\sum_{|m_3|+1\leq i<j\leq |m_3|+m_4}\biggl\|\biggl[\Ad(V_1\cdots V_{i}).(\eta_l^+(t_i)\pm\eta_l^-(t_i)),\\\Ad(V_1\cdots V_{j}).(\eta_l^+(t_j)\pm\eta_l^-(t_j))\biggr]\biggr\|\end{multline*}
\begin{multline*}w_{\g,\eta}^{(\g)}(\pm)(t)=\sum_{1\leq i\leq m_1<j\leq m_1+|m_2|}\biggl\|\biggl[\Ad(U_1\cdots U_{i}).(\g_l^+(s_i)\pm\g_l^-(s_i)),\\\Ad(U_1\cdots U_{j}).(\eta_l^+(s_j)\pm\eta_l^-(t_j))\biggr]\biggr\|\end{multline*}
\begin{multline*}w_{\g,\eta}^{(\eta)}(\pm)(t)=\sum_{1\leq i\leq |m_3|<j\leq |m_3|+m_4}\biggl\|\biggl[\Ad(V_1\cdots V_{i}).(\g_l^+(t_i)\pm\g_l^-(t_i)),\\\Ad(V_1\cdots V_{j}).(\eta_l^+(t_j)\pm\eta_l^-(t_j))\biggr]\biggr\|\end{multline*}
and similarly the elements $W_{\g^\pm,\g^\pm}(t),\ldots, W_{\g,\eta}^{(\eta)}(\pm)(t)$ of ${\cal E}^+$ defined by the same equalities except that the symbol $\|\cdot\|$ is removed.
We also define the quantities $w_{\g^\pm,\g^\pm},\ldots,W_{\g,\eta}^{(\eta)}(\pm)$ obtained by integrating on ${\T}$ with respect to $dt$ the corresponding functions of $t$. 

\begin{lemme}\lab{c1}
\begin{multline*}\b_{k-1}^2\max(\|\pa \g_k\|_{L^1({\T})},\|\pa \eta_k\|_{L^1({\T})})\leq \\2\frac{\b_{l-1}}{\b_{k-1}}M'+\biggl(w_{\g,\g}^{(\eta)}(-)+w_{\g,\eta}^{(\eta)}(-)+w_{\eta,\eta}^{(\eta)}(-)\biggr).\end{multline*}
\end{lemme}
\begin{proof}

This is a simple consequence of proposition~\ref{prop:1.3.4} applied to formulae~(\ref{gkl1}),~(\ref{ekl1}) and of lemma~\ref{u:2} (we have also used the fact that from lemma~\ref{m_1...} $|m_4|+|m_3|\leq 2(\b_{l-1}/\b_{k-1})$). 
\end{proof}
\begin{prop}\lab{c1'}
$$\lim_{k\to\infty}\b_{k-1}^2\max(\|\pa \g_k\|_{L^1({\T})}),\|\pa \eta_k\|_{L^1({\T})})=0$$
\end{prop}
\begin{proof}
In view of formulas~(\ref{gkl}), (\ref{ekl}) we get
\begin{align*}&f_k^\mp\geq |m_1|f_l^\mp+|m_2|e_l^\pm\\
&e_k^\pm\geq |m_3| f_l^\mp+|m_4| e_l^\pm\end{align*}
and consequentely from lemma~\ref{m_1...} and the definition of $u_i^\pm$ 
we get
$$\b_{k-1}u_k^\pm\geq \b_{l-1}u_l^\pm.$$
Also from proposition~\ref{prop:4.1} and lemmas~\ref{u},~\ref{u:1}
$$w_{\g^\mp,\g^\mp}^{(\g)}+w_{\g^\mp,\eta^\pm}^{(\g)}+w_{\eta^\pm,\eta^\pm}^{(\g)}\leq\frac{4}{\d}\biggl(f_k^{\mp}-(m_1f_l^\mp+|m_2|e_l^\pm) \biggr)^{1/2}.\biggl(2M(m_1+|m_2|)\biggr)^{3/2}$$ 
$$w_{\g^\mp,\g^\mp}^{(\eta)}+w_{\g^\mp,\eta^\pm}^{(\eta)}+w_{\eta^\pm,\eta^\pm}^{(\eta)}\leq\frac{4}{\d}\biggl(e_k^{\pm}-(|m_3|f_l^\mp+m_4e_l^\pm) \biggr)^{1/2}.\biggl(2M(|m_3|+m_4)\biggr)^{3/2}$$ 
and in view of~(\ref{gkl1}),~(\ref{ekl1}) 
$$w_{\g,\g}^{(\g)}(+)+w_{\g,\eta}^{(\g)}(+)+w_{\eta,\eta}^{(\g)}(+)\leq\frac{4}{\d}\biggl(f_k-(m_1f_l+|m_2|e_l) \biggr)^{1/2}.\biggl(2M(m_1+|m_2|)\biggr)^{3/2}$$ 
$$w_{\g,\g}^{(\eta)}(+)+w_{\g,\eta}^{(\eta)}(+)+w_{\eta,\eta}^{(\eta)}(+)\leq\frac{4}{\d}\biggl(e_k-(|m_3|f_l+m_4e_l) \biggr)^{1/2}.\biggl(2M(|m_3|+m_4)\biggr)^{3/2}$$ 
Using the fact that $\sqrt{a}+\sqrt{b}\leq 2\sqrt{a+b}$ we get
\begin{multline*}
\sqrt{\a_k}\biggl(w_{\g^\mp,\g^\mp}^{(\g)}+w_{\g^\mp,\eta^\pm}^{(\g)}+w_{\eta^\pm,\eta^\pm}^{(\g)}\biggr)+\biggl(w_{\g^\mp,\g^\mp}^{(\eta)}+w_{\g^\mp,\eta^\pm}^{(\eta)}+w_{\eta^\pm,\eta^\pm}^{(\eta)}\biggr)\leq \\\frac{8}{\d}\biggl(4M\b_{k-1}^{-1}\biggr)^{3/2}\biggl(\a_k\biggl(f_k^\mp-(m_1f_l^\mp+|m_2|e_l^\pm\biggr)+\biggl(e_k^\pm-(|m_3|f_l^\mp+m_4e_l^\pm\biggr)\biggr)^{1/2}\end{multline*}
and using  lemma~\ref{m_1...} and the definition of $u_i,~u_i^\pm$  
\begin{multline*}
\sqrt{\a_k}\biggl(w_{\g^\mp,\g^\mp}^{(\g)}+w_{\g^\mp,\eta^\pm}^{(\g)}+w_{\eta^\pm,\eta^\pm}^{(\g)}\biggr)+\biggl(w_{\g^\mp,\g^\mp}^{(\eta)}+w_{\g^\mp,\eta^\pm}^{(\eta)}+w_{\eta^\pm,\eta^\pm}^{(\eta)}\biggr)\leq \\\frac{8}{\d}\biggl(4M\b_{k-1}^{-1}\biggr)^{3/2}\biggl(u_k^\pm-\frac{\b_{l-1}}{\b_{k-1}}u_l^\pm\biggr)^{1/2}.\end{multline*}
In the same way
\begin{multline*}\sqrt{\a_k}\biggl(w_{\g,\g}^{(\g)}(+)+w_{\g,\eta}^{(\g)}(+)+w_{\eta,\eta}^{(\g)}(+)\biggr)+\biggl(w_{\g,\g}^{(\eta)}(+)+w_{\g,\eta}^{(\eta)}(+)+w_{\eta,\eta}^{(\eta)}(+)\biggr)\leq \\\frac{8}{\d}\biggl(4M\b_{k-1}^{-1}\biggr)^{3/2}\biggl(u_k-\frac{\b_{l-1}}{\b_{k-1}}u_l\biggr)^{1/2}\end{multline*}

Now a simple calculation shows that ($*=\g,\eta$)
$$W_{\g,\g}^{(*)}(-)(t)+W_{\g,\g}^{(*)}(+)(t)=2(W^{(*)}_{\g^-,\g^-}(t)+W^{(*)}_{\g^+,\g^+}(t)),$$
$$W_{\eta,\eta}^{(*)}(-)(t)+W_{\eta,\eta}^{(*)}(+)(t)=2(W^{(*)}_{\eta^-,\eta^-}(t)+W^{(*)}_{\eta^+,\eta^+}(t)),$$
$$W_{\g,\eta}^{(*)}(-)(t)+W_{\g,\eta}^{(*)}(+)(t)=2(W^{(*)}_{\g^-,\eta^+}(t)+W^{(*)}_{\g^+,\eta^-}(t)),$$
from which follows, using previous estimates, 
\begin{multline*}\sqrt{\a_k}\biggl(w_{\g,\g}^{(\g)}(-)+w_{\g,\eta}^{(\g)}(-)+w_{\eta,\eta}^{(\g)}(-)\biggr)+\biggl(w_{\g,\g}^{(\eta)}(-)+w_{\g,\eta}^{(\eta)}(-)+w_{\eta,\eta}^{(\eta)}(-)\biggr)\leq \\\frac{40}{\d}\b_{k-1}^{-2}(4M)^{3/2}(\bar v_{k,l})^{1/2}\end{multline*}
where $\bar v_{k,l}=\max(\bar u_k-\bar u_l,{\bar u}_k^+-{\bar u}_l^+,{\bar u}_k^--{\bar u}^-_l)$.

In view of lemma~\ref{c1} this implies that
$$\b_{k-1}^2\|\pa\eta_k\|_{L^1({\T})}\leq 2(\frac{\b_{k-1}}{\b_{l-1}})M'+\frac{40}{\d}(4M)^{3/2}(\bar v_{k,l})^{1/2}.$$
Taking $l=2[[k/2]/2]$ gives the conclusion of the proposition.
\end{proof}

\section{Convergence of the renormalization scheme}

Referring to formulae~(\ref{ekl}), (\ref{ekl1}) we set
\begin{multline}\r_{k,l}^\pm=\int_{{\T}}\biggl\|\biggl[\g_l^\mp(t_{|m_3|-1}),\Ad(C^{(l)}(t_{|m_3|})^{-1})).\g_l^\mp(t_{|m_3|})\biggr]\biggr\|dt+\\\int_{{\T}}\biggl\|\biggl[\g_l^\mp(t_{|m_3|}),\eta_l^\pm(t_{|m_3|+1})\biggr]\biggr\|dt+\\\int_{{\T}}\biggl\|\biggl[\eta_l^\pm(t_{|m_3|+1}),\Ad(A^{(l)}(t_{|m_3|+1})).\eta_l^\pm(t_{|m_3|+2})\biggr]\biggr\|dt\end{multline}
and 
\begin{multline}\r_{k,l}=\int_{{\T}}\biggl\|\biggl[(\g_l^-+\g_l^+)(t_{|m_3|-1}),\Ad(C^{(l)}(t_{|m_3|})^{-1})).(\g_l^-+\g_l^+)(t_{|m_3|})\biggr]\biggr\|dt+\\\int_{{\T}}\biggl\|\biggl[(\g_l^-+\g_l^+)(t_{|m_3|}),(\eta_l^++\eta_l^-)(t_{|m_3|+1})\biggr]\biggr\|dt+\\\int_{{\T}}\biggl\|\biggl[(\eta_l^++\eta_l^-)(t_{|m_3|+1}),\Ad(A^{(l)}(t_{|m_3|+1})).(\eta_l^++\eta_l^-)(t_{|m_3|+2})\biggr]\biggr\|dt\end{multline}
where $t_{|m_3|-1}=t+(m_4-1)\b_l-2\b_{l-1}$, $t_{m_3}=t+(m_4-1)\b_l-\b_{l-1}$, $t_{m_3+1}=t+(m_4-1)\b_l$, $t_{m_3+2}=t+(m_4-2)\b_l$.
Clearly
$$\r_{k,l}^\pm\leq w^{(\eta)}_{\g^\mp,\g^\mp}+w^{(\eta)}_{\g^\mp,\eta^\pm}+w^{(\eta)}_{\eta^\pm,\eta^\pm}$$
$$\r_{k,l}\leq  w^{(\eta)}_{\g,\g}(+)+w^{(\eta)}_{\g,\eta}(+)+w^{(\eta)}_{\eta,\eta}(+)$$
Since Lebesgue measure is invariant by translation  we can also write
\begin{multline}\r_{k,l}^\pm=\int_{{\T}}\biggl\|\biggl[\g_l^\mp(t-2\b_{l-1}),\Ad(C^{(l)}(t-\b_{l-1})^{-1})).\g_l^\mp(t-\b_{l-1})\biggr]\biggr\|dt+\\\int_{{\T}}\biggl\|\biggl[\g_l^\mp(t-\b_{l-1}),\eta_l^\pm(t)\biggr]\biggr\|dt+\\\int_{{\T}}\biggl\|\biggl[\eta_l^\pm(t),\Ad(A^{(l)}(t)).\eta_l^\pm(t-\b_l)\biggr]\biggr\|dt\end{multline}
and 
\begin{multline}\r_{k,l}=\int_{{\T}}\biggl\|\biggl[(\g_l^-+\g_l^+)(t-2\b_{l-1}),\Ad(C^{(l)}(t-\b_{l-1})^{-1})).(\g_l^-+\g_l^+)(t-\b_{l-1})\biggr]\biggr\|dt+\\\int_{{\T}}\biggl\|\biggl[(\g_l^-+\g_l^+)(t-\b_{l-1}),(\eta_l^++\eta_l^-)(t)\biggr]\biggr\|dt+\\\int_{{\T}}\biggl\|\biggl[(\eta_l^++\eta_l^-)(t),\Ad(A^{(l)}(t)).(\eta_l^++\eta_l^-)(t-\b_l)\biggr]\biggr\|dt\end{multline}
so that these quantities do not depend on $k$ (but we do not remove the index $k$).
We shall then define the functions $r_{k,l}(t),r_{k,l}^\pm(t):{\T}\to [0,\infty)$ by the same formulae as above but where the integration  symbol is removed. 
Finally we introduce the corresponding tilded variables: 
$${\ti \g}_j^\pm(t)=\b_{j-1}\g_j^\pm(\b_{j-1}t),\qquad {\ti \eta}_j^\pm(t)=\b_{j-1}\eta_j^\pm(\b_{j-1}t)$$
and the corresponding
\begin{multline}{\ti r}_{k,l}^\pm(t)=\biggl\|\biggl[{\ti \g}_l^\mp(t-2),\Ad(C^{(l)}(t-1)^{-1})).{\ti \g}_l^\mp(t-1)\biggr]\biggr\|+\\\biggl\|\biggl[{\ti \g}_l^\mp(t-1),{\ti \eta}_l^\pm(t)\biggr]\biggr\|+\\\biggl\|\biggl[{\ti \eta}_l^\pm(t),\Ad(A^{(l)}(t)).{\ti \eta}_l^\pm(t-\a_l)\biggr]\biggr\|\end{multline}
\begin{multline}{\ti r}_{k,l}(t)=\biggl\|\biggl[({\ti \g}_l^-+{\ti \g}_l^+)(t-2),\Ad(C^{(l)}(t-1)^{-1})).({\ti \g}_l^-+{\ti \g}_l^+)(t-1))\biggr]\biggr\|+\\\biggl\|\biggl[({\ti \g}_l^-+{\ti\g}_l^+)(t-1)),({\ti \eta}_l^++{\ti\eta}_l^-)(t)\biggr]\biggr\|+\\\biggl\|\biggl[({\ti \eta}_l^++{\ti\eta}_l^-)(t),\Ad(A^{(l)}(t)).({\ti \eta}_l^++{\ti\eta}_l^-)(t-\a_l)\biggr]\biggr\|\end{multline}
We have
$$\int_0^{\b_{l-1}^{-1}}{\ti r}_{k,l}(t)dt=\b_{l-1}\r_{k,l},\qquad \int_0^{\b_{l-1}^{-1}}{\ti r}_{k,l}(t)dt=\b_{l-1}\r_{k,l}.$$
If for some $k$ the numbers  $\a_k,\a_{k+1}$ are in the interval $((1/5),(1/4))$ we have
$$\b_{k-1}^2(\r_{k+2,k}+\r_{k+2,k}^++\r_{k+2,k}^-)\leq 125\bar v_{k+2,k}$$
Having this in mind and the fact that
$$\lim_{k\to\infty}\b_{k-1}^2(\|\pa\eta_k\|_{L^1({\T})}+\|\pa\g_k\|_{L^1({\T})})=0$$
we introduce the quantity:
$$\e_k=\b_{k-1}^2\biggl(\|\pa\eta_k\|_{L^1({\T})}+\|\pa\g_k\|_{L^1({\T})}+\r_{k+2,k}+\r_{k+2,k}^++\r_{k+2,k}^-\biggr).$$
and the functions $f:{\T}\to [0,\infty)$
$$f(t)=\|\pa\eta_k(t)\|+\|\pa\g_k(t)\|+r_{k+2,k}(t)+r_{k+2,k}^+(t)+r_{k+2,k}^-(t),$$
and $\ti f:[0,\b_{k-1}^{-1}]\to [0,\infty)$
$$\ti f(t)=\|\pa{\ti \eta}_k(t)\|+\|\pa{\ti\g}_k(t)\|+{\ti r}_{k+2,k}(t)+{\ti r}_{k+2,k}^+(t)+{\ti r}_{k+2,k}^-(t).$$
We have
$$\int_0^{\b_{k-1}^{-1}}\ti f(t)dt=\b_{k-1}^{-1}\e_k.$$
Consequentely there exists an integer $0\leq j\leq [\b_{k-1}^{-1}]$ such that
$$\int_j^{(j+10)}\ti f(t)dt\leq 20\e_k.$$
Since $\a\in\S$ we can say that there exixts an infinite set ${\cal K}$ of integers  $k$ for which the numbers $\a_{k},\a_{k+1}$ are in $((1/5),(1/4))$ and consequentely $\e_k$ goes to zero as $k$ goes to infinity while staying in ${\cal K}$.
We have to keep in mind that in all the previous construction a shift parameter $\u\in{\T}$ was in fact present and that  the variables we called $\g_k,\g_k^\pm,\ldots$ corresponding to $A(\cdot)$ were in fact variables $\g_{k,\u},\g_{k,\u}^\pm,\ldots$ corresponding to $A(\cdot-\u)$. 
In other words there exists a sequence $\u_{k}$ ($k\in{\cal K}$) such that 
\begin{multline}\|\pa{\ti \eta}_{k,\u_k}\|_{L^1(-5,5)}+\|\pa{\ti\g}_{k,\u_k}\|_{L^1(-5,5)}+\sum_{h\in{\cal B}_{k,\u_k}}\|\ti h\|_{L^1(-5,5)}+\\\sum_{\ti h^+\in{\cal B}^+_{k,\u_k}}\|{\ti h}^+\|_{L^1(-5,5)}+\sum_{{\ti h}^-\in{\cal B}^-_{k,\u_k}}\|{\ti h}^-\|_{L^1(-5,5)}\leq 2\e_{k}
\end{multline}
 goes to zero as $k\to\infty,k\in{\cal K}$ where 
\begin{align*}{\cal B}_{k,\u_k}^\pm=&\biggl\{\biggl[{\ti \g}_{k,\u_k}^\mp(t-2),\Ad(C_{\u_k}^{(k)}(t-1)^{-1})).{\ti \g}_{k,\u_k}^\mp(t-1)\biggr],\quad \biggl[{\ti \g}_{k,\u_k}^\mp(t-1),{\ti \eta}_{k,\u_k}^\pm(t)\biggr],\\&\biggl[{\ti \eta}_{k,\u_k}^\pm(t),\Ad(A_{\u_k}^{(k)}(t)).{\ti \eta}_{k,\u_k}^\pm(t-\a_k)\biggr]\biggr\}\end{align*}
\begin{align*}{\cal B}_{k,\u_k}=\biggl\{\biggl[({\ti \g}_{k,\u_k}^-+{\ti \g}_k^+)(t-2),\Ad(C_{\u_k}^{(k)}(t-1)^{-1})).({\ti \g}_{k,\u_k}^-+{\ti \g}_{k,\u_k}^+)(t-1))\biggr],\\\biggl[({\ti \g}_{k,\u_k}^-+{\ti\g}_{k,\u_k}^+)(t-1)),({\ti \eta}_{k,\u_k}^++{\ti\eta}_{k,\u_k}^-)(t)\biggr],\\\biggl[({\ti \eta}_{k,\u_k}^++{\ti\eta}_{k,\u_k}^-)(t),\Ad(A_{\u_k}^{(k)}(t)).({\ti \eta}_{k,\u_k}^++{\ti\eta}_{k,\u_k}^-)(t-\a_k)\biggr]\biggr\}\end{align*}

On the other hand we know that the $C^s$-norms of ${\ti\g}_{k,\u_k}$, ${\ti\eta}_{k,\u_k}$, ${\ti\g}_{k,\u_k}^\pm$, ${\ti \eta}_{k,\u_k}^\pm$ on $[-5,5]$ are bounded by some constant $M_s$ which does not depend on $k$: the same is hence true for the elements of ${\cal B}_{k,\u_k}$, ${\cal B}_{k,\u_k}^\pm$, and so for $\ti f_{\u_k}$. The convexity inequalities imply that
\begin{prop}\lab{bbb}$$\lim_{k\in{\cal K}}\max(\|\pa\ti\g_{k,\u_k}\|_{C^s(-5,5)},\|\pa\ti\g_{k,\u_k}\|_{C^s(-5,5)})=0$$
and
$$\max_{{\ti h}\in{\cal B}_{k,\u_k},{\ti h}^\pm\in{\cal B}_{k,\u_k}^\pm} (\|\ti h\|_{C^s(-5,5)},\|{\ti h}^\pm\|_{C^s(-5,5)})=0$$
\end{prop}

\section{Back to the circle}
Applying the preceding result in case $s=1$ and with $k$ large enough in ${\cal K}$ (but fixed) we can say that for $t\in (-5,5)$
$$\pa(\ti C^{(k)})(\ti C^{(k)})^{-1}=\ti\g_k(0)+O_t(\nu_k),\qquad \pa(\ti A^{(k)})(\ti A^{(k)})^{-1}=\ti\eta_k(0)+O_t(\nu_k),$$
and
$$\ti C^{(k)}(t)=e^{t\ti\g_k(0)}\ti C^{(k)}(0)+O_t(\nu_k),\qquad\ti A^{(k)}(t)=e^{t\ti\eta_k(0)}\ti A^{(k)}(0)+O_t(\nu_k),$$
where $t\mapsto O_t(\nu_k)$ is some small 
smooth function with $C^{s_0}$-norm on the interval $(-5,5)$  less than some constant times $\nu_k$.
From proposition~\ref{bbb} we can then write
$$\|[{\ti \g}_k^\mp(0),{\ti\eta}_k^\pm(0)]\|=O(\nu_k)$$
$$\Ad({\ti C}^{(k)}(0)).{\ti \g}_k^\mp(0)={\ti \g}_k^\mp(0)+O(\nu_k)$$
$$\Ad({\ti A}^{(k)}(0)).{\ti \eta}_k^\pm(0)={\ti\eta}_k^\pm(0)+O(\nu_k)$$
Moreover, by construction ${\ti \g}_k^\mp(0)$ and ${\ti\eta}_k^\pm(0)$ are in ${\cal E}_{\d'}$ and have norms  bounded away from zero by some constant. Hence  the previous inequalities tell us that the matrices ${\ti C}^{(k)}(0)$, ${\ti A}^{(k)}(0)$ are nearly elliptic with axis ${\ti \g}_k^-(0)$, ${\ti \g}_k^+(0)$,${\ti\eta}_k^+(0)$ and ${\ti\eta}_k^-(0)$ which are thus almost colinear.
Since ${\ti\g}_k(0)={\ti\g}_k^+(0)-{\ti\g}_k^-(0)$ and ${\ti\eta}_k(0)={\ti\eta}_k^+(0)-{\ti\eta}_k^-(0)$ this implies
$$[{\ti \g}_k(0),{\ti\eta}_k(0)]=O(\nu_k)$$
and
$$\Ad({\ti C}^{(k)}(0)).{\ti \g}_k(0)={\ti \g}_k(0)+O(\nu_k)$$
$$\Ad({\ti C}^{(k)}(0)).{\ti \eta}_k(0)={\ti \eta}_k(0)+O(\nu_k)$$
$$\Ad({\ti A}^{(k)}(0)).{\ti \eta}_k(0)={\ti\eta}_k(0)+O(\nu_k)$$
$$\Ad({\ti A}^{(k)}(0)).{\ti \g}_k(0)={\ti \g}_k(0)+O(\nu_k)$$
Another consequence is that the directions of ${\ti\g}_k(0)$ and ${\ti\eta}_k(0)$ are inside ${\cal E}_{\d'}$
On the other hand we know that $(1,{\ti C}^{(k)}(\cdot))$ and $(\a_k,{\ti A}^{(k)}(\cdot))$ commute so that
$${\ti A}^{(k)}(t+1){\ti C}^{(k)}(t)={\ti C}^{(k)}(t+\a_k){\ti A}^{(k)}(t)$$ 
and in view of the preceding formulae this gives
$$e^{\a_k{\ti \g}_k(0)-{\ti\eta}_k(0)}=Id+O(\nu_k)$$
A simple computation (trigonalize the matrix in a complex hermitian base) then shows that
$\a_k{\ti\g}_k(0)-{\ti \eta}_k(0)+O(\nu_k)$ is an elliptic matrix: there exists $P$ a 2 by 2 invertible matrix of bounded norm (independently of $k$) and an integer $r$ such that
$$\a_k{\ti\g}_k(0)-{\ti \eta}_k(0)+O(\nu_k)=P\begin{pmatrix}0&-2\pi r_k\\2\pi r_k&0\end{pmatrix}P^{-1}.$$

We claim:

\begin{prop}\lab{prop:9.1}For $k\in{\cal K}$ large enough, the initial ${\Z}^2$-action $((1,Id),(\a,A))$ is conjugated to a ${\Z}^2$-action $(1,Id),(\a_k,U_k))$ such that
$$\lim_{k\in{\cal K}}\|U_k(\cdot)-E_{r}(\cdot)\|_{C^{s_0}}=0,$$
where $r$ is an integer the absolute value of which is equal to the absolute value of the degree of $A(\cdot)$.
\end{prop}
\begin{proof}

Observe that since $\ti\eta_k(0),\ti\g_k(0)$ are bounded the same is true for $r_k$. Now if $\nu_k$ is small enough and if we normalize the action we obtain a perturbation of a system $(\a_k,E_{r_k}(\th))$ which has to be of degree $r_k$. But the degree of an action is, up to the sign, invariant by conjugation.  
\end{proof}

Before concluding this section we give the following inequality which will be important in the next section:

\begin{prop}\lab{prop:9.2}If for all $\th\in{\T}$ $L(A(\th))\in{\cal E}_\d$ then the following inequality is  satisfied:
$$2\pi r\geq \int_{\T}N(L(A(\th)))d\th.$$
\end{prop}
\begin{proof}

We have
\begin{align*}\a_k\ti\g_k-\ti\eta_k&=(\a_k{\ti\g}_k^++{\ti\eta}_k^-)-(\a_k{\ti\g}_k^-+{\ti\eta}_k^+)\\
&={\ti u}_k^--{\ti u}_k^++O(\nu_k).\end{align*}
But the assumption of the proposition implies that we can make the choice $\g_0^+=\g_0^-=\eta_0^-=0$ and consequentely formulae~(\ref{gkl}), (\ref{ekl}) with $l=0$ and $k$ even imply $\g_k^+=0$ and $\eta_k^-=0$ that is ${\ti u}_k^-=0$. Since $|\a_k{\ti \g}_k(0)-{\ti\eta}_k(0)|$ converges to $2\pi r$ and since $u_k^+$ is increasing this gives the desired conclusion since
$${\ti u}_0^+=\int_{{\T}}N(LA(\th))d\th.$$  
\end{proof}
\section{Proof of the main theorem: cocycles homotopic to the identity}
We can now give a proof of our main theorem part A in case the cocycle $A$ is homotopic to the identity: For $\a\in \S$ we choose $s_0$ the least order of differentiability  for which the local theorem applies. Then we renormalize the action. Since $r_k=0$ this means that $\|U_k(\cdot)-Id\|_{C^{s_0}}\leq \e_{s_0}$ if $k\in{\cal K}$ is large enough (this is a consequence of the convexity inequalities). Now the fibered rotation number of $(\a_k,U_k)$ is diophantine w.r.t $\a_k$ since the rotation number is, the way we have defined it for actions, invariant by conjugation. Consequentely the local theorem applies to $(\a_k,U_k)$ which means that the action generated by $((1,Id),(\a,A))$ and hence the cocycle $(\a,A)$ is reducible.
\section{The case of the Schr\"odinger equation: density of positive Lyapunov exponents}\lab{sec:schrod}
In  that case the initial cocycle $A$ is of the form
$$A_V(\th)=\begin{pmatrix}V(\th)&1\\
                         -1&0\end{pmatrix}.$$
Notice that such a cocycle is homotopic to the identity. Reasonning by contradiction the same way  we did to derive corollary~\ref{dple} from the Main Theorem A, we can assume that  for some $\l$ arbitrarily small  $(\a,A_{V+\l})$ has $C^0$-bounded fibered products and a non zero rotation number rational w.r.t $\a$ and thus is  conjugated to  $\pm Id$: there exists a smooth map $B:{\T}\to SL(2,{\R})$ 
$$B(\th)=\begin{pmatrix}a(\th)&c(\th)\\b(\th)&d(\th)\end{pmatrix}$$
such that (setting $\ti V(\cdot)=V(\cdot)+\l$, $\ti A(\cdot)=A_{V+\l}(\cdot)$)
\begin{equation}\pm \ti A(\th)=B(\th+\a)^{-1}IdB(\th).\lab{10:0}\end{equation}
Notice that the degree of $B$ is then non zero.
Then we should have (we treat the case $\pm=+$)
\begin{equation}a(\th+\a)d(\th)=b(\th+\a)c(\th).\lab{10:1}\end{equation}
\begin{equation}d(\th+\a)c(\th)-c(\th+\a)d(\th)=1\lab{10:1'}\end{equation}
On the other hand  from~(\ref{10:0})
$$B'(\th+\a)\begin{pmatrix}\ti V(\th)&1\\-1&0\end{pmatrix}+B(\th+\a)\begin{pmatrix}{\ti V}'(\th)&0\\0&0\end{pmatrix}=B'(\th)$$
from which we deduce
$$b'(\th+\a)=d'(\th),\qquad a'(\th+\a)=c'(\th)$$
and consequentely there exixt real numbers $e,f$ such that
\begin{equation}b(\th+\a)=d(\th)+f,\qquad a(\th+\a)=c(\th)+f.\lab{10:2}\end{equation} 
Equalities~(\ref{10:1}) and (\ref{10:2})
imply
\begin{equation}ed(\th)=fc(\th)\lab{10:3}\end{equation}
and from~(\ref{10:1'}) we get
$$ef=0.$$
If for example $e\ne 0$ then $f=0$ and  equality~(\ref{10:3}) would imply $d(\th)=0$ and equality~(\ref{10:2}) gives $b(\th+\a)=0$. This contradicts $\det B(\th)=1$.
On the other hand if both $e$ and $f$ are zero
$$b(\th+\a)=d(\th),\qquad a(\th+\a)=c(\th).$$
Let us  now perturb $\ti V(\cdot)$ by some (small) smooth function of the form $\e w(\cdot)$ ($\e$ is small):
\begin{align*}B(\th+\a)&\begin{pmatrix}V(\th)+\e w(\th)&1\\-1&0\end{pmatrix}B(\th)^{-1}\\&=Id+\e w(\th)B(\th+\a)\begin{pmatrix}1&0\\0&0\end{pmatrix}B(\th)^{-1}\\
&=Id+\e \begin{pmatrix}w(\th)a(\th+\a)d(\th)&- w(\th)c(\th)a(\th+\a)\\w(\th)b(\th+\a)d(\th)&-w(\th)b(\th+\a)c(\th)\end{pmatrix}\\
&=Id+\e\begin{pmatrix}w(\th)c(\th)d(\th)&-w(\th)c(\th)^2\\w(\th)d(\th)^2&-w(\th)d(\th)c(\th)\end{pmatrix}.
\end{align*}
Let us introduce
$$\l=\int_{{\T}}w(\th)c(\th)d(\th)d\th,\qquad \u=\int_{{\T}}w(\th)d(\th)^2d\th,\qquad \nu=\int_{{\T}}w(\th)c(\th)^2d\th.$$
We want to find $w\in C^\infty({\T},{\R})$ such that
\begin{equation}\lab{hyperb}\l^2-\u\nu>0.\end{equation}
Since  the degree of $B(\cdot)$ is non-zero and is the index of $\th\mapsto (c(\th),d(\th))$, the functions $c(\cdot), d(\cdot)$ vanish but  not  at the same points (since $B$ is invertible). Let us denote by $x$, $y$ points such that $c(x)=0$, $d(y)=0$ ($x\ne y$). Now if $\chi$ is a positive smooth function with compact support and mass 1, we define,
$$w(\th)=\d^{-1}\chi(\d^{-1}(x-\th))-\d^{-1}\chi(\d^{-1}(y-\th)),$$
and we have for $\d$ small enough
$$\l\approx 0,\qquad \u\approx d(x)^2>0,\qquad \nu\approx -c(y)^2<0,$$
and so~(\ref{10:5}) is satisfied.

Let us denote 
$$A_\e(\th)=\begin{pmatrix}V(\th)+\e w(\th)&1\\-1&0\end{pmatrix}:$$
we get
$$B(\th+\a)A_\e(\th)B(\th)^{-1}=(Id+\e\bar A)+\e F(\th),$$
where $Id+\e\bar A$ is hyperbolic since~(\ref{hyperb}) is satisfied.
It is now a standard fact from perturbation theory that one can find smooth maps $Y_n,F_n:{\T}\to SL(2,{\R})$ and constants $\bar A_n\in SL(2,{\R})$ such that
\begin{multline*}\biggl(I+\sum_{k=1}^n\e^k Y_k(\th+\a)\biggr)\biggl((Id+\e\bar A)+\e F(\th)\biggr)=\\\biggl((Id+\sum_{k=1}^n\e^k\bar A_k)+\e^{n+1} F_n(\th)\biggr)\biggl(I+\sum_{k=1}^n\e^k Y_k(\th)\biggr),\end{multline*}
with (for example)
$$|\bar A_n-\bar A|\leq C_n(F),\qquad |F_n|_2\leq C_n(F).$$
The constants $C_n(F)$ which depend on the function $F$ (and its derivatives) can increase very fast and the above series do not converge in general (due to small divisors phenomenon) but this is not our problem here.
For $n$ fixed ($n=10$ forexample) and $\e$ small enough the constant matrix 
$$Id+\sum_{k=1}^n\e^k\bar A_k$$
is hyperbolic and satisfies some fixed cone condition. Hence the same is true for 
$$Id+\sum_{k=1}^n\e^k\bar A_k+\e^n F_n(\th)$$
which is thus hyperbolic.
This shows that $(\a,(Id+\e\bar A+\e F))$ is hyperbolic and hence has positive Lyapounov exponent.
\section{Cocycles non-homotopic to the identity}
In that section we assume that the initial system we are dealing with is non-homotopic to the identity and has degree $r\in {\Z}-\{0\}$ and our aim is to prove the following theorem the proof of which occupies the whole section.
\begin{theo}\lab{theo:12.1}Assume $\a\in\S$, $A:{\T}\to SL(2,{\R})$ is of degree $r$ and assume that the fibered products $A_n(\cdot)$ ($n\in{\Z}$) along $\th\mapsto\th+\a$ are $C^0$-bounded; then  $(\a,A(\cdot))$ is smoothly conjugated to a system $(\a,E_r(\cdot))$ where $|r|=\deg(A)$ and
$$E_r(\th)=\begin{pmatrix}\cos(2\pi r\th)&-\sin(2\pi r\th)\\\sin(2\pi r\th)&\cos(2\pi r\th)\end{pmatrix}$$
\end{theo}

With the notation of proposition~\ref{prop:9.1} for $k\in{\cal K}$ large enough $(\a,A)$ is conjugated to $(\a_k,U_k)$ with $\a_k\in CD(\g,\s)$ and 
$$U_k(\th)=E_r(\th)e^{F(\th)},$$
where $|r|$ is $|deg(A)|$ and
$F(\cdot)$ is of small $C^{s_0}$-norm.
We are then in what we can call a local (perturbative) situation: to what extent  a system $C^\infty$-close to $E_r(\cdot)$  is conjugated to $E_r(\cdot)$ or to any sytsem not too difficult to understand. Stated this way the problem is still too difficult to be solved. Nevertheless, as we shall see, one can obtain some useful information.

\subsection{Normal form}
Since $SL(2,{\R})$ and $SU(1,1)$ are isomorphic it is indifferent to work in either one of these two spaces, but for easier computations we shall work with $SU(1,1)$ representation.
In the whole subsection we assume that $r$ is a fixed integer, $\a\in CD(\g,\s)$ and that $A:{\T}\to SU(1,1)$ is smooth and close to $E_r(\cdot)$:
$$A(\th)=E_r(\th)e^{U(\th)},$$
with $F\in C^\infty({\T},su(1,1)),$
$$\|U\|_{s_0}\leq \e_0,$$
and where 
$$E_r(\th)=\begin{pmatrix}e^{2\pi i r\th}&0\\0& e^{-2\pi i r\th}\end{pmatrix}.$$
We now want to find some small $Z:{\T}\to su(1,1)$, $Z=\{y,z\}$ such that
$$e^{Z(\th+\a)}E_r(\th)e^{U(\th)}e^{-Z(\th)}=E_r(\th)e^{V(\th)},$$
where $V(\cdot)$ is the simplest possible. To do that we solve the linearized equation
$$\Ad(E_{-r}(\th)).Z(\th+\a)+U(\th)-Z(\th)=V(\th),$$
which amounts to solving the system
\begin{align}y(\th+\a)-y(\th)&=f_1(\th)-f(\th)
\\e^{-4\pi i r\th}\nu(\th+\a)-\nu(\th)&=g_1(\th)-g(\th)
\end{align} 
where we have written $U=\{f,g\}$, $V=\{f_1,g_1\}$. 
In turns out that
\begin{itemize}
\item the first equation has a solution if $\hat {f}_1(0)=\hat f(0)$
\item the second equation has a solution if $g_1$ is  some trigonometric polynomial of the form
$$P_r(\th)=\sum_{k=-(2r-1)}^{0}a_ke^{2\pi ik\th},$$
where $k_0$ can chosen to be any integer (a convenient choice for our purpose is $k_0=2r$).
\item if this is so, then the corresponding solutions are less regular than the initial data: there is a loss of differentiability (for the first equation this is the classical loss of derivatives due to small divisors).  
\end{itemize}
More precisely let us introduce the truncature and remainder operator up to order $N\in{\N}$:
$$T_NU(\th)=\sum_{|k|\leq N}\hat U(k)e^{2\pi k\th},\qquad R_NU(\th)=\sum_{|k|>N}\hat U(k)e^{2\pi k\th},$$
where $\hat U(k)$ are the Fourier coefficients of $U:{\T}\to su(1,1)$.
One can then prove the following theorem the proof of which is very similar to that of~\cite{K_annsmaths}, theorem 8.1 and which we shall then omit:
\begin{prop}\lab{theo:5.3}Let $\a$ be in $CD(\g,\s)$, $r\in {\N}-\{0\}$ , $N\in{\N}$. Then there exist $\e_0>0$, $s_0\in{\N}$, $a>0$ such that for any $U\in C^\infty({\R}/{\Z},su(1,1))$ satisfying,
\begin{equation}N^a\|U\|_{s_0}\leq\e_0,\lab{petit}\end{equation}
there exist functions $Z,F\in C^\infty({\R}/{\Z},su(1,1))$ such that ($h=\{1,0\}$),
$$e^{Z(\th+\a)}e^{2\pi  r\th h}e^{U(\th)}e^{-Z(\th)}=e^{2\pi  r\th h}e^{\G_{2r}U(\th)+F(\th)},$$
where $\G_{2r}U$ is of the form
$$\G_{2r}U=\{\hat f(0),\sum_{k=-(2r-1)}^0  \rho_ke^{2\pi i k\th}\},$$
and satisfies the following inequality:
$$\forall s\geq 0,\ \|\G_{2r}U\|_s\leq C_{s,r} \|U\|_s.$$
Moreover, for any $s'\geq s\geq 0$,
$$\|F\|_s\leq C_{s,r}\|U\|_s\|U\|_0N^a+C_{s',r}N^{-(s'-s-a)}\|U\|_{s'};$$
and,
$$\|Z\|_s\leq C_{s,r} N^a \|U\|_s,$$
where the  constant $C_{s,r}$ and $a$ are positive.
\end{prop}

Before beginning the proof of theorem~\ref{theo:12.1} we prove two lemmas.
\begin{lemme}\lab{hab}Assume that $\a\in ((1/5,1/4]$. There exists a constant $\e_1>0$ for which the following holds: If the fibered products of $(\a,E_r(\cdot)e^{U(\cdot)})$ are $C^0$-bounded then (with the notation $U=\{t,\nu\}$)
$$\int_{\T}|\nu(\th)|d\th\leq C\|\pa U\|_{C^0},$$
provided $\|U\|_{C^0}\leq \e_1$ (here $C$ is some positive constant independent on $U$).
\end{lemme}
\begin{proof}

Let us define 
$$v=\int_{\T}\nu(\th)d\th,\qquad \e=\biggl(\int_{\T}|\nu(\th)|^2d\th\biggr)^{1/2}$$
and assume that for $\l>0$ one has
%$$\int_{\T}|\pa_\th U(\th)|^2d\th<\l\e.$$
$$\|\pa_\th U(\th)\|_{C^0}<\l\e.$$
Then
$$\nu(\th)=v+O(\l\e).$$ 
Let us introduce
$$\g_0^\pm(\cdot)\equiv 0,\qquad \eta_0^-(\cdot)\equiv 0,\qquad \eta_0^+(\cdot)=(\pa_\th A(\cdot))A(\cdot)^{-1},$$
$$A^{(1)}(\cdot)=A(\cdot-4\a)^{-1}A(\cdot-3\a)^{-1}A(\cdot-2\a)^{-1}A(\cdot-\a)^{-1},$$
and
$$\g_1(\cdot)^\pm=\eta_0^\mp(\cdot)$$
then for $\th\in{\T}$
\begin{multline*}\eta_1^\pm(\th)=\Ad\biggl(A(\th-4\a)^{-1}\biggr).\eta_0^\mp(\th-4\a)+\\\Ad\biggl(A(\th-4\a)^{-1}A(\th-3\a)^{-1}\biggr).\eta_0^\mp(\th-3\a)+\\\Ad\biggl(A(\th-4\a)^{-1}A(\th-3\a)^{-1}A(\th-2\a)^{-1}\biggr).\eta_0^\mp(\th-2\a)+\\\Ad\biggl(A(\th-4\a)^{-1}A(\th-3\a)^{-1}A(\th-2\a)^{-1}A(\th-\a)^{-1}\biggr).\eta_0^\mp(\th-\a)\end{multline*}
with $A(\th)=E_r(\th)e^{U(\th)}$. Since the fibered products of $(\a,A)$ are bounded we get using   proposition~\ref{prop:9.2} we must have
\begin{equation}2\pi r\geq \int_0^1N({\ti\eta}_1^-(t))dt+\int_0^{\a_1}N({\ti\g}_1^+(t))dt.\lab{rare}\end{equation}
If we define
$$\phi(\th)=\eta_0^+(\th)+\Ad\biggl(A(\th)^{-1}\biggr).\eta_0^+(\th+\a),$$
we have
\begin{equation}N(\eta_1^-(\th))\geq N(\phi(\th-4\a))+N(\phi(\th-2\a)).\lab{aube}\end{equation}
Using the fact that
$$L(e^{U(\th)})=\pa_\th U(\th)+O(\l\e^2)$$
and 
$$\Ad(e^U).h=h+[U,h]+\frac{1}{2}[U,[U,h]]+O_3(U),$$
we get after some simple computation similar to that of~\cite{K_annsmaths} proposition 7.1 (but in a $su(1,1)$ setting),
$$\phi(\th)=\{M_1(\th),M_2(\th)\}$$
with
\begin{align*}M_1(\th)&=4\pi r(1+\e^2)+\pa_\th a(\th)+\pa_\th a(\th+\a)+O(\l\e^2)\\
M_2(\th)&=4\pi irv+O(\l\e)\end{align*}
and from that
$$N(\phi(\th))=4\pi r\sqrt{1+(\e^2/2)}+\pa_\th a(\th)+\pa_\th a(\th+\a)+O(\l\e^2).$$
Using~(\ref{aube}) and  the previous estimate,  integration over ${\T}$ and the fact that $\a\a_1=1-4\a$ gives after some calculation
$$\int_0^1N({\ti\eta}_1^-(t))dt+\int_0^{\a_1}N({\ti\g}_1^+(t))dt=2\pi r+{2\pi r\ov 5}\e^2+O(\l\e^2),$$
and~(\ref{rare}) provides us with a contradiction as soon as $\l$ is small enough.
%$$\Ad(e^U).h=\{1,0\}+2\{0,-i\nu\}+2\{|\nu|^2,t\nu\}+O_3(U),$$
%and consequentely 
%\begin{align*}N(A_2(\th))^2&=(2+2|\nu|^2)^2-4|\nu|^2(1+t^2)+O_3(U)\\
%&=4+4|\nu|^2+O_3(U).\end{align*}
%Thus
%$$N(A_2(\th))=2+|\nu|^2+O_3(U).$$
%Integrating w.r.t $\th$ on ${\T}$ and using~(\ref{rare}) gives the conclusion of the proposition.

%For the proof we need the following result which is due to M.R Herman and in an improved form to A. Avila and J. Bochi (we refer to the very elegant paper~\cite{A-B}): if for a matrix $A$ in $SL(2,{\R})$  we set (here $\|A\|$ denotes the operator norm)
%$${\cal N}(A)=\log\biggl(\frac{\|A\|+\|A^{-1}\|}{2}\biggr),$$
%and if we denote by $R_\th$ the rotation matrix  of angle $2\pi\th$ in the real plane,
%then for any sequence of matrices $A_1,\cdots ,A_n$ in $SL(2,{\R})$ one has
%$$\int_{\T}{\cal N}(A_nR_\th\cdots A_1R_\th)d\th=\sum_{j=1}^n{\cal N}(A_j).$$
%In the $SU(1,1)$ representation we set for $B\in SU(1,1)$,
%${\cal N}(B)={\cal N}(P^{-1}BP)$ and it is then clear that
%$$\int_{\T}{\cal N}\biggl(E_r(\th+(n-1)\a)e^{U(\th+(n-1)\a)}\cdots E_r(\th)e^{U(\th)}\biggr)d\th=\sum_{j=1}^n{\cal N}(e^{U(\th+j\a)}).$$
%As a consequence if for every $\th$ one has $|\nu(\th)|\geq (1/2)|\hat\nu(0)|>0$ then ${\cal N}(e^{U(\th+j\a)})$ is, uniformly in $n$, bounded from below and the fibered Lyapunov exponent of $(\a,E_r(\cdot)e^{U(\cdot)})$ is positive. But this condition is clearly implied by inequality~({\ref{cim}) for some fixed $\l>0$. 

\end{proof}
We also give the following normalization lemma the proof of which can be found in~\cite{K_annsmaths} Lemma 8.2:

\begin{lemme}\lab{cor:5.*}If $U=\{t,\nu\}\in C^{\infty}({\R}/{\Z},su(1,1))$ satisfies for some $\e_0$ independent of $U$,
$$\|t-\int_{{\T}^1}t(\th)d\th\|_0+\|\nu\|_0\leq \e_0,$$
then there exists $\th_1\in{\R}/{\Z}$ such that,
$$e^{2\pi r(\th-\th_1)h}e^{U(\th-\th_1)}=e^{2\pi r\th h}e^{U_1(\th)},$$
with $U_1(\cdot)=\{t_1(\cdot),\nu_1(\cdot)\}$ satisfying,
$$\hat t_1(0)=\int_{{\T}^1}t_1(\th)d\th =0,$$
and,
$$\|t_1\|_s+\|\nu_1\|_s\leq C_s(\|t-\int_{{\T}^1}t(\th)d\th\|_s+\|\nu\|_s).$$
\end{lemme}

The key theorem is then the following:

\begin{theo}\lab{theo:5.4} Let $\a\in{\R}/{\Z}$, $r\in{\N}-\{0\}$. There exists $\e_1>0$ such that if $U=\{t,\nu\}\in C^\infty({\R}/{\Z},su(1,1))$ satisfies,
$$\hat t(0)=0,$$
and,
$$\|U\|_{C^2}\leq \e_1,$$
and if $(\a,E_r(\cdot)e^{U(\cdot)})$ has $C^0$-bounded fibred-products
then:
for any $s\geq 0$,
\begin{equation}\|\Lambda_{2r}U\|_s\leq C_{s,r}.\|U-\Lambda_{2r}U\|_s,\lab{100}\end{equation}
where we have defined
$$\Lambda_{2r}\{t,\nu\}=\{\hat t(0),\sum_{k=-(2r-1)}^0\hat\nu(k)e^{2i\pi k\th}\},$$
(in our case $\hat t(0)=0$).
\end{theo}
\begin{proof}

If $A(\th)=e^{2\pi r\th h}e^{U(\th)}$
we have
\begin{align*}L(A(\th))&=2\pi r h+\Ad(e^{2\pi r h\th}).L(e^{U(\th)})\\
&=2\pi r h+\Ad(e^{2\pi r h\th}).({e^{\ad(U(\th))}-\Id\ov \ad(U(\th))}).\pa_\th U(\th).\end{align*}
Now,
$$({e^{\ad(U(\th))}-\Id\ov \ad(U(\th))}).\pa_\th U(\th)=\pa_\th U(\th)+{1\ov 2}\ad(U(\th)).\pa_\th U(\th)+Q_1(\th),$$
where $Q_1(\th)=\{p_1(\th),q_1(\th)\}$ is a converging sum of terms of degree at least 3 in $(U,\pa_\th U)$, more precisely for some $C>0$,
$$\forall\th\in{\T}^1,\ |Q_1(\th)|\leq C(|U(\th)|^2|\pa_\th U(\th)|).$$
We now write,
$$U(\th)=\begin{pmatrix}it(\th)&\nu(\th)\\\bar \nu(\th)&-it(\th)\end{pmatrix}=\{t(\th),\nu(\th)\},$$
$$\pa_\th U(\th)=\{\pa_\th t(\th),\pa_\th \nu(\th)\},$$
and,
$${1\ov 2}[U,\pa_\th U]=\{-{\rm Im}((\pa_\th \nu).\bar\nu),-i(t\pa_\th\nu+\nu\pa_\th t)\}.$$
Hence,
\begin{align*}L(A(\th))&=\Ad(e^{2\pi rh\th}).\biggl(2\pi rh+\pa_\th U+{1\ov 2}[U,\pa_\th U]+Q_1(\th)\biggr)\\&=\{2\pi r+\pa_\th t-{\rm Im}((\pa_\th\nu).\bar\nu)+p_1(\th),\pa_\th\nu-i(t\pa_\th\nu+\nu\pa_\th t)+q_1(\th)\},\end{align*}
which implies,
$$q(L(A))=|2\pi r+\pa_\th t-{\rm Im}((\pa_\th\nu).\bar\nu)+p_1(\th)|^2+|\pa_\th\nu-i(t\pa_\th\nu+\nu\pa_\th t)+q_1(\th)|^2,$$
and consequently,
\begin{multline*}q(L(A))=(2\pi r)^2+2(2\pi r)\pa_\th t-2(2\pi r){\rm Im}(\bar\nu\pa_\th\nu)+(\pa_\th t)^2+\\
-|\pa_\th\nu|^2+\ti Q_2(\th),\end{multline*}
where,
$$\forall\th\in{\T}^1,\ |\ti Q_2(\th)|\leq C(|U(\th)|^2|\pa_\th U(\th)|+|U(\th)|.|\pa_\th U(\th)|^2).$$
As a consequence,
\begin{multline*}N(L(A))=(2\pi r)+\pa_\th t-{\rm Im}(\bar\nu\pa_\th\nu)+\frac{1}{4\pi r}\biggl((\pa_\th t)^2
-|\pa_\th\nu|^2\biggr)+\\\ti Q_2(\th),\end{multline*}
where $Q_2(\cdot)$ satisfies an inequality similar to the one satisfied by $\ti Q_2(\cdot)$.
We now integrate this inequality over the circle ${\R}/{\Z}$ ($\hat t(0)=0$):
\begin{multline*}\int_{{\R}/{\Z}}N(L(A(\th)))d\th=2\pi r-\int_{{\R}/{\Z}}{\rm Im}(\bar \nu\pa_\th\nu)d\th+\\
+\frac{1}{4\pi r}\biggl(\int_{{\R}/{\Z}}(\pa_\th t(\th))^2d\th-\int_{{\R}/{\Z}}|\pa_\th \nu(\th)|^2d\th\biggr)+\int_{{\T}^1}|Q_2(\th)|d\th.\end{multline*}
If we write the Fourier series expansion $\nu(\th)=\sum_{k\in{\Z}}\nu_k e^{2\pi ik \th}$, we get,
$$\int_{{\T}^1}|\pa_\th\nu|^2d\th=(2\pi)^2\sum_{k\in {\Z}}|k|^2|\nu_k|^2,\qquad\int_{{\T}^1}{\rm Im}(\bar\nu\pa_\th\nu)d\th=2\pi\sum_{k\in {\Z}}k|\nu_k|^2.$$
and consequentely,
\begin{multline*}\int_{{T}^1}N(L(A(\th)))d\th\geq \\2\pi r+\frac{\pi}{r}\sum_{k=-(2r-1)}^{-1}|\nu_k|^2|k|(2r-|k|)-\int_{{\T}^1}(\pa_\th t(\th))^2d\th+\\
-3\sum_{k\in{\Z}-\{0,\ldots,-(2r-1)\}}|k|^2|\nu_k|^2-\int_{{\T}^1}|Q_2(\th)|d\th.\end{multline*}
The inequality satisfied by $Q_2(\cdot)$ can be integrated over ${\R}/{\Z}$ and one can prove (see \cite{K_annsmaths} lemma 8.3 for details) that 
$$\int_{\T}|Q_2(\th)|d\th=O(\e_1)\int_{\T}|U(\th)|^2d\th,$$
hence
$$\int_{\T}|Q_2(\th)|d\th=O(\e_1)\sum_{k\in{\Z}}(|\hat t(k)|^2+|\hat\nu_k|^2).$$
Finally:
\begin{multline*}\int_{{T}^1}N(L(A(\th)))d\th\geq 2\pi r+(\frac{\pi}{r}-O(\e_1))\sum_{k=-(2r-1)}^{-1}|\nu_k|^2-\\
-(3+O(\e_1))\|\pa_\th(U-\Lambda_{2r}U)\|^2_{L^2}+O(\e_1)|\nu_0|^2,\end{multline*}
and if $\e_1$ is small enough
\begin{multline*}\int_{{T}^1}N(L(A(\th)))d\th\geq 2\pi r+\frac{\pi}{2r}\sum_{k=-(2r-1)}^{-1}|\nu_k|^2\\
-4\|\pa_\th(U-\Lambda_{2r}U)\|^2_{L^2}+O(\e_1)|\nu_0|^2,\end{multline*}
and from~\ref{prop:9.2} we get
$$\frac{\pi}{2r}\sum_{k=-(2r-1)}^{-1}|\nu_k|^2\leq
4\|\pa_\th(U-\Lambda_{2r}U)\|^2_{L^2}+O(\e_1)|\nu_0|^2.$$
This in turn implies,
$$\|\pa_\th(\Lambda_{2r}U)\|_0^2\leq 50r^3\|\pa_\th(U-\Lambda_{2r}U)\|^2_{L^2}+O(\e_1)|\nu_0|^2$$
By lemma~\ref{hab} we know that 
\begin{align*}|\nu_0|^2&\leq \l\|\pa_\th U\|_0^2\\
&\leq \l(\|\pa_\th(U-\Lambda_{2r}U)\|_0+\|\pa_\th\Lambda_{2r}U\|_0)^2\\
&\leq C_r \|\pa_\th(U-\Lambda_{2r}U)\|_0^2+O(\e_1)|\nu_0|^2\end{align*}
and hence
$$\|\Lambda_{2r}U\|_{L^2}\leq C_r\|\pa(U-\Lambda_{2r}U)\|_{L^2}$$
which clearly is the conclusion of the proposition if one observes that
$$\|\Lambda_{2r}U\|_{C^s}\leq (2r)^{s+1}\|\Lambda_{2r}U\|_{L^2}.$$
\end{proof}

\subsection{Proof of theorem~\ref{theo:12.1}}
The proof of theorem~\ref{theo:12.1} follows from the iterated use of the following proposition (with for example $N_n=\exp((3/2)^n)$) the proof of which is very similar to Proposition 9.1 of\cite{K_annsmaths}. 
\begin{prop}Assume that $A(\cdot)$ is non homotopic to the identity and the fiberd products of $(\a,A)$ are $C^0$-bounded. Then, there exist $\e_1>0$, $a>0$ and an integer $s_0>0$ such that for any $U\in C^\infty({\T}^1,sl(2,{\R}))$ satisfying,
$$\|U\|_{s_0}N^a\leq \e_1,$$
the following is true: there exist $Z_1,U_{2}\in C^\infty({\T}^1,su(2))$ such that,
$$(\a,e^{2\pi r\th h}e^{U(\th)})\Conj(e^{Z_1}) (\a,e^{2\pi r\th h}e^{U_{2}(\th)}),$$
and for which the following facts are true: for any $s'\geq s\geq 0$,
$$\|U_{2}\|_s\leq C_{s,r}N^a\|U\|_s\|U\|_0+C_{s',r}N^{-(s'-s-a)}\|U\|_{s'},$$
$$\|Z_1\|_s\leq C_{s,r} N^a\|U\|_s.$$
\end{prop}

The final arguments to conclude can be found in~\cite{K_annsmaths}, \cite{K_ast}.

\newpage
\section{Appendix}
\subsection{Complexification}

We denote by ${\bf H}_\pm$ the Poincar\'e half-planes:
$${\bf H}_\pm=\{w\in{\C},\ \pm\Im w>0\},$$
and endow them with the Poincar\'e metric:
$$\<w_1,w_2\>_z=\Ree {w_1\bar w_2\ov |\Im z|^2}.$$
The open unit disk ${\bf D}=\{z\in{\C},|z|<1\}$, the complement of its closure and ${\bf H}_\pm$ are conformally equivalent and we can then define a Poincar\'e metric on each of these spaces. A conformal equivalence between ${\bf D}$ and ${\bf H}_-$ is given for example by the map from ${\bf H}_-$ to ${\bf D}$ $w\mapsto (w-i)(i+w)^{-1}$ with inverse $z\mapsto i(1+z)(1-z)^{-1}$.

There is a natural action of $SL(2,{\C})$ on $\bar{\C}={\C}\cup\{\infty\}$ namely for $m\in\bar{\C}$
$$\begin{pmatrix}a&b\\c&d\end{pmatrix}.m=\frac{am+b}{cm+d},$$
and it thus induces  actions for the  subgroups $SL(2,{\R})$ and $SU(1,1)$.  
%There are   natural action of $SL(2,{\R})$ and $SU(1,1)$ on ${\bf C}$ namely
%$$\begin{pmatrix}a&b\\c&d\end{pmatrix}.m=\frac{am+b}{cm+d},\qquad \bm u&\bar v\\ v&\bar u\em .\t=\frac{u\t+\bar v}{ v\t+\bar u},$$
%with $a,b,c,d\in{\R}$ satisfying $ad-bc=1$ and $u,v\in{\C}$ satisfying $|u|^2-|v|^2=1$. 
The spaces ${\bf H}_\pm$ are preserved by the $SL(2,{\R})$-action, the spaces ${\bf D}$, ${\C}\cup\{\infty\}-\bar {\bf D}$ are preserved by the $SU(1,1)$-action and on these spaces the correponding actions preserve the Poincar\'e metric.  The converse is true: if the action of $V\in SL(2,{\C})$ on ${\C}$ sends the disk ${\bf D}$ into itself then either the unit circle is preserved and $V$ is in $SU(1,1)$ or the action of $V$ on ${\bf D}$ contracts the Poincar\'e metric of ${\bf D}$. For any such $V$ one can define a function $\r_V(\cdot):\bar{\bf D}\to {\C}$ such that
$$V.\begin{pmatrix}\t\\1\end{pmatrix}=\r_V(\t)\begin{pmatrix}\ti\t\\1\end{pmatrix}$$
and $\r_V(\t)$ is continuous and never zero. Since  $\bar{\bf D}$ is simply connected
 one can find a lift $\phi_V(\cdot)$
\begin{align*}\bar{\bf D}&\to {\R}\\
\t\mapsto &\phi_z(\t)\end{align*}
such that $e^{i\phi_V(\t)}=\r_V(\t).$ For such a $V$ contracting the Poincar\'e metric of the disk there exists a unique fixed point $p_V\in{\bf D}$ for the action of $V$. But the $SU(1,1)$-matrix
$$U=\begin{pmatrix}u&v\\\bar v&\bar u\end{pmatrix}$$
where 
$$u=\frac{1}{\sqrt{1-|p_V|^2}},\qquad v=-p_Vu,$$
sends $p_V$ to 0; therefore since $UVU^{-1}$ contracts the Poincar\'e metric on the disk it is of the form
$$UVU^{-1}=\begin{pmatrix}z&0\\0&z^{-1}\end{pmatrix},$$
with $z\in {\bf D}-\{0\}$. We shall denote by $D_z$ the diagonal matrix in the preceding formula. Since any such $D_z$ contracts the Poincar\'e metric of the disk we have proved that any matrix $V$ contracting the Poincar\'e metric of the disk is conjugated to a matrix $D_z$, the conjugacy being is $SU(1,1)$.
 
Given $U\in SU(1,1)$
$$\bm u&\bar v\\ v&\bar u\em,\qquad (|u|^2-|v|^2=1)$$
and $z\in\bar{\bf D}-\{0\}$
we define $U_z=UD_z$. It is  a standard fact that for  $0<|z|<1$ the map $D_z$  sends ${\bf D}$   into itself and $D_z^{-1}$  sends ${\C}\cup\{\infty\}-\bar {\bf D}$ into itself while contracting the Poincar\'e metrics. For any $\t\in\bar{\bf D}$ the complex number $\r_z(\t)=\r_{U_z}(\t)$  satisfies 
$$\r_z(\t)=zv\t+z^{-1}\bar u.$$
The function $\r_z(\t)$ is never zero continuous w.r.t $\t$ and holomorphic in $z\in{\bf D}-\{0\}$. Since we shall need it when we deal we rotation numbers we give the following lemma:
\begin{lemme}\lab{l1}For every $z\in\bar{\bf D}_0$, $\t_1\in{\bf D}$ and $\t_2\in\bar{\bf D}$ 
$$|\phi_z(\t_1)-\phi_z(\t_2)|<{1\ov 2}.$$
\end{lemme}
\begin{proof}

Assume on the contrary that there exist $\t_1\in{\bf D}$, $\t_2\in\bar {\bf D}$ such that the inequality is not satisfied. Then there are two points $\t_1^*,\t_2^*$ in the line segment $[\t_1,\t_2]$ such that 
$$\phi_z(\t_1^*)-\phi_z(\t_2^*)={1\ov 2}$$
and $\{\t_1,\t_2\}\ne\{\t_1^*,\t_2^*\}$.
If 
$$U_z(x)=\bm u& \bar v\\ v&\bar u\em D_z,$$
with $|u|^2-|v|^2=1$ this implies that there exists a positive number $t$ such that
$$vz\t^*_1+\bar u z^{-1}=-t(vz\t^*_2+\bar u z^{-1}),$$
that is
$$ vz(\t^*_1+t\t^*_2)=-\bar u z^{-1}(1+t),$$
and therefore
$$|v||z||\t^*_1+t\t^*_2|=|u||z^{-1}|(1+t).$$
Since $|u|>0$, $|v|\leq |u|$ and $|z|\leq 1$ we get
$$1+t\leq |\t_1^*+t\t^*_2|;$$
but at least one of the points $\t_1^*,\t_2^*$ is of module less than 1. This is a contradiction.
\end{proof}

A corollary is

\begin{cor} \lab{l2}If $V\in SL(2,{\C})$ preserves or contracts the Poincar\'e metric of the disk then for every  $\t_1,\t_2\in\bar{\bf D}$
$$|\phi_{V}(\t_1)-\phi_{V}(\t_2)|\leq 1.$$
\end{cor}
\begin{proof}

If $V$ preseves the Poincar\'e metric of the disk it is in $SU(1,1)$ and the conclusion follows from the preceding lemma. 

Otherwise, we know that $V$ is of the form
$$V=UD_zU^{-1},\qquad U\in SU(1,1),\quad 0<|z|<1.$$
But for any $\t_1,\t_2\in\bar{\bf D}$ there exist $\ti\t_1,\ti\t_2\in\bar {\bf D}$ such that ($i=1,2$)
$$\phi_V(\t_i)=\phi_{U^{-1}}(\t_i)+\phi_{UD_z}(\ti\t_i),$$
and the preceding lemma provides the conclusion of the lemma.
\end{proof}

The same discussion applies for the action of $SL(2,{\R})$ on each of the Poincar\'e half-planes. Just to set some notations we define for $A\in SL(2,{\R})$ and $0<|z|<1$ the matrix $A_z=AC_z$ where
$$C_z=\begin{pmatrix}e(z)&-f(z)\\
f(z)&e(z)\end{pmatrix}\in SL(2,{\C})$$
with $$e(z)={1\ov 2}(z+z^{-1}),\qquad f(z)={1\ov 2i}(z-z^{-1}),$$
so that when $z=e^{i\b}$ ($\b\in{\R}$) $C_z$ is just a real rotation matrix.
Notice that $C_z=P^{-1}D_zP$ where  
$$P=\bm -1&-i\\ -1 &i\em$$
is the matrix giving the group isomorphism
\begin{align*}SL(2,{\R})&\to SU(1,1)\\
A&\mapsto PAP^{-1}\end{align*}

For $0<|z|<1$, $C_z$ sends  ${\bf H}_-$ into itself and $C_z^{-1}$ send ${\bf H}_+$  into itself while contracting the Poincar\'e metrics.

\subsection{The complex rotation number}
In that section we assume that $X$ is a compact metric space endowed with a probability measure $\u$,  that $T:X\to X$ is a $\u$-ergodic  homeomorphism and that $A:{X}\to SL(2,{\R})$ is continuous. The case we are mainly concerned with is $X={\T}$, $\u$ is the Haar measure on ${\T}$ and  $T$ is the minimal translation $x\mapsto x+\a$. As usual $(T,A)$ is the map defined on the product $X\times SL(2,{\R})$ by $(T,A)(x,y)=(Tx,A(x)y)$.

We  define 
$$A_z(x)=\bm a_z(x)&b_z(x)\\ c_z(x)&d_z(x)\em =A(x)C_z,$$
which is in $SL(2,{\C})$ and, in the $SU(1,1)$ representation, 
$$U_z(x)=PA_z(x)P^{-1}=U(x).D_z$$
with
$$D_z=\bm z&0\\0&z^{-1}\em .$$
For every $x\in X$, $m\in{\bf H}_-$ (resp. $\t\in{\bf D}$) and $0<|z|\leq 1$ there exist  complex numbers $\g(z,x,m)$ and  $\r(z,x,\t)$) such that 
$$A_z(x)\bm m\\ 1\em=\g(z,x,m)\bm \ti m\\1\em,\qquad U_z(x)\bm \t\\1\em =\r(z,x,\t)\bm\ti\t\\ 1\em .$$ 
In particular 
$$zv(x)\t+z^{-1}\bar u(x)=\r(z,x,\t).$$
The function $\r(z,x,\t)$ is well defined on $0<|z|\leq 1$, never zero, continuous w.r.t $(x,\t)$, holomorphic w.r.t $z$ ($0<|z|\leq 1$); notice also that since $z\r(z)$ is bounded it admits a extension through $z=0$. 
For convenience we introduce the simply connected open  domain ${\bf D}_0$ which is the set of $\r e^{i\th}$, $0<\r<1$, $\th\in (-\pi,\pi)$.

\subsubsection{Using the hyperbolic structure: the $m$ and $\t$ functions}\lab{uhs}
We define the map $F_z$:
\begin{align*}X\times {\bf D}&\to X\times{\bf D}\\
(x,\t)&\mapsto (Tx,U_z(x).\t)
\end{align*}
and we shall denote by the same letter the map obtained by replacing ${\bf D}$ with ${\bf H}_-$ and $U_z$ with $A_z$. This maps acts on the Banach space of contnuous sections: to a continuous $\t(\cdot):X\to {\bf D}$ one associates ${\cal F}\t(x)=U_z(T^{-1}x)\t(T^{-1}x)$. 
We have seen that for every $x\in{X}$ and every $0<|z|<1$ the matrix $A(x)C_z$  contracts the metric induced by the Poincar\'e metric on ${\bf H}_-$ and similarly for each $0<|z|<1$,  $U_z(x)$ contracts the Poincar\'e metric on the disk ${\bf D}$. An application of Picard fixed point theorem to the contracting mapping ${\cal F}$ (defined on the Banach  space of continuous sections) shows that for every such $z$ there are  (unique) continuous sections over ${\bf T}$ 
$$m_-(z,\cdot):{\bf T}\to {\bf H}_-,\qquad \t_-(z,\cdot):{\bf T}\to {\bf D}$$
such that for every $x\in {X}$
\begin{equation}A(x)C_z\begin{pmatrix}m_-(z,x)\\1\end{pmatrix}=\g(z,x,m_-(z,x))\begin{pmatrix} m_-(z,Tx)\\1\end{pmatrix} ,\lab{e:54}\end{equation}
and
\begin{equation}U(x)D_z\begin{pmatrix} \t_-(z,x)\\ 1\em =\r(z,x,\t_-(z,x))\bm \t_-(z,Tx)\\1\end{pmatrix}.\lab{e:55}\end{equation}
It is a classical fact that the functions $m_-(z,x)$ and $\t_-(z,x)$ are holomorphic w.r.t $z$ as long as $0<|z|<1$ .
%and therefore  $z\t(x,z)$, $zm(x,z)$ extend holomorphically through $z=0$ (for $\t$ it is clear since it is bounded).

Using Fatou's theorem it then appears that for each $x\in{X}$ and almost every $\b\in [0,2\pi]$ the radial limits
$$\lim_{r\to 1}m_-(re^{i\b},x),\qquad \lim_{r\to 1}\t_-(re^{i\b},x),$$
exist for (Lebesgue) almost every $\b$ and an application of Fubini's theorem allows us to claim that for Lebesgue a.e $\b$ there is a set of $x\in{X}$ of full Haar measure (depending on $\b$) such that the preceding limits exist. Therefore, passing to the limit shows that for a.a $\b$ equations~(\ref{328}),(\ref{330}) with $z=e^{i\b}$ are still satisfied for $\u$-a.a $x\in X$. We still denote by $\t(e^{i\b},x)$ the limit of $\t(re^{i\b},x)$ as $r\to 1$ whenever it exists.

Similarly, by considering the inverse maps $(A(x)C_z)^{-1}$, $(U(x)D_z)^{-1}$ one defines sections $m_+(z,\x)$, $\t_+(z,x)$ such that $m_+(z,x)\in{\bf H}_+$, $\t_+(z,x)\in \bar {\C}-{\bar D}$ and which satisfy
\begin{equation}(A(T^{-1}x)C_z)^{-1}\begin{pmatrix}m_+(z,x)\\1\end{pmatrix}=\g(z,x,m_+(z,x))\bm m_+(z,T^{-1}x)\\1\end{pmatrix},\lab{e:56}\end{equation}
and
\begin{equation}(U(T^{-1}x)D_z)^{-1}\begin{pmatrix} \t_+(z,x)\\ 1\end{pmatrix} =\r(z,x,\t_+(z,x))\begin{pmatrix} \t_+(z,T^{-1}x)\\1\end{pmatrix}.\lab{e:57}\end{equation}
All the preceding discussion applies to $m_+$, $\t_+$.
\subsubsection{The complex rotation number}
We define it for $z$ in the simply connected domain ${\bf D}_0$ as
\begin{equation}\z(z)=\int_{x\in X}\Log \left(\r(z,x,\t_-(z,x))\right)d\u(x).\lab{zeta}\end{equation}
Notice that this definition holds for a.a $z\in\pa{\bf D}_0$. Also, an application of the dominated convergence theorem  in~(\ref{zeta}) with $z=re^{i\b}$ insures that for a.a $\b\in [0,2\pi]$
$$\lim_{r\to 1} \z(re^{i\b})=\z(e^{i\b}).$$

For complex $z$, the Lyapunov exponent of $(T,A_z)$ is the limit
$$\l(T,A_z)=\lim_{n\to \pm\infty}\int_X\|A_{z,n}(x)\|d\u(x),$$
which exists and is also equal for $\u$-a.a $x$ to
$$\lim_{n\to \pm\infty}\|A_{z,n}(x)\|.$$
This is a non negative number.
The map $z\to \l(T,A_z)$ is u.s.c. and therefore, if for $z_0$, $\l(T,A_{z_0})=0$ then $z\mapsto\l(T,A_z)$ is continous at $z_0$.
The following theorem identifies the real part of $\z$ with the Lyapunov exponent:
\begin{theo}For every $z\in{\bf D}_0$ and almost all $z\in\pa{\bf D}_0$ one has
$${\rm Re}\z(z)=\l(\a,A(\cdot)R_z).$$
\end{theo}
\begin{proof}

Equation~(\ref{330}) shows that for every $z\in{\pa D}_0$ (resp. a.e $z\in\pa{\bf D}_0$) one has for every $x\in X$ (resp. $\u$-a.e $x\in X$) 
\begin{equation}\|U_{z,n}(x)\begin{pmatrix}\t(z,x)\\1\end{pmatrix}\|=|\r_n(z,x,\t(z,x))| \|\begin{pmatrix}\t(z,T^nx)\\1\end{pmatrix}\|\lab{347}\end{equation}
where
$$\r_n(z,x,\t(z,x))=\r(z,x,\t(z,x))\cdots\r(z,T^nx,\t(z,T^nx)),$$
and by Birkhoff's ergodic theorem applied to the $L^\infty(\u)$ function $x\mapsto \Log(\r(z,x,\t(z,x)))$ we get
\begin{equation}{\rm Re}\z(z)=\lim_{n\to\infty}{1\ov n}\sum_{k=0}^{n-1}\Log|\r(z,T^kx,\t(z,T^kx))|.\lab{347'}\end{equation}
thus,
$$\lim_{n\to\infty}{1\ov n}\Log\|U_{z,n}\begin{pmatrix}\t(z,x)\\1\end{pmatrix}\|={\rm Re}\z(e^{i\b}).$$
Since for $0<|z|<1$ the section $\t(\cdot)$ is attracting, the limit in~(\ref{347'}) is positive for $z\in{\bf D}_0$ and hence for a.e $z\in\pa{\bf D}_0$ one has ${\rm Re} \z(z)\geq 0$ (just take the limit $r\to 1$ in~(\ref{zeta}) with $z=re^{i\b}$ and apply dominated convergence).  Hence, 
$$\lim_{n\to\infty}{1\ov n}\Log\|U_{z,n}(x)\begin{pmatrix}\t(z,x)\\1\end{pmatrix}\|=\lim{1\ov n}\Log\|U_{z,n}(x)\|,$$
and we have proved that  ${\rm Re}\z(z)=\l(T,U_z)$ for every $z\in{\bf D}_0$ and a.e $z\in\pa{\bf D}_0$.

The proof of the proposition is complete.
\end{proof}

\subsubsection{The geometric rotation number}
We now intend to give an interpretation of $\Im \z(z)$ for $z\in\pa{\bf D}_0$ in terms of classical rotation number. This requires more assumption on the map $T$ and on $A$.

\bn{\bf NB:} {\it Throughout this subsection we assume that $T$ is uniquely ergodic and that $A:X\to SL(2,{\R})$ is homotopic to the identity.}

It is then clear that there exists a continuous lift $\phi$
\begin{align*}\phi:\bar {\bf D}_0\times X\times \bar{\bf D}&\to {\R}\\
(z,x,\t)&\mapsto \phi_z(x,\t)\end{align*}
such that
$$e^{2\pi i\phi_z(x,\t)}=\frac{\r(z,x,\t)}{|\r(z,x,\t)|}.$$

\bigskip Let us still denote by $F_z$ the map introduced in~\ref{uhs}
and set for $n\geq 0$
\begin{align*}&\r_n(z,x,\t)=\r(z,(x,\t))\cdots\r(z,F_z^{n-1}(x,\t)),\\
&\phi_{z,n}(x,\t)=\phi_{z,n}(x,\t)+\cdots+\phi_{z,n}(F_z^{n-1}(x,\t))\end{align*}

We still have
$$e^{2\pi i\phi_{z,n}(x,\t)}=\frac{\r_n(z,x,\t)}{|\r_n(z,x,\t)|},$$
and
$$U_{z,n}(x)\bm \t\\ 1\em =\r_n(z,x,\t)\bm \ti \t\\1\em.$$
We have 
\begin{lemme}\lab{l3}For any $n\in{\N}$, $z\in\bar{\bf D}_0$, $\t_1,\t_2\in\bar{\bf D}$
$$|\phi_{z,n}(\t_1)-\phi_{z,n}(\t_2)|\leq 1.$$
\end{lemme}
\begin{proof}

Since $U_{z,n}(x)$ is a product of matrices preserving or contracting the Poincar\'e metric on the disk it has the same property and therefore corollary~\ref{l2} applies
\end{proof}
We can now prove
\begin{theo}Assume $T$ is uniquely ergodic and $A$ homotopic to the identity. For every $z\in\pa{\bf D}_0$, $x\in X$, $\t\in\bar{\bf D}$ the following limit
$$\lim_{n\to\infty}\frac{1}{n}\phi_{z,n}(x,\t)$$
exists, does not depend on $x$ and $\t$ and is uniform in $(x,\t)\in X\times\bar{\bf D}$. We call it the geometric rotation number of $(T,A_z)$ and denote it by $\r_g(T,A_z)$. Moreover, as a consequence of uniform convergence, this geometric rotation number is continuous w.r.t $z\in\bar{\bf D}_0$. Also, when $z\in\pa{\bf D}_0$ this number coincide with the usual notion of fibered rotation number.
\end{theo}
\begin{proof}

We follow the method in~\cite{He}. 
Let $\nu$ be some invariant probability measure  on $X\times \bar{\bf D}$ for the mapping $F_z$. Since $\u$ is uniquely ergodic $\nu$ projects down to $\u$. By Birkhoff ergodic theorem, for a set of full $\nu$-measure $(x,\t)$ one has the following equality
$$\lim_{n\to\infty}{1\ov n}\sum_{k=0}^{n-1}\phi_z(F_z^k(x,\t))=h(x,\t).$$
By lemma~\ref{l3} if the sum in the preceding equality converges to $h(x_0,\t_0)$ for some  $(x_0,\t_0)$ then it converges to the same limit $h(x_0)$ for any $(x_0,\t)$, $\t\in\bar{\bf D}$; the set of $x\in X$ for which $h(x)$ exists is invariant by $T$, $h(\cdot)$ is also $T$ invariant on $B$: the measure $\u$ being uniquely ergodic this means that there is a set of full $\u$-measure $B$ such that the preceding equality holds for $(x,\t)\in B\times \bar{\bf D}$ and the limit is independent of $(x,\t)\in B\times\bar{\bf D}$. Hence  
$$\lim_{n\to\infty}{1\ov n}\sum_{k=0}^{n-1}\phi_z(F_z^k(x,\t))=\int_{X\times \bar{\bf D}}\phi_z(x,\t)d\u(x,\t).$$
Consequentely, if $\nu_1$, $\nu_2$ are two invariant measures on $X\times \bar{\bf D}$ for the mapping $F_z$
$$\int_{X\times \bar{\bf D}}\phi_z(x,\t)d\nu_1(x,\t)=\int_{X\times \bar{\bf D}}\phi_z(x,\t)d\nu_2(x,\t).$$
This in turn implies (cf. Lemme, p. 487 \cite{He}) that
$${1\ov n}\sum_{k=0}^{n-1}\phi_z(F_z^k(x,\t))$$
converges uniformly in $(x,\t)\in X\times \bar{\bf D}$ to the common value $C$ of the preceding integrals.

The other items of the theorem are clear.
\end{proof}

We can now identify the imaginary part of the complex rotation number with the geometric rotation number.
\begin{theo}For any $z\in{\bf D}_0$ (and a.e $z\in\pa{\bf D}_0$)
$$\Im\z(z)=\r_g(T,A_z).$$
\end{theo}
\begin{proof}

By Birkhoff's ergodic theorem applied to the $L^\infty(\u)$ function $x\mapsto \Im\Log(\r(z,x,\t(z,x)))$ we get
\begin{equation}{\rm Im}\z(z)=\lim_{n\to\infty}{1\ov n}\sum_{k=0}^{n-1}\Im\Log\r(z,T^kx,\t(z,T^kx)),\lab{air}\end{equation}
that is
\begin{align*}{\rm Im}\z(z)&=\lim_{n\to\infty}{1\ov n}\sum_{k=0}^{n-1}\phi_z(T^kx,\t(z,T^kx))\\
&=\lim_{n\to\infty}{1\ov n}\phi_{z,n}(x,\t_z(x))\end{align*}
and the conclusion of the theorem follows.
\end{proof}

\subsection{Behavior at the boundary}
\subsubsection{Computation of ${\Ree w(\l)}$}
For $|z|<1$,  $m\in{\bf H}_-$   we have
$$C_z\begin{pmatrix}m\\1\end{pmatrix}=(f(z)m+e(z))\begin{pmatrix}{e(z)m-f(z)\ov f(z)m+e(z)}\\1\end{pmatrix},$$
We thus have  that $\Im \hat m<0$ where
$$\hat m={e(z)m-f(z)\ov f(z)m+e(z)}:$$
a simple calculation shows that (we remove the index $z$ to $e(z),f(z)$):
\begin{equation}\Im \hat m=\frac{(|e|^2+|f|^2)\Im m+(1+|m|^2)(1/4)(|z|^2-|z|^{-2})}{|fm+e|^2}\lab{fr}\end{equation}
and therefore $\hat m\in{\bf H}_-$ if $m\in{\bf H}_-$.

Moreover for $x\in{X}$ the matrix 
$$A(x)=\begin{pmatrix}a(x)&b(x)\\c(x)&d(x)\end{pmatrix}$$ is in $SL(2,{\R})$ and it acts also on ${\bf H}_-$ and its action on $\hat m$ is:
$$\begin{pmatrix}a(x)&b(x)\\c(x)&d(x)\end{pmatrix}\bm \hat m\\ 1\em=(c\hat m+d)\begin{pmatrix}\ti m\\1\end{pmatrix},$$
and we have 
$$\Im(\frac{a\hat m+b}{c\hat m+d})=\frac{\Im \hat m}{|c\hat m+d|^2}.$$
Consequentely 
$$A(x)C_z\begin{pmatrix}m\\1\end{pmatrix}=\g\begin{pmatrix}\ti m\\1\end{pmatrix}$$
with
$$\g_-=(f(z)m+e(z))(c(x)\hat m(z)+d(x))$$
and
\begin{align*}\Im \ti m(z,x)&=\frac{\Im \hat m(z)}{|c(x)\hat m(z)+d(x)|^2}\\
&=\frac{(|e(z)|^2+|f(z)|^2)\Im m+(1+|m|^2)(1/4)(|z|^2-|z|^{-2})}{|c(x)\hat m(z)+d(x)|^2|f(z)m+e(z)|^2}\end{align*}

We then get 
\begin{equation}\frac{\Im \ti m(z,x)}{\Im m}=\frac{1}{|\g(z,x)|^2}.\left(|e(z)|^2+|f(z)|^2+\frac{1}{4}(|z|^2-|z|^{-2})\frac{1+|m|^2}{\Im m}\right).\end{equation}

We now set for $|z|<1$
$$z=\frac{\l-i}{i+\l},\qquad \l=i\frac{1+z}{1-z}$$
so that $\l\in{\bf H}_+$.
We get after a simple calulation
$$|e(z)|^2+|f(z)|^2=1+16\frac{|\Im \l|^2}{(1+|\l|^2)^2}(1+O(|\Im\l|^2))$$
and similarly
$${1\ov 4}\left(|z|^2-|z|^{-2}\right)=-2\Im\l\frac{1+|\l|^2}{|1+\l^2|^2}.$$
Since,
\begin{equation}A(x)C_z\begin{pmatrix}m_-(z,x)\\1\end{pmatrix}=\g_-(z,x)\begin{pmatrix}m_-(Tx)\\1\end{pmatrix}\lab{ef}\end{equation}
We thus have
\begin{equation}\frac{\Im \ti m(z,x)}{\Im m}=\frac{1}{|\g(z,x)|^2}\left(1-2\Im\l\frac{1+|\l|^2}{|1+\l^2|^2}\frac{1+|m|^2}{\Im m}+O((\Im \l)^2)\right)\end{equation}

The previous computations show that
\begin{multline}\frac{\Im m_-(z,Tx)}{\Im m_-(z,x)}=\frac{1}{|\g(z,x)|^2}\left(1-2\Im\l\frac{1+|\l|^2}{|1+\l^2|^2}\frac{1+|m_-(z,x)|^2}{\Im m_-(z,x)}+O((\Im \l)^2)\right)\lab{rap}\end{multline}
where the constant in $O(\cdot)$ is universal.
Now we use the following lemma:
\begin{lemme}\lab{fatou}Assume that for any small $\e>0$ we are given $\u$-measurable functions  $f_\e(\cdot):X\to [0,\infty]$   such that for $\u$-a.e $x\in X$ $f_\e(x)$ converges pointwise to $f_0(x)$. Then
$$\liminf_{\e\to 0}{1\ov \e}\int_X\Log(1+\e f_\e(x))d\u(x)\geq\int_X f_0(x)d\u(x).$$
\end{lemme}
\begin{proof}

This is just an application of Fatou's lemma to the sequence of function ${1\ov\e}\Log(1+\e f_\e(x))$.
%From Jensen's inequality ($\u$ is a probability measure) and $\Log(1+t)\leq t$ ($t\geq 0$) we have
%$${1\ov \e}\int_X\Log(1+\e f_\e(x))d\u(x)\leq \int_X f_\e(x))d\u(x).$$
%On the other hand 
%Applying Fatou's lemma:
%$$\int_X f_0(x)d\u(x)\leq \liminf_{\e\to 0}\int{1\ov\e}\Log(1+\e f_\e(x))d\u(x),$$
%and the proof of the lemma is complete.
\end{proof}

It is clear that for each $z\in{\bf D}_0$ the function $-m_-(z,x)$ is positively bounded from below and  thus  taking the $\Log$ in~(\ref{rap}) 
 integrating with respect to $d\u(x)$ and the fact that $\u$ is $T$ invariant ($z\in{\bf D}_0$):
\begin{align*}0&=\int_{X}\Log(-\Im m_-(z,Tx))d\u(x)-\int_{X}\Log(-\Im m_-(z,x))d\u(x)\\
&=-2\Ree w(\l)\\
&+\int_X\Log\left(1-2\Im\l\frac{1+|\l|^2}{|1+\l^2|^2}\frac{1+|m_-(z,x)|^2}{\Im m_-(z,x)}+O((\Im \l)^2)\right)d\u(x)
,\end{align*}  
where 
$$w(\l)=\int_X\Log(\g(\l,x))d\u(x).$$
Now we take the limit $\Im \l\to 0$ and we get using the previous lemma that for a.a $\l_0$ (resp a.a $z_0\in\pa{\bf D}_0$):
$$\liminf_{\Im\l\to 0}(1+(\Ree \l)^2)\frac{\Ree \ti w(\l)}{\Im\l}\geq I_-(\l_0),$$
with
$$I_-(\l)=\int_{X}-\frac{1+|m_-(z,x)|^2}{\Im m_-(z,x)}d\u(x)$$

The same kind relation holds for $m_+$  but we have to justify it (in view of the preceding lemma we have to deal with nonnegative functions). First observe that
$$A_z(T^{-1}x)^{-1}\begin{pmatrix}m_+(x)\\1\end{pmatrix}=\g_+(z,T^{-1}x)^{-1}\begin{pmatrix}m_+(z,T^{-1}x)\\1\end{pmatrix}$$
and setting
$$\u(z,x)\begin{pmatrix}n_+(z,x)\\1\end{pmatrix}=C_z\begin{pmatrix}m_+(z,x)\\1\end{pmatrix},\qquad \ti A(x)=A(T^{-1}x)^{-1}$$ 
one has
$$\ti A(x)C_{(1/z)}\begin{pmatrix}n_+(z,x)\\1\end{pmatrix}=\g_+(z,x)^{-1}\frac{\u(z,T^{-1}x)}{\u(z,x)}\begin{pmatrix}n_+(z,x)\\1\end{pmatrix}.$$
Thus, the same computations that we have already done apply provided we change $z$ in $z^{-1}$: if we define as previously 
$$n=n(x),\qquad \hat n=C_z^{-1}.n,\qquad \ti n=\ti A(x)C_z^{-1}.n$$
we have
\begin{equation}\frac{\Im \ti n}{\Im n}=|\g_+(z,x)|^2.\frac{|\u(z,x)|^2}{|\u(z,Tx)|^2}\left(|e(z)|^2+|f(z)|^2-\frac{1}{4}(|z|^2-|z|^{-2})\frac{1+|n|^2}{\Im n}\right).\end{equation}
But a simple computation shows that
$${1+|n|^2}=\frac{(|e|^2+|f|^2)(1+|m_+(z,x)|^2)+(|z|^2-|z|^{-2})\Im m_+(z,x)}{|fm_+(z,x)+e|^2}$$
and in view of~(\ref{fr}) (with $m$ replaced by $n$):
\begin{align*}\frac{1+|n|^2}{\Im n}&=\begin{pmatrix}(|e|^2+|f|^2)&(|z|^2-|z|^{-2})\\(|z|^2-|z|^{-2})(1/4)&(|e|^2+|f|^2)\end{pmatrix}.\frac{1+|m_+(z,x)|^2}{\Im m_+(z,x)}\\
&=(1+O(\Im\l)\frac{1+|m_+(z,x)|^2}{\Im m_+(z,x)}\end{align*}
Since
$$\int_X \Log(\g_-(z,x))d\u(x)=-\int\Log(\g_+(z,x))d\u(x)$$
(which is a consequence of the fact that matrices $A_z(\cdot)$ are of detreminant 1)
we can now proceed as before and obtain (use the $T$-invariance of $\u$):
$$\liminf_{\Im\l\to 0}(1+(\Ree \l)^2)\frac{\Ree \ti w(\l)}{\Im\l}\geq I_+(\l_0),$$
with
$$I_+(\l)=\int_{X}\frac{1+|m_-(z,x)|^2}{\Im m_-(z,x)}d\u(x)$$
Finally for almost every $z_0\in\pa{\bf D}_0$:
\begin{equation}2\liminf_{\Im\l\to 0}(1+(\Ree \l)^2)\frac{\Ree \ti w(\l)}{\Im\l}\geq I(\l_0)\end{equation}
with 
\begin{equation}\lab{eq:I}I(\l_0)=\int_{X}\left(\frac{1+|m_+(z_0,x)|^2}{\Im m_+(z_0,x)}-\frac{1+|m_-(z_0,x)|^2}{\Im m_-(z_0,x)}\right)d\u(x).\end{equation}

\subsubsection{Computation of $\pa_z w(z)$} Let us now compute $\pa_z w(z)$. If we take the derivative of~(\ref{ef}) w.r.t $z$ we get using the fact that
$$\pa_ze(z)={i\ov z}f(z),\qquad \pa_zf(z)=-{i\ov z}e(z)$$
the equality
\begin{multline}\lab{w1}{i\ov z}A(x)C_z\begin{pmatrix}1\\-m_-(z,x)\end{pmatrix}+A(x)C_z\begin{pmatrix}1\\0\end{pmatrix}\pa_zm_-(z,x)=\\\pa_z\g(z,x)\begin{pmatrix}m_-(z,Tx)\\1\end{pmatrix}+\g(z,x)\begin{pmatrix}1\\0\end{pmatrix}\pa_zm_-(z,Tx).\end{multline}
On the other hand we have the following decomposition in the basis 
$$\begin{pmatrix}m_+(z,x)\\1\end{pmatrix},\qquad\begin{pmatrix}m_-(z,x)\\1\end{pmatrix}$$
\begin{multline*}\begin{pmatrix}1\\-m_-(z,x)\end{pmatrix}=\frac{1+m_-(z,x)m_-(z,x)}{m_+(z,x)-m_-(z,x)}\begin{pmatrix}m_+(z,x)\\1\end{pmatrix}-\\\frac{1+m_-(z,x)m_+(z,x)}{m_+(z,x)-m_-(z,x)}\begin{pmatrix}m_-(z,x)\\1\end{pmatrix}\end{multline*}

and
\begin{multline*}\begin{pmatrix}1\\0\end{pmatrix}=\frac{1}{m_+(z,x)-m_-(z,x)}\begin{pmatrix}m_+(z,x)\\1\end{pmatrix}-\\\frac{1}{m_+(z,x)-m_-(z,x)}\begin{pmatrix}m_-(z,x)\\1\end{pmatrix};\end{multline*}
notice that the last equation is also true with $Tx$ in place of $x$.

Using this and~(\ref{e:54}), (\ref{e:56}),  we can write equation~(\ref{w1}) in the basis 
$$\begin{pmatrix}m_+(z,x)\\1\end{pmatrix},\qquad\begin{pmatrix}m_-(z,x)\\1\end{pmatrix}$$
and we get, projecting on the second vector of this basis and using the fact that the sections $m_\pm(z,x)$ are $(T,A_z)$ invariant, the identity
\begin{multline}{i\ov z}\g_-(z,x)-\frac{1+m_-(z,x)m_+(z,x)}{m_+(z,x)-m_-(z,x)}+\g_-(z,x)\frac{\pa_zm_-(z,x)}{m_+(z,x)-m_-(z,x)}=\\-\pa_z\g_-(z,x)+\g_-(z,x)\frac{\pa_zm_-(z,Tx)}{m_+(z,Tx)-m_-(z,Tx)}\end{multline}

that is,
\begin{multline}{i\ov z}\frac{1+m_-(z,x)m_+(z,x)}{m_+(z,x)-m_-(z,x)}+\frac{\pa_zm_-(z,x)}{m_+(z,x)-m_-(z,x)}=\\-\frac{\pa_z\g_-(z,x)}{\g_-(z,x)}+\frac{\pa_zm_-(z,Tx)}{m_+(z,Tx)-m_-(z,Tx)}.\end{multline}
Integrating w.r.t $x$ and using the $T$-invariance of $\u$ we get
\begin{equation}\pa_z\Log \g_-(z)={i\ov z}\int_{X}\frac{1+m_-(z,x)m_+(z,x)}{m_+(z,x)-m_-(z,x)}d\u(x).\end{equation}

If we put $\ti w(\l)=w(z)$ we get
$$\pa_\l\ti w(\l)=\pa_zw(z).\frac{2 i}{(\l+i)^2};$$
hence
\begin{align}J&={1\ov 2}(\l-i)(\l+i)\pa_\l\ti w(\l)\\
&=\left({1\ov 2}(1+(\Ree \l)^2)+O(\Im \l)\right)\pa_\l \ti w(\l),\end{align}
where
$$J=\int_{X}\frac{1+m_-(z,x)m_+(z,x)}{m_-(z,x)-m_+(z,x)}d\u(x).$$
We hence have (after some computations)
\begin{align*}&\Im \pa_\l\ti w(\l)=\\&2\int_{X}\frac{(1+|m_-(z,x)|^2)\Im m_+(z,x)-(1+|m_+(z,x)|^2)\Im m_-(z,x)}{|m_-(z,x)-m_+(z,x)|^2}d\u,\end{align*}
and we observe that the integrand is nonnegative.
\subsubsection{Consequences}
There is a very nice identity (cf.~\cite{C-J}, \cite{Ko}) which is the key point in our discussion of the complex rotation number:

\begin{multline}I-4\Im J=\int_{X}\left(\frac{1+|m_+(\l,x)|^2}{\Im m_+(\l,x)}-\frac{1+|m_-(\l,x)|^2}{\Im m_-(\l,x)}\right)\\\left(\frac{(\Ree m_-(\l,x)-\Ree m_+(\l,x))^2+(\Im m_-(\l,x)+\Im m_+(\l,x))^2}{|m_-(\l,x)-m_+(\l,x)|^2}\right)d\u\end{multline}

\begin{theo}If there exists an interval $(-\d,\d)$ such that for every $\b\in(-\d,\d)$ one has $\l(T,A(\cdot)R_\b)=0$ then the fibered products of $A(\cdot)$ are bounded and in fact the cocycle $(T,A(\cdot))$ is $C^0$-conjugated to a cocycle with values in $SO(2,{\R})$. Moreover, if $X={\T}$, $T$ is an irrational translation $\th\mapsto \th+\a$ and $A(\cdot)$ homotopic to the identity, $\b$ can be chosen in $(-\d,\d)$ so that the fibered rotation number of $(\a,A(\cdot)R_\b)$ be diophantine w.r.t $\a$. 
\end{theo}
\begin{proof}

We set $K=\{z=e^{i\b},\b\in(-\d,\d)\}\subset\pa{\bf D}_0$ and $\ti K$ the corresponding set of $\l$'s in the real line.
Notice that the condition implies that the map $z\mapsto \l(T,A(\cdot)C_z)$ is continuous on $K$ (the Lyapunov exponent is u.s.c). Since $\l(T,A_z)$ is the real part of $w(\l)$ for $\l\in{\bf H}_{0,+}$ this implies that the limits $\lim_{b\to 0}w(a+ib)$ exists and by a classical version of Schwarz reflection principle, the function $w$ extends holomorphically through the segment (or the arc) $\ti K$. In particular, $w(\cdot)$ admits a  derivative at $\l=\l_0\in{\cal I}$ which implies:
$$\Im\pa_\l\ti w(\l)|_{\l=\l_0}=-\lim_{\Im\l\to 0}\frac{\Ree w(\l)}{\Im\l} $$
which is the identity
$$I-4\Im J\leq 0,$$
and hence $I-4J=0$ since the integrand is nonegative.
The consequences of this equality are discussed in~\cite{C-J} so we will not repeat the argument (which uses Fubini theorem and a generalized Schwarz reflection principle). The result is that for $\u$-a.e $x\in X$ the maps $\l\to m_\pm(\l,x)$ extend holomorhicaly through the line segment $\ti K$ and $m_+(\l,x)=\bar m_-(\l,x)$. Now, since $\u$ is $T$ ergodic, $\u$-a.e $x\in X$ is transitive for $T$ which means that the orbit of $T^nx$ ($n\in\N$) is dense in $X$ and thus we can choose a transitive point $x_0$ in $X$ so that $m_\pm(z,x_0)$ extend holomorphicaly through $K$. For any $x\in X$ there is a sequence $x_n=T^n x_0$ that converges to $x$. An application of Montel's theorem (normal family argument) shows that one can extract a subsequence $x_{n_k}$ so that in a neighborhood of $z_0\in  K$ the sequence $m_\pm(z,x_{n_k})$ converges uniformly to some holomorphic map $\ti m_\pm(z,x)$. Since for $z$ in a compact set of ${\bf D}_0$, $\lim_{n\to\infty}m_\pm(z,x_n)= m_\pm(z,x)$ the map $\ti m_\pm(z,x)$ is a holomorphic extension of $m_\pm(z,x)$ through $K$. Hence for any $x\in X$, $m_\pm(z,x)$ extends to a holomorphic function through $K$. Since $\pm m_\pm(z,x)$ are invariant for $(T,A_z)$ and continuous, we see that for $z\in\pa{\bf D}_0$ $m_+(z,x)-m_-(z,x)$ is never zero at some point  (otherwise this quantity will be identically zero w.r.t $x$ and a look at~(\ref{eq:I}) shows that it is not possible). Therefore for any $z\in K$, the cocycle $(T,A_z)$ can be $C^0$-conjugated to some (possibly non constant) cocycle with values in $SO(2,{\R})$.

The last part of the theorem is a consequence of the fact that the rotation number ${\Ree}\z(z)$ cannot be constant w.r.t $z\in\pa{\bf D}_0$ since otherwise $\z(\cdot)$ would be constant w.r.t $z$ in ${\bf D}_0$, which is impossible with respect to the formula expressing $\pa_\l w(\l)$ in terms of $m_\pm(z,x)$.
\end{proof}

\subsection{The case of the discrete  Schr\"odinger equation}
When $V:X\to{\R}$ is continuous and
$$\ti A(x)=\begin{pmatrix}V(x)&1\\-1&0\end{pmatrix},$$
the complexification we have to use is different and we shall define for $\l\in {\bf H}_+$:
$$\ti A_\l(x)=\begin{pmatrix}V(x)-\l&1\\-1&0\end{pmatrix}.$$
Notice that the plane rotation matrix $Q$ of angle $\pi/2$ conjugates $\ti A_\l(x)$ to $A(x)={\ti A}_\l(x)^{-1}$
$$\begin{pmatrix}0&-1\\1&V(x)-\l\end{pmatrix}=Q\begin{pmatrix}V(x)+\l&1\\-1&0\end{pmatrix}Q^{-1},$$
so that we can work with the matrix $A_\l(x)$ instead of $\ti A_\l(x)$.

For $\Im\l<0$, $A_\l(x)$ preserves ${\bf H}_-$ and $A_\l(x)^{-1}$ preserves ${\bf H}_+$ while contracting the Poincar\'e metric so that we can define similarly invariant sections for $(T,A(\cdot)$ namely $m_\pm(\l,x)$ with values in ${\bf H}_\pm$. With the notations previously defined one finds:
$$\frac{\Im m_-(\l,Tx)}{\Im m_-(\l,x)}=\frac{1}{\g(\l,x)^2}(1-\frac{\Im\l}{\Im m_-(\l,x)})$$
and if we introduce $n_+(\l,x)=-1/m_+(\l,x)$
$$\frac{\Im n_+(\l,T^{-1}x)}{\Im n_+(\l,T^{-1}x)}={\g(\l,x)^2}(1+\Im\l\frac{|n_+(\l,x)|^2}{\Im n_+(\l,x)}).$$
As before one takes the logarithm, integrates w.r.t $\u$, uses lemma~\ref{fatou} and the following equality (in that order)
$$\frac{|n_+|^2}{\Im n_+}=\frac{1}{\Im m_+}$$
to get
$$4\liminf_{\Im\l\to 0^+}\frac{\Ree \ti w(\l)}{\Im\l}\geq I(\l_0),$$
with
$$I(\l)=\int_{X}\frac{1}{\Im m_+(\l,x)}d-\frac{1}{\Im m_-(\l,x)}d\u(x).$$
The computation of $\pa_\l w(\l)$ is done like previously (with the matrix $\ti A_\l(x)$) and we get
$$\pa_\l w(\l)=\int_X\frac{1}{m_+(\l,x)-m_-(\l,x)}d\u(x)=-J$$
Now it is not difficult to see that $I-4J$ is 
\begin{multline}I-4\Im J=\int_{X}\left(\frac{1}{\Im m_+(\l,x)}-\frac{1}{\Im m_-(\l,x)}\right)\\\left(\frac{(\Ree m_-(\l,x)-\Ree m_+(\l,x))^2+(\Im m_-(\l,x)+\Im m_+(\l,x))^2}{|m_-(\l,x)-m_+(\l,x)|^2}\right)d\u\end{multline}
The rest of the proof is then the same as in the previous section.

\pagebreak

\end{document}